\documentclass{article}[12pt]
\usepackage{amsmath,amsfonts,amsthm}
\usepackage{amssymb}
\usepackage{float}
\usepackage{epsfig}
\usepackage{tikz}
\usetikzlibrary{decorations.pathreplacing}
\usepackage{subfigure}
\usepackage{graphics}
\usepackage{graphicx}
\usepackage{color}
\usepackage{xcolor}
\usepackage{soul}
\usepackage[colorlinks=true, citecolor=blue, linkcolor=blue, urlcolor=blue]{hyperref}
\usepackage[makeroom]{cancel}

\newcommand{\parens}[1]{\left(#1\right)}

\newtheorem{corollary}{Corollary}[section]
\newtheorem{definition}{Definition}[section]

\def\tcred{\textcolor{black}}

\usepackage{latexsym}
\usepackage{graphics}
\usepackage{amsmath}
\usepackage{amsthm}
\usepackage{xspace}
\usepackage{amssymb}
\usepackage{color}
\usepackage{cite}
\usepackage{latexsym}
\usepackage{graphics}
\usepackage{amsmath}
\usepackage{amsthm}
\usepackage{xspace}
\usepackage{amssymb}
\usepackage{epsfig}
\usepackage{comment}

\newcommand{\beql}[1]{\begin{equation}\label{#1}}
\newcommand{\eeql}{\end{equation}}
\newcommand{\eqn}[1]{(\ref{#1})}

\newcommand{\R}{\mathbb{R}}

\newcommand{\pr}{\mathbb{P}}
\newcommand{\E}{\mathbb{E}}
\newcommand{\Q}{\mathbb{Q}}

\newcommand{\cs}{{\cal S}}

\newcommand{\ch}{{\cal H}}

\newcommand{\ca}{{\cal A}}

\newcommand{\Z}{\mathbb{Z}}

\def\tcp{\textcolor{black}}
\def\tcgr{\textcolor{black}}
\def\tcr{\textcolor{black}}
\definecolor{deepblue}{RGB}{0,70,140}
\def\tcb{\textcolor{black}}
\def\tcbl{\textcolor{black}}
\def\tcblue{\textcolor{black}}
\definecolor{orange}{rgb}{1,0.45,0}
\def\tcorange{\textcolor{black}}
\newcommand{\edit}[1]{\textcolor{black}{#1}}
\def\tcor{\textcolor{black}}
\def\tr{\textcolor{black}}

\newtheorem{thm}{Theorem}[section]
\newtheorem{lem}[thm]{Lemma}
\newtheorem{prop}[thm]{Proposition}
\newtheorem{cor}[thm]{Corollary}

\theoremstyle{remark}
\newtheorem{remark}[thm]{Remark}
\RequirePackage{graphicx}

\begin{document}

\title{Asymptotic optimality of dynamic first-fit packing on the half-axis
}

\author
{
Philip A. Ernst \\
Department of Mathematics\\
Imperial College London\\
London SW7 2AZ, UK\\
\texttt{p.ernst@imperial.ac.uk}\\
\and
Alexander L. Stolyar \\
ISE Department
and Coordinated Science Lab\\
University of Illinois at Urbana-Champaign\\
Urbana, IL 61801, USA\\
\texttt{stolyar@illinois.edu}
\and
Jixin Wang \\
Department of Mathematics\\
Imperial College London\\
London SW7 2AZ, UK\\
\texttt{jixin.wang23@imperial.ac.uk}\\
}

\date{September 23, 2026}

\maketitle

\begin{center}
\textit{We dedicate this paper to our colleague, mentor, and friend, Professor Larry Shepp (1936-2013)}
\end{center}

\begin{abstract}
 
 We revisit a classical problem in dynamic storage allocation. Items arrive in a linear storage medium, modeled as a half-axis, at a Poisson rate $r$ and depart after an independent exponentially distributed unit mean service time. The arriving item sizes (lengths) are assumed to be independent and identically distributed (i.i.d.) from a common distribution $H$. A widely employed algorithm for allocating the items is the ``first-fit''  discipline, namely, each arriving item is placed in the left-most vacant interval large enough to accommodate it. In a seminal 1985 paper, Coffman, Kadota, and Shepp \cite{CKS86} proved that in the special case of unit length items (i.e. degenerate $H$), 
\edit{ as $r\to\infty$, the first-fit algorithm is asymptotically optimal in the following sense: the steady-state ratio of expected ``empty space'' (gaps between items) to expected occupied space tends towards $0$. In a sequel to \cite{CKS86}, the authors of \cite{CKS85} conjectured that the first-fit discipline is also asymptotically optimal for non-degenerate $H$.}

In this paper we provide the first proof of first-fit asymptotic optimality for non-degenerate distributions $H$ of item sizes. Our main result is for the case 
when $H$ is concentrated on countably many positive real sizes forming an increasing sequence that is either finite or goes to infinity, with the average item size being finite. We prove that under the first-fit discipline, as $r \rightarrow \infty$, the steady-state packing configuration (scaled down by $r$) converges in distribution to the limiting packing configuration with smaller items on the left, larger items on the right, and with no gaps between.
In particular, this proves asymptotic optimality of first-fit in the \tcred{following sense: if $P$ is the expected occupied space, then in steady-state the empty space (scaled down by $r$) in $[0,P]$ vanishes.}

\end{abstract}




\section{Introduction}

We revisit a classical problem in dynamic storage allocation. Suppose that one possesses a linear storage medium, 
such as an optical disk, a parking lot, or a warehouse space, and utilizes it to place ``items'' arriving as a Poisson stream of rate $r$.  \edit{The arriving items
may be of different types, characterized by an item size (or length).} Assume that item
sizes are independent and identically distributed (i.i.d.) from a common distribution $H$. Items are further assumed to depart after an i.i.d. exponential, unit mean ``service'' time (see, e.g., \cite{Aldous,CKS85,CKS86}).
The medium is modeled as a half-axis. 
It is assumed that the system is large scale and that our interest is in the asymptotic regime with $r\to\infty$.
The task of the system administrator is to find a protocol for assigning incoming items. A widely employed algorithm is the ``first-fit'' discipline, which dictates that one should place each arriving item in the left-most vacant interval large enough to accommodate it. 

In a seminal 1985 paper, Coffman et al. \cite{CKS86} proved that for the case of unit length items (i.e. degenerate $H$), the first-fit assignment is asymptotically optimal in the following sense: in the steady-state, the ratio of expected ``empty space'' (gaps between items) to expected occupied space tends to $0$ as $r\to\infty$. In a sequel to \cite{CKS86}, the authors of \cite{CKS85} conjectured the first-fit discipline to be asymptotically optimal for non-degenerate $H$. \edit{Another form of the first-fit asymptotic optimality conjecture is as follows:  under first-fit, as $r \rightarrow \infty$, the steady-state packing configuration (scaled down by $r$) converges in distribution to a limiting packing configuration with no empty space.}
The authors of \cite{CKS85} employed Monte Carlo simulations to empirically confirm the plausibility of their conjecture, but did not offer formal proof.  


In this paper we provide the first proofs of first-fit asymptotic optimality for non-degenerate distributions $H$. We begin with the case when items can be of sizes 1 and 2. We prove that under the first-fit discipline, as $r \rightarrow \infty$, the steady-state packing configuration (scaled down by $r$) converges in distribution to \edit{the limiting packing configuration with smaller items on the left, larger items on the right, and with no gaps between.
Therefore, we prove not only the asymptotic optimality (in the sense of vanishing scaled empty space), but also the specific structure of the limiting 
packing configuration.}
In Sections \ref{section62} and \ref{sec55} we then extend this result to the far more general case in which there are countably many item types with positive real sizes forming an increasing sequence, $\alpha_i, i=1,2,\ldots$, that is either finite or goes to infinity, and with the average item size being finite, i.e. $\sum_{i\geq1}\alpha_i p_i<\infty$, where $p_i$ is the probability of size $\alpha_i$.

We conclude this section with a brief description of related literature. As discussed above, the setting of \cite{CKS85} and \cite{CKS86} is that of a Poisson arrival process with first-fit assignment and an exponential unit mean service time distribution. The authors of \cite{CKS86} also prove that $\E(R) - P = o(r)$, where $R$ is the distance from the origin to the right endpoint of the rightmost item and $P$ is the average occupied space. Aldous \cite{Aldous} considers the same model in the context of the process of parking cars in a parking lot, where each arriving car parks in the lowest-numbered available space. The author provides a lower bound on the asymptotic expected wasted space for first-fit in the case of equally sized items. Coffman and Leighton \cite{CL86} achieve further progress in this direction by providing lower bounds on the wasted space for any online discipline (including first-fit). However, the authors do not prove upper bounds on the wasted space for first-fit with unequally sized items. The present paper closes this gap in the following sense. First, for the model with item sizes 1 and 2 we prove (see Theorem~\ref{conj1}) that the expected empty space in the interval $[0,P]$ is $o(r)$, where in this setting $P=(p_1+2 p_2)r$. We then generalize this result (see Theorem \ref{thm-ff2}) to the case of countably many item types with positive real sizes forming an increasing sequence that is either finite or goes to infinity, with the average item size being finite. As far as the broader literature on dynamic storage allocation is concerned: there is an extensive literature on the efficiency of the first-fit algorithm for more general service time distributions (see \cite{CFL90,CFL91,SK10b,SK11,SK12}). Related references which consider the efficiency of more general allocation disciplines include \cite{CL86,CKLS86,K00,K04,PL07,SK10}.



The remainder of the paper is organized as follows. In Section~\ref{sec-model-result} we formally define the special model and state our main result (Theorem~\ref{conj1}) for it. Section~\ref{sec-proof-main} contains the proof of Theorem~\ref{conj1}, whose basic outline is given and discussed in Section~\ref{sec-proof-outline}. In Section \ref{section62}, we generalize to the case in which there are countably many item types with positive real sizes forming an increasing sequence that is either finite or goes to infinity, with the average item size being finite. The main result for the general model is stated in Theorem \ref{thm-ff2} and its proof is given in Section \ref{sec55}. 


{\bf Basic notation.} The notation $\mathbb{N}_+$, the notation $\Z_+$, the notation $\Q_+$, and the notation $\R_+$ are used for the sets of positive integers, non-negative integers, non-negative rational numbers, and non-negative real numbers, respectively.
The indicator function of a condition or event $A$ is denoted $I(A)$. The notations $\stackrel{P}{\rightarrow}$ 
and $\Rightarrow$ 
denote convergence in probability and in distribution, respectively.

\section{Model: special case with item sizes $1$ and $2$}
\label{sec-model-result}

We begin by formulating the model in which items are of sizes 1 and 2. Customers (or, items) of two types, $i=1,2$, arrive in the system as Poisson processes of rates $p_i r$, where $p_i>0$, $p_1+p_2=1$, and $r$ is a scaling parameter. 
A type-$i$ customer, which we also call an $i$-customer (or $i$-item),
has ``size" $\alpha_i$. Here, $\alpha_1=1$ and $\alpha_2=2$. Arriving customers are placed for ``service" on the half-axis $\R_+$. 
We shall exclusively consider the first-fit discipline: an $i$-customer is placed into the left-most available (not occupied by other customers) interval $[x,x+\alpha_i)$.
Each customer departs after an independent exponentially distributed, unit mean, service time.

For any $r$, the system is stochastically stable, i.e. it has a unique stationary distribution, in which the numbers of $i$-customers are independent Poisson with means $p_i r$.
If we are interested in the steady-state of the process, then, without loss of generality, we can and will consider the process with states such that all customers occupy intervals with integer end points. Then the process describing the system evolution is a continuous-time countable Markov chain which is irreducible and positive recurrent. Specifically, we will define 
a system state as $S= (F_i(\cdot), i=1,2)$, where $F_i(x), x\ge 0$ is the number of $i$-customers located completely to the left of point $x$, i.e. completely within $[0,x)$ (equivalently, completely within $[0,x]$). Note that each $F_i(x)$ is a function of the continuous argument $x\in [0,\infty)$. This does not change the fact that the state space is countable.
The state space $\cs$ of the Markov chain $S(t), ~t \ge 0,$ does not depend on $r$. 

\edit{The process state at time $t$ is $S(t) = (F_i(\cdot;t), i=1,2)$, where $F_i(x;t), x\ge 0,$ is the ``$F_i(x)$ at time $t$.''}
Let $S(\infty) = (F_i(\cdot;\infty), i=1,2)$ be a random element whose distribution is the stationary distribution of 
the process $S(t)$.

Our first main result is given by Theorem \ref{conj1} below.

\begin{thm}
\label{conj1}
As $r\to\infty$, $(1/r)(F_1(p_1 r;\infty), F_2(p_1 r + 2 p_2 r;\infty)) \stackrel{P}{\rightarrow} (p_1,p_2)$.
\end{thm}

\begin{figure}[H]
\centering
        \includegraphics[scale=.4]{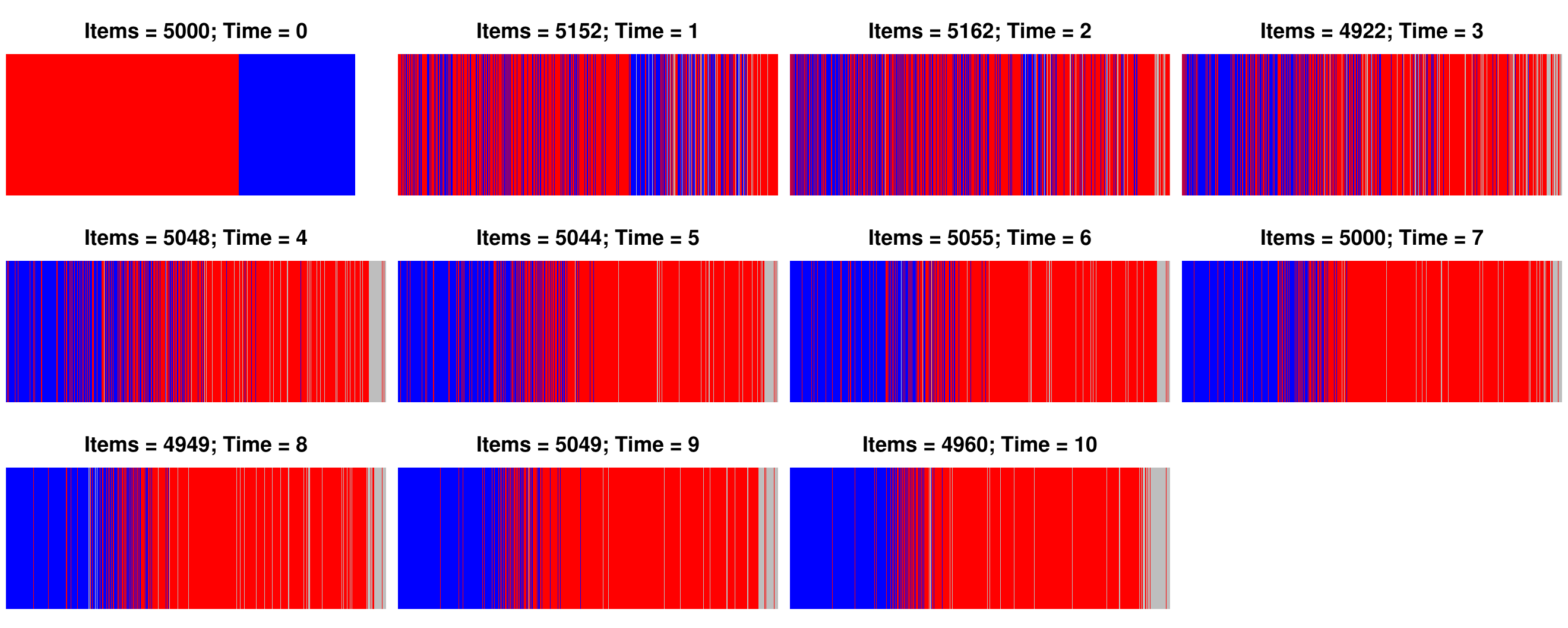}
    \caption{A simulation of the system state evolution for $r=5000$, $p_1=p_2=1/2$. The 1-items appear in blue and the 2-items appear in red. Grey represents empty space.}
    \label{fig-sim}
\end{figure}

The result says that the steady-state of the system with large $r$ is such that ``almost all'' 1-customers concentrate on the left in the interval $[0,p_1 r)$, ``almost all'' 2-customers concentrate immediately to the right
in the interval $[p_1 r, (p_1 + 2 p_2)r)$, 
and that this configuration is ``almost optimal'' in that the empty space in $[0, (p_1 + 2 p_2)r)$ is $o(r)$.
A very high level intuition for the result is that smaller customers (items) naturally have a better chance to take positions further to the left. This intuition is illustrated in Figure~\ref{fig-sim}, in which we show a simulation of the transient behavior of the system with $r=5000, p_1=p_2=1/2$. The system's initial state is ``the opposite'' of its steady-state: we have $p_2 r = 2500$ 2-items on the left and $p_1 r = 2500$ 1-items immediately to the right, with no gaps between items. Figure~\ref{fig-sim} displays system snapshots at regular time intervals. One may observe that the distribution of 1-items fairly quickly ``moves left'' and almost concentrates at the left end of the half-axis, while the distribution of 2-items moves right and almost concentrates immediately to the right of 1-items.

This high level intuition, however, is insufficient to adequately explain why the limiting steady-state configuration is {\em exactly} as specified in Theorem~\ref{conj1}, with ``all'' smaller customers on the left and with ``no'' empty space in 
$[p_1 r, (p_1 + 2 p_2)r)$. In Section~\ref{sec-proof-outline} we provide an outline for the steps of the proof of Theorem~\ref{conj1}. We also offer intuition for each step of the proof.

\section{Proof of Theorem~\ref{conj1}}
\label{sec-proof-main}

\subsection{Outline and basic intuition}
\label{sec-proof-outline}

The proof of Theorem~\ref{conj1} will consist of the following steps, which we state as propositions.

\begin{prop}
\label{prop1}
As $r\to\infty$, $(1/r)(F_1(p_1 r;\infty) + 2 F_2(p_1 r;\infty)) \stackrel{P}{\rightarrow} p_1$. (Equivalently, for any $0<y<p_1$,
$(1/r)(F_1(y r;\infty) + 2 F_2(y r;\infty)) \stackrel{P}{\rightarrow} y$.)
\end{prop}

Proposition~\ref{prop1} states that, as $r\to\infty$, the steady-state total empty space in the interval $[0, p_1 r)$, {\em scaled down by a factor  $1/r$}, vanishes. This is a rather simple property, due to the fact that the new empty space in the interval $[0, y r)$, $y < p_1$, is ``created'' at most at the rate $yr$, while it is -- when positive -- ``eliminated'' at the rate at least $p_1 r$, due to the arrival of 1-items alone.

\begin{prop}
\label{prop2}
As $r\to\infty$, $(1/r) F_2(p_1 r;\infty) \stackrel{P}{\rightarrow} 0$. (Equivalently, for any $0<y<p_1$,
$(1/r) F_2(y r;\infty) \stackrel{P}{\rightarrow} 0$.)
\end{prop}

Proposition~\ref{prop2} states that, as $r\to\infty$, not only does the steady-state scaled total empty space in the interval $[0, p_1 r)$ vanish, but so does the scaled total space occupied by 2-items. This property is more involved. 
The key intuition for why the number of 2-items in $[0,yr)$ in the steady-state cannot be $O(r)$ is as follows. A ``typical'' 1-item departure will only create a size 1 ``hole'' (empty space) which cannot be ``filled in'' by an arriving 2-item, and so will ``typically'' be filled in by a 1-item arrival. On the other hand, a size 2 hole created by a 2-item departure will have a positive probability of receiving a 1-item arrival, precluding the hole being filled in by another 2-item. Thus, if the number of 2-items in $[0,yr)$ were to be $O(r)$, this number would have a negative steady-state drift, which is of course impossible.

\begin{prop}
\label{prop3}
As $r\to\infty$, $(1/r)((F_1(\infty;\infty)-F_1(p_1 r;\infty)), (F_2(\infty;\infty)-F_2(p_1 r;\infty))) \stackrel{P}{\rightarrow} (0,p_2)$.
 \end{prop}
 
Proposition~\ref{prop3} is a corollary of Proposition~\ref{prop1} and Proposition~\ref{prop2}, because we know that the steady-state total number of $i$-items, scaled by $1/r$, converges to $p_i$. Thus, Proposition~\ref{prop3} does not require a separate proof. The proposition states that, in the limit as $r\to\infty$, the scaled numbers of $i$-items are such that ``all'' 1-items are located to the left of the point $p_1 r$ and ``all'' 2-items 
are located to the right of the point $p_1 r$.

{\color{black}
\begin{prop}
\label{prop4}
As $r\to\infty$, $(1/r) F_2(p_1 r + 2 p_2 r;\infty) \stackrel{P}{\rightarrow} p_2$. Equivalently, 
\beql{eq-prop4-cond}
(1/r) F_2(y r;\infty) \stackrel{P}{\rightarrow} (y-p_1)/2, ~\mbox{for any}~ p_1 \le y< p_1+2 p_2.
\eeql
 \end{prop}
 }
 
Proposition~\ref{prop4} completes the proof of Theorem~\ref{conj1}, because it states that the scaled amount of empty space in the interval $[p_1 r, p_1 r +2 p_2 r)$ vanishes. This step of the proof deals with the issue of possible space fragmentation to the right of the interval $[0,p_1 r)$. Indeed, Proposition~\ref{prop3} only tells us that in the steady-state ``all'' 2-items are located to the right of the interval $[0,p_1 r)$, which in principle allows for $O(r)$ unoccupied space  
in $[p_1 r, yr)$. The proof first shows that the unoccupied space in $[p_1 r, yr)$ consists of ``only'' size 1 holes; the proof then shows that if the number of such holes were to be $O(r)$, then this number would have a negative steady-state drift due to $O(r)$-rate of merger of neighboring holes.
 
\subsection{Preliminaries}

For any fixed $r$ the following properties hold. Consider a fixed and finite $x \ge 0$. Denote by $S_x$ the projection of state $S \in \cs$, which only describes the occupancy configuration in the interval $[0,x)$: $S_x = (F_i(\xi), ~\xi \le x, ~i=1,2)$. 
\edit{Denote by $\ch_x$ the set of bounded real-valued functions of $S \in \cs$, which depend only on the projection $S_x$.
(For the special model with two item sizes considered here, the set of possible values of $S_x$ is finite and, therefore, {\em any} function in $\ch_x$ is automatically bounded.)
 Given the structure of our model, it is easy to observe the following: any function $h \in \ch_x$ is within the domain of the (infinitesimal) generator $\ca$ of the Markov process $S(t)$, even though the value $\ca h (S)$ may depend not only on $S_x$.
 (This is because changes of the state projection $S_x$,
which may occur with probabilities $O(\Delta)$ within a small time interval $\Delta>0$, may depend on the occupancy configuration to the right of point $x$.
For example, suppose $x$ is an integer and interval $[x-2,x-1)$ is occupied. Then, $S_x$ does not depend on whether or not 
interval $[x-1,x+1)$ is empty or is occupied by a 2-item. The value of $\ca h (S)$, however, does depend on that, because in the former case it has to account for a possible arrival of a new 1-item into $[x-1,x)$, while in the latter case it does not.)
Indeed, for $h \in \ch_x$, the form of the generator $\ca h$ is determined by possible transitions (and their rates) due to a single item arrival or departure, which change $S_x$ - the occupancy configuration in $[0,x)$. Specifically, for a given $S$ and $x$, let $\Gamma_x(S)$  be the finite (possibly empty) set of items located completely within $[0,x)$, 
let $S^{-\gamma}$ be the state resulting from the item $\gamma$ departure,
and let $S^{+i}$ be the state resulting from the arrival of an $i$-item; note that, unless an $i$-item arrival changes the state projection $S_x$,
i.e. $S^{+i}_x \ne S_x$, $h(S^{+i}) = h(S)$ holds. Note that for each $S$ there is only a finite set of those transitions that may change $S_x$.
Then, for a given $S$ and $h \in \ch_x$,
$$
 \ca h (S)= \sum_{\gamma \in \Gamma_x(S)} [h(S^{-\gamma}) - h(S)] + \sum_{i=1,2} I\{S^{+i}_x \ne S_x\} p_i r [h(S^{+i}) - h(S)].
$$ 
}
We conclude that, for any $x$ and for any $h \in \ch_x$, $\E h (S(\infty))$ is finite, and 
\beql{eq-bar}
\E \ca h (S(\infty))=0,
\eeql
where $S(\infty)$ is the random system state in the stationary regime. We will use this fact repeatedly in what follows.

We next introduce some common notation and conventions employed throughout the proof. We often fix $y>0$ and for each $r$ we consider the projection $S_{\lfloor yr \rfloor}(\infty)$ (or sometimes $S_{\lfloor yr \rfloor +2}(\infty)$) of the stationary system random state $S(\infty)$.
To avoid clogging notation, we will assume, without loss of generality, that $yr$ happens to be an integer. This allows us to write $yr$ instead of $\lfloor yr \rfloor$ in expressions like $[0,\lfloor yr \rfloor)$, $S_{\lfloor yr \rfloor}(\infty)$, etc. and does not cause any problems. To simplify exposition, we will denote by $Y=F_1(yr;\infty)$ and $Z=F_2(yr;\infty)$ the total number of 1-items and 2-items, respectively, (completely) in $[0,yr)$ in the steady-state. Similarly, we denote by $X$ the total occupied space  in $[0,yr)$ in the steady-state (this quantity is ``almost equal'' to $Y + 2 Z$, but may be larger by $1$ due to the possibility of a 2-item occupying $[yr-1, yr+1)$ and not being counted by $Z$. Also note that $X$ is a projection 
of $S_{yr+1}(\infty)$). We denote by $D$ the total number of 2-items that the interval $[0,yr)$ can potentially (completely) fit into its empty space ($D$ is a projection of $S_{yr+2}(\infty)$). A contiguous ``maximal'' empty interval $[k, k+\ell)$,
with integer $k\ge 0$ and $\ell \ge 1$, will be called a size $\ell$ hole, or simply an $\ell$-hole (``maximal'' here means that, for $k\ge 1$, the intervals $[k+\ell,k+\ell+1)$ and $[k-1,k)$  are occupied). Throughout the paper we adopt the convention that when we speak about items or holes {\em in} an interval, we mean {\em completely in.} We denote by $G$ the number of odd-size holes in the interval $[0,yr)$ in the steady-state. We emphasize that the random vector $(Y,Z,X,G,D)$ is a projection
of $S_{yr+2}(\infty)$. We will introduce additional projections of $S_{yr+2}(\infty)$ and other additional notation 
later as necessary.


\subsection{Proof of Proposition~\ref{prop1}}

Let us fix $0< y < p_1$ and, as explained earlier, consider the projection $S_{yr+2}(\infty)$ of the stationary system random state $S(\infty)$. Using the squared empty space $(yr-X)^2$ in $[0,yr)$ as a Lyapunov function, we will use a standard drift argument to show that 
\beql{eq-lyap-drift}
\E (yr-X) \le c < \infty ~~\mbox{for all large $r$}.
\eeql
(This is stronger than we need for Proposition~\ref{prop1}, but we will need \eqn{eq-lyap-drift} later.)

\edit{
We can write
$$
\ca (yr-X)^2 = \sum_{\tau} \mu_\tau [(yr-X+b_\tau)^2 - (yr-X)^2] = \sum_{\tau} \mu_\tau [2(yr-X)b_\tau + b_\tau^2],
$$
where $\tau$ is a transition (changing $X$) due to a single item arrival or departure, $\mu_\tau$ is the rate of this transition, and $b_\tau$ is the increment of 
$yr-X$ due to the transition. The total rate $\sum_{\tau} \mu_\tau$ of all relevant transitions $\tau$ is at most $2r$ (the rate of arrivals is at most $r$ and the rate of departures is at most $yr < r$), and $|b_\tau| \le 2$, so we have
$$
\ca (yr-X)^2 \le \sum_{\tau} \mu_\tau [2(yr-X)b_\tau] + 8r.
$$
Note that if there is any empty space in $[0,yr)$, i.e. $yr-X > 0$, 
we have $\mu_\tau = p_1 r$ for the transition $\tau$, with $b_\tau=-1$, due to a 1-item arrival. We also have $\sum_{\tau: b_\tau = 1} \mu_\tau \le Y+1$ and 
$\sum_{\tau: b_\tau = 2}  \mu_\tau = Z$, for the transitions corresponding to item departures.
So, we obtain
$$
\ca (yr-X)^2 \le 2 (yr-X) [-p_1r + Y  +1 + 2 Z] + 8r \le 2 (yr-X) [-(p_1-y)r/2] + 8 r,
$$
where the second inequality holds for all large $r$. 
Taking expectations on both sides, and recalling \eqn{eq-bar}, we obtain \eqn{eq-lyap-drift} with 
$c=8/(p_1-y)$. $\Box$
}

\subsection{Proof of Proposition~\ref{prop2}}\label{sec34}

Again, we fix $0< y < p_1$ and consider the projection $S_{yr+2}(\infty)$ of $S(\infty)$. 
The general intuition for Proposition~\ref{prop2} is given in Section \ref{sec-proof-outline},
but as the key technical element of this proof, we will consider the drift of the quantity $Z+D$, which counts,
within the interval $[0,yr)$, both the 2-items currently present in the system and the capacity to fit new 2-items.

\tcred{Within this proof it is convenient to consider the quantity $\hat G$, which includes the number $G$ of odd-size holes (completely) within $[0,yr)$,
plus possibly an odd-size empty interval adjacent to an occupied slot \tcorange{or point $0$} on the left and the boundary point $yr$ on the right; clearly, $G \le \hat G \le G+1$.}

Firstly, it follows from Proposition~\ref{prop1} that
\begin{align}\label{eq3}
\lim_{r\rightarrow \infty} \E \hat G/r = 0.
\end{align}
We have 
\begin{equation}\label{eq4}
\begin{aligned}
\ca (Z+D) 
\leq &\,2\hat G  - I(\hat G=0,D>0)\cdot p_1r ~.
\end{aligned}
\end{equation}
This is because $Z+D$ can only increase 
(by at most $1$) due to departure of an item adjacent to an odd-sized hole
\tcred{or the additional empty interval counted into $\hat G$,}
and it decreases, at least, when \tcred{$\hat G=0$} and a 1-item arrives into 
a space potentially available for a 2-item arrival.
Applying expectations to both sides and employing \eqn{eq-bar}, 
we obtain
\begin{align}
\label{eq6}
 \lim_{r\rightarrow \infty}\mathbb{P}(\hat G=0,D>0)=0.   
\end{align}
The remainder of the proof is by contradiction. Suppose, for the sake of contradiction, that Proposition~\ref{prop2} does not hold.
This means that for some $\delta>0$ and $\epsilon>0$, there exists a sequence of $r$ increasing to infinity, along which
\beql{eq-contra}
\pr\{Z \ge r\delta\} \ge \epsilon;
\eeql
for the rest of this proof we consider this sequence of $r$.

\edit{
We see from \eqn{eq-lyap-drift} that for an arbitrarily small $\eta >0$, there exists 
a sufficiently large $K>0$ such that $\pr\{yr - X \le K\} \ge 1-\eta$ for all large $r$.
Let us choose $\eta>0$ sufficiently small (and the corresponding finite $K>0$), so that, 
in view of \eqn{eq-contra}, we have
$$
\pr\{Z \ge r\delta, ~ yr - X \le K \} \ge \epsilon_1 = \epsilon - \eta >0.
$$
Consider the stationary version of the process $S(t)$ over a fixed time interval of order $1/r$ length, say $1/r$ to be specific.
Recall that item arrivals and departures affecting the packing configuration in $[0,yr+2)$ occur at a rate of at most $2r$, 
with the overall arrival rate of $1$-items being $p_1 r$. Then, uniformly in $r$ and all initial states, in the time interval of length $1/r$, with probability 
at least $\delta_1>0$, there will be no departures from the interval $[0,yr+2)$, no $2$-item arrivals, and at least $K$ $1$-item arrivals.
We obtain that, in steady-state, for all large $r$
$$
\pr\{Z \ge r\delta, ~ X=yr \} \ge \epsilon_2 = \epsilon_1 \delta_1 >0.
$$
Consider again the stationary version of the process over 
a $1/r$-length interval. Uniformly in $r$, for any initial state such that $Z \ge r\delta$,
in the time interval of length $1/r$, with probability 
at least $\delta_2>0$, there will be no arrivals into the system, no $1$-item departures from the interval $[0,yr)$, 
no departure of a $2$-item possibly covering point $yr$,
and there will be at least one $2$-item departure from the interval $[0,yr)$.
We obtain that, for all large $r$, in steady-state,
\tcred{$\pr\{\hat{G}=0, D>0\} \ge \epsilon_2 \delta_2 >0$},
which contradicts \eqn{eq6}. This concludes the proof.
$\Box$
}

{\color{black}
\subsection{Proof of Proposition~\ref{prop4}}\label{sec35}

The purpose of this section is to prove Proposition \ref{prop4}, namely \eqn{eq-prop4-cond}.
The set of those $y \in [p_1, p_1 + 2p_2]$ for which property \eqn{eq-prop4-cond} holds has the form $[p_1,y^*]$
(because, clearly, it is closed and we already know that \eqn{eq-prop4-cond} holds for $y=p_1$).
So, it suffices to prove the following ``induction step'': 
{\em if \eqn{eq-prop4-cond} holds for some $y\in [p_1, p_1 + 2p_2)$, then there exists a sufficiently small $\epsilon>0$ (which may depend on $y$),
such that \eqn{eq-prop4-cond} holds for $y+\epsilon$ as well.}

Let us fix $y\in [p_1, p_1 + 2p_2)$ and assume that \eqn{eq-prop4-cond} holds for it. Let us fix a small $\epsilon >0$;
how small, will be specified later.

We need some additional notation. 
Denote by $G_1$ the number of odd-size holes in $[0, p_1 r)$ in steady-state;
we already know that $G_1/r \stackrel{P}{\longrightarrow} 0$ as $r\to\infty$.
By $D^{\epsilon}$ we denote the total number of $2$-items that an empty space in $[0,(y+\epsilon)r)$ can potentially fit.
(Note that $D^0=D$.)
Let 
$G^{\epsilon}$
denote the number of odd-size holes 
in $[p_1 r, (y+\epsilon)r)$. 

Let the random variable $U^{i, \epsilon}$, where $i\in \{1,2, \ldots\} \cup \{\infty\}$, count the number of 
the pairs of odd-size holes in $\left[p_1 r, (y+\epsilon) r \right)$ in the steady-state, which satisfy the following two conditions:
\begin{enumerate}
    \item Between these two odd-size holes there are only 2-items and even-sized holes. 
    \item The right hole is size 1 and the distance between the holes does not exceed $2i$. 
\end{enumerate}

The general direction of the proof of the induction step is as follows. 
Recalling Propositions~\ref{prop1}, \ref{prop2}, \ref{prop3}, and \eqn{eq-prop4-cond}, 
it suffices to show that, as $r\to\infty$, both $D^{\epsilon}/r \stackrel{P}{\longrightarrow} 0$ and $G^{\epsilon}/r \stackrel{P}{\longrightarrow} 0$.
(Because the space in $[p_1 r, (y+\epsilon)r)$ which is {\em not} occupied by $2$-items is upper-bounded by $2D^\epsilon+G^\epsilon+F_1((y+\epsilon)r)-F_1(p_1 r)+3$.)
We will first show (in Lemma~\ref{lem1}) that $D^{\epsilon}/r \stackrel{P}{\longrightarrow} 0$. 
Given that, and $G_1/r \stackrel{P}{\longrightarrow} 0$, in order to prove
$G^\epsilon/r \stackrel{P}{\longrightarrow} 0$ it 
will suffice to prove
 that $U^{\infty, \epsilon}/r \stackrel{P}{\longrightarrow} 0$, which we do in
Lemma \ref{lem2}. 

To improve exposition in the proofs of Lemmas~\ref{lem1} and \ref{lem2} (and later in the paper)
we will use notation $\widetilde o(r)$ for a random variable whose expectation is $o(r)$ as $r\to\infty$.
We also note that the steady-state rate of items' departures (and then arrivals as well), which cover or are adjacent to a given point, is upper-bounded by $2$, and therefore is $\widetilde o(r)$. Using this fact, in the rate estimates we will not explicitly 
account for the ``boundary'' departures/arrivals covering or adjacent to points $p_1 r$, or $yr$, or $(y+\epsilon)r$ -- they will be accounted for
by the $\widetilde o(r)$ term in the corresponding generator expression (or $o(r)$ in the generator expectation).

\medskip

\begin{lem}
\label{lem1}
There exists $\epsilon>0$ (which may depend on $y$), such that 
$\lim_{r\rightarrow \infty} \E D^\epsilon/r = 0$.
\end{lem}
\renewcommand{\proofname}{Proof of Lemma~\ref{lem1}}
\begin{proof}
The intuition is simple. When $D^\epsilon>0$, it has a strong negative part of the drift, $-p_2 r$, due to 2-item arrivals.
Given \eqn{eq-prop4-cond}, interval $[0,p_1 r)$ is ``tightly packed'' (up to $o(r)$ slots) by $1$-items and interval $[p_1 r, yr)$ 
is ``tightly packed'' (up to $o(r)$ slots) by $2$-items. As a result, 
 the positive part of the $D^\epsilon$ drift is, roughly speaking, upper-bounded by 2-item departures from $[p_1 r, yr)$, which ``typically'' increase $D^\epsilon$ by only $1$, 
 plus any item departures from $[yr, (y+\epsilon)r)$; the upper bound is then, roughly, $(y-p_1)r/2 + 2 \epsilon r$, which is smaller than $p_2 r$
 when $\epsilon$ is sufficiently small. We will use a Lyapunov function of the form $[D^\epsilon - m]^+$. The details are as follows. 

 When $D^\epsilon > 0$, its downward drift corresponding to the arrival rate of a size-$2$ item is at least $p_2 r$.  
     
When we consider 
     an increase of $D^\epsilon$, note the following: an item of size-$1$ can increase $D^\epsilon$ (by at most $1$) upon its departure only if it is adjacent to a hole; a departure of a size-$2$ item can increase $D^\epsilon$ by at most $2$; a departure of a size-$2$ item can increase $D^\epsilon$ by more than $1$
     only if the item is adjacent to a hole.
     
     The rate of increase of $D^\epsilon$ can be decomposed 
     into three parts.  The first part corresponds to the increase of available space for size‐$2$ items in the interval $[0,p_1 r)$. 
      The rate of increase of $D^\epsilon$ contributed by this first component is upper-bounded by
$$
\Lambda_1^{\left(D^\epsilon \right)}:=2\biggr(p_1 r - F_1 (p_1 r;\infty)\biggr).
$$
The second component corresponds to the increase in available space for size-$2$ items  in the interval $[ry,r(y+\epsilon))$. The rate of increase of $D^\epsilon$ contributed by this second component can be bounded above by twice the maximum number of items in this interval, which is
$$
\Lambda_2^{\left(D^\epsilon \right)}:=2 \epsilon r.
$$
The third component corresponds to the increase in available space for size-$2$ items in the interval $[p_1 r, yr)$.  
The rate of increase of $D^\epsilon$ contributed by this third component can be upper bounded by
$$
   \Lambda_3^{\left(D^\epsilon\right)}:= \frac{\left(y-p_1\right)r}{2} + 2 [(y - p_1)r - 2 (F_2(yr;\infty) -F_2(p_1 r;\infty))],
  $$    
  where the second term accounts for possible departures of $1$-items and $2$-items adjacent to holes.
  
We denote $\Lambda^{(D^\epsilon)}:=\Lambda_1^{\left(D^\epsilon\right)}+\Lambda_2^{\left(D^\epsilon\right)}+\Lambda_3^{\left(D^\epsilon\right)}$.
Using \eqn{eq-prop4-cond}, we see that 
\beql{eq-Lam-converges}
\Lambda^{(D^\epsilon)}/r \stackrel{P}{\rightarrow} (y-p_1)/2 + 2 \epsilon,
\eeql
and because $\Lambda^{(D^\epsilon)}/r$ is uniformly bounded in $r$, 
\beql{eq-Lam-converges-mean}
\E \Lambda^{(D^\epsilon)}/r \to (y-p_1)/2 + 2 \epsilon.
\eeql
At this point we fix a small $\epsilon>0$, such that
\beql{eq-eps-choice}
(y-p_1)/2 + 2 \epsilon < p_2.
\eeql

For a fixed integer $m>0$, we consider the Lyapunov function $\left(D^\epsilon-m\right)^{+}$. We can write
\begin{equation*}
    \mathcal{A}\left(D^\epsilon-m\right)^{+} \leq-p_2 r\cdot I\left(D^\epsilon>m\right)+\Lambda^{\left(D^\epsilon\right)}\cdot I\left(D^\epsilon\geq m-1\right)
    +\widetilde o(r),
\end{equation*}
where $\widetilde o(r)$ accounts for possible ``boundary'' departures/arrivals.
Since $\mathbb{E}\left[\mathcal{A}\left(D^\epsilon-m\right)^{+}\right]=0$, it follows that
\begin{equation}\label{eq:dm777}
    0 \leq -p_2\cdot \mathbb{P}\left(D^\epsilon>m\right)+\frac{1}{r} \mathbb{E}\left[\Lambda^{\left(D^\epsilon\right)} \cdot I\left(D^\epsilon \geq m-1\right)\right] +o(r)/r.
\end{equation}
Recall \eqn{eq-eps-choice} and
fix any number $z \in ((y-p_1)/2 + 2 \epsilon, p_2)$.
 We then have
\begin{equation}\label{eq:dm2777}
        \frac{1}{r} \mathbb{E}\left[\Lambda^{\left(D^\epsilon\right)} \cdot I\left(D^\epsilon \geq m-1\right)\right]\leq  \frac{1}{r} \mathbb{E}\left[\Lambda^{\left(D^\epsilon\right)} \cdot I\left(\frac{\Lambda^{(D^\epsilon)}}{r} > z \right)\right]+ z \cdot \mathbb{P}\left(D^\epsilon\geq m-1\right) +o(r)/r.
\end{equation}
Combining inequalities \eqref{eq:dm777} and \eqref{eq:dm2777} yields
\begin{equation*}
\begin{aligned}
        (p_2 -z)\cdot \mathbb{P}\left(D^\epsilon>m\right)\leq  z \cdot \mathbb{P}\left(m-1 \leq D^\epsilon\leq m\right)+\frac{1}{r} \mathbb{E}\left[\Lambda^{\left(D^\epsilon\right)} \cdot I\left(\frac{\Lambda^{\left(D^\epsilon\right)}}{r} > z \right)\right]+ o(r)/r.
\end{aligned}
\end{equation*}
For any $N>1$, there exists an integer $m$ with $1 \leq m \leq N$ such that
$$
\mathbb{P}\left(m-1 \leq D^\epsilon \leq m\right) \leq \frac{2}{N} .
$$
Consequently, we have that for any $r$ and $N>1$, 
\begin{equation}\label{dmn777}
    (p_2 - z) \cdot \mathbb{P}\left(D^\epsilon>N\right)\leq z \cdot\frac{2}{N}+\frac{1}{r} \mathbb{E}\left[\Lambda^{\left(D^\epsilon\right)} \cdot I\left(\frac{\Lambda^{(D^\epsilon)}}{r} > z \right)\right]+o(r)/r.
\end{equation}
Hence for any fixed $\eta>0$, 
\begin{equation}\label{dmn2777}
    (p_2 - z) \cdot \mathbb{P}\left(D^\epsilon/r>\eta \right)\leq z \cdot\frac{2}{\eta r}+\frac{1}{r} \mathbb{E}\left[\Lambda^{\left(D^\epsilon\right)} \cdot I\left(\frac{\Lambda^{(D^\epsilon)}}{r} > z \right)\right]+o(r)/r.
\end{equation}
Letting $r\to\infty$, using \eqn{eq-Lam-converges}-\eqn{eq-Lam-converges-mean} and the bounded convergence theorem, we see that the right-hand side of \eqn{dmn2777} vanishes.
Therefore,
$$
 \frac{D^\epsilon}{r}\xrightarrow{P}  0,  \ \text{as} \ r \rightarrow \infty.
$$
Applying bounded convergence once again, we obtain
\begin{equation*}
\lim_{r\rightarrow\infty}\frac{\mathbb{E}\left(D^\epsilon\right)}{r}=0.
\end{equation*}
\end{proof}
\medskip

\begin{lem}
\label{lem2}
For the $\epsilon>0$ specified in Lemma~\ref{lem1},
 $\lim_{r\rightarrow \infty} \E U^{\infty, \epsilon}/r = 0$.
\end{lem}

\renewcommand{\proofname}{Proof of Lemma~\ref{lem2}}
\begin{proof}

We first show that $\lim_r \mathbb{P}(G_1=0)=0$. This is very intuitive. We know that, as $r\to\infty$, the interval $[0,p_1 r)$ is ``completely'' occupied by 1-items. Therefore, 1-items leave the interval $[0,p_1 r)$ at the rate ``equal'' to $p_1 r$. But then 1-items should also enter the interval $[0,p_1 r)$ at the maximal possible rate $p_1 r$. This is only possible if 
$\lim_r \mathbb{P}(G_1 > 0)=1.$ The details are as follows.

Consider $\ca G_1 $. Note that $G_1$ increases 
when a 1-item in interval $[0, p_1 r)$, which is not adjacent to an odd-size hole, 
exits the system, at a rate of at least $F_1(p_1r; \infty) - 2G_1$. The quantity $G_1$ may decrease in two different ways: (i) an item in $[0, p_1 r)$ adjacent to an odd-size hole leaves the system, with the rate not exceeding $2G_1$; and (ii) a new 1-item enters $[0, p_1 r)$ and occupies a place within an odd-sized hole, with rate at most $I(G_1>0)p_1r$. Then, to obtain a bound on $\mathbb{P}(G_1>0)$, we can write
$$
0 = \E \ca G_1 \geq \mathbb{E}[F_1(p_1r;\infty)-2G_1]- 2 \mathbb{E} G_1- \mathbb{P}(G_1>0) p_1 r +o(r),
$$
where $o(r)$ accounts for possible ``boundary'' departures/arrivals.
Invoking Propositions \ref{prop1} and \ref{prop2}, we obtain
$$
\liminf_{r\to\infty}  \mathbb{P}(G_1>0) \geq  \liminf_{r\to\infty} \frac{1}{p_1 r}\biggr\{\mathbb{E}[F_1(p_1r;\infty)-2G_1]- 2 \mathbb{E} G_1 \biggr\}= \frac{1}{p_1}(p_1-0)=1,
$$
and then
\begin{equation}
\label{eqc}
 \lim_{r\to\infty}  \mathbb{P}(G_1=0)=0.   
\end{equation}

We next demonstrate, by induction, that, for any fixed $i \ge 1$,  
\beql{eq-Ui-limit}
\E U^{i,\epsilon}/r \to 0.
\eeql
The intuition for the induction base, i.e. \eqn{eq-Ui-limit} with $i=1$, is that if $U^{1,\epsilon} = O(r)$ were to hold, then the drift of $G^\epsilon$ would be dominated by the negative $O(r)$ drift due to the departures of single 2-items separating pairs of odd holes; this would give $G^\epsilon$ a negative steady-state drift, which is impossible.  The intuition for the induction step from $i\ge 1$ to $i+1$ is that if $U^{i+1,\epsilon} = O(r)$ while $U^{i,\epsilon} = o(r)$, then the drift of 
$U^{i,\epsilon}$ would be dominated by the positive $O(r)$ drift due to 2-item departures adjacent to the left holes of the relevant odd-hole pairs; this would give $U^{i,\epsilon}$ a positive steady-state drift, which is impossible. 
The details are as follows.

As we consider increments of different quantities, say $U^{i,\epsilon}$, let us adopt the following convention. If a system transition, say due to an item departure, simultaneously ``creates'' some hole-pairs counted into $U^{i,\epsilon}$, while ``destroying'' some others, we will count the corresponding increase and decrease of $U^{i,\epsilon}$ separately. For example, if one such hole-pair is created while two others are destroyed, we 
count this as an increase by $1$ and a simultaneous decrease by $2$, with the resulting increment $-1$.

Let $i=1$ and consider $\ca G^\epsilon$.
Note that $G^\epsilon$ can only increase in one of two ways (besides item arrivals overlapping/bordering boundary point $p_1 r$): (i) a 1-item leaves the interval $[p_1 r, (y+\epsilon) r)$, with a rate not exceeding $F_1((y+\epsilon) r;\infty) - F_1(p_1 r;\infty)$; and (ii) a 1-item arrives in the interval $[p_1 r, (y+\epsilon) r)$ and occupies one position in an even-size hole, with a rate upper-bounded by $I(G_1=0)\times p_1r$. There is at least one way for $G^\epsilon$ to decrease: when the sole 2-item between a pair of holes, accounted for by the quantity $U^{1, \epsilon}$, departs, which occurs at the rate equal to $U^{1, \epsilon}$. We can write
$$
0 = \E \ca G^\epsilon \leq \mathbb{P}(G_1=0)\times p_1r+\mathbb{E}[F_1((y+\epsilon) r;\infty)-F_1(p_1 r;\infty)]- \mathbb{E}U^{1, \epsilon} +o(r).
$$
Using (\ref{eqc}) and Proposition \ref{prop3}, we obtain
\begin{equation}\label{equ1}
    \limsup_{r\to\infty} \frac{\mathbb{E}[U^{1, \epsilon}]}{r} \leq  \limsup_{r\to\infty}\mathbb{P}(G_1=0)\times p_1+ \limsup_{r\to\infty} \frac{1}{r}\mathbb{E}[F_1((y+\epsilon) r;\infty)-F_1(p_1 r;\infty)]=0.
\end{equation}

To execute the induction step from $i$ to $i+1$, consider $\ca U^{i,\epsilon}$.
Note that $U^{i,\epsilon}$ can increase in at least one way: when 
there is a pair of odd-size holes, accounted for by the quantity $U^{i+1,\epsilon}$, with the distance between the holes equal to exactly $2(i+1)$,
and the left-most 2-item in the interval separating the two holes departs -- the rate of this is 
$U^{i+1, \epsilon}-U^{i, \epsilon}$. The random variable $U^{i, \epsilon}$ can decrease in one of three ways.
The first is a departure of a 2-item, adjacent 
to a hole belonging to one of the hole-pairs counted into $U^{i, \epsilon}$; the number of such items is at most $4 U^{i, \epsilon}$, and such a departure may destroy at most $3$ hole-pairs counted in $U^{i, \epsilon}$.
The second occurs when a 1-item arrives and occupies a position in the left hole, with an arrival rate upper-bounded by $I(G_1=0)p_1 r $. The third arises when a 1-item next to the two odd-size holes exits the system at a rate upper-bounded by $F_1((y+\epsilon) r; \infty) - F_1(p_1r; \infty)$. 
The above allows us to write
$$
0= \E \ca U^{i,\epsilon} \geq \mathbb{E}\left[U^{i+1, \epsilon}-U^{i, \epsilon}\right]- 3 \cdot 4 \cdot \mathbb{E}U^{i,\epsilon}-\mathbb{E}[F_1((y+\epsilon) r;\infty)-F_1(p_1 r;\infty)]-\mathbb{P}(G_1=0) p_1 r +o(r),
$$
where $o(r)$ accounts for possible ``boundary'' departures/arrivals.
Using (\ref{eqc}) and Proposition \ref{prop3}, we have
$$ 
    \limsup_{r\to \infty}\frac{\mathbb{E}[U^{i+1, \epsilon}]}{r} \leq 13 \cdot \limsup_{r\to \infty} \frac{\mathbb{E} U^{i, \epsilon}}{r}+\limsup_{r\to\infty}\mathbb{P}(G_1=0)p_1+ \limsup_{r\to\infty} \frac{1}{r}\mathbb{E}[F_1((y+\epsilon) r;\infty)-F_1(p_1 r;\infty)] = 0,
$$ 
which completes the induction step and proves \eqn{eq-Ui-limit}.

Note that in the interval $[0, (y+\epsilon) r)$, the number of non-overlapping segments exceeding a length of $2i$ will not surpass $(y+\epsilon) r/(2i)$. This gives $U^{\infty, \epsilon}- U^{i, \epsilon}\leq (y+\epsilon) r/(2i)$. Thus, for any $i \ge 1$,
$$
 \limsup_{r\to \infty}\frac{\mathbb{E}[U^{\infty, \epsilon}]}{r}\leq  \limsup_{r\to \infty}\frac{\mathbb{E}[U^{i, \epsilon}]}{r}+ \frac{y+\epsilon }{2i} = \frac{y+\epsilon }{2i}.
$$
Letting $i\to \infty$ now gives
$$
    \lim_{r\to \infty}\frac{\mathbb{E}[U^{\infty,\epsilon}]}{r}= 0.
$$

\end{proof}

This completes the proof of Proposition \ref{prop4}.
}

Note that the asymptotic regime that we consider is such that, in essence, we are studying the hydrodynamic scaling of the process. In other words, for each $r$, we are interested in the rescaled process
\beql{eq-hd}
f_i^r(x;t) \doteq (1/r) F_i(rx;t), ~~x \ge 0, ~t\ge 0, ~~i=1,2,
\eeql
and, more specifically, in the limit of its stationary distributions as $r\to\infty$. This can naturally lead to an approach to the proof of Theorem~\ref{conj1} based on studying the dynamics of the hydrodynamic limit process $(f_i(\cdot;\cdot))$
obtained as a limit of $(f^r_i(\cdot;\cdot))$ as $r\to\infty$. In fact, this approach does indeed work for 
the proof of Theorem~\ref{conj1}, namely for establishing Propositions~\ref{prop1}-\ref{prop4}. However, in our case, we can choose Lyapunov functions in such a way that it is ``good enough'' to work with instantaneous process drift, described by the process generator. For example, in the proof of Proposition~\ref{prop2}, we use the Lyapunov function $Z+D$, which counts both the present 2-items and the capacity to fit new 2-items, as opposed to the perhaps more ``natural'' Lyapunov function $Z$, which counts only 2-items. We choose to work with the generator of the process, rather than hydrodynamic limits, because this approach (when it works) makes the proofs generally shorter. 

\bigskip

\bigskip

\section{The general model}\label{section62}

The purpose of this section is to generalize the result in Theorem \ref{conj1} to the setting in which we have countably many item types with positive real sizes forming an increasing sequence that is either finite or goes to infinity, with the average item size being finite. The main result appears as Theorem \ref{thm-ff2} and is proven in Section \ref{sec55}.
 
\subsection{A generalization of the model in Section \ref{sec-model-result}}
\noindent We will now consider a generalization of the model in Section \ref{sec-model-result}. We shall again consider a system in which arriving items are placed on the half-line $\mathbb{R}_+$. This time, however, the items may be of countably many types, indexed by $i=1,2,\ldots$. A size-$\alpha_i$ item (an $i$-item, type-$i$ item) has size $\alpha_i$, which is a real positive number. We assume that item sizes form a strictly increasing sequence $0<\alpha_1<\alpha_2<\cdots$, which is either finite or converges to infinity. \tcorange{In the finite-type case, all indices and sums over item types are restricted to the available types.}
Items of size $\alpha_i$ arrive according to independent Poisson processes with rates $p_i r$, where  \tcorange{$p_i>0$,} $\sum_{i=1}^\infty p_i=1$ and $r$ is a scaling parameter. We further assume that the expectation of item size is finite, i.e. $
\sum_{j\geq 1}\alpha_j p_j \;= M \;< \infty.$ As before, we exclusively consider the first-fit discipline: upon arrival, an $i$-item is placed in the leftmost interval $[x, x+\alpha_i)$ that is not occupied by any other item. Each item departs after an independent exponentially distributed, unit mean, service time.

Note that if we are interested in properties of stationary distributions, then, without loss of generality, 
\edit{we can consider the process with a countable state space. Indeed, 
a renewal cycle  (from empty state to empty state), viewed as a sequence of items, with specified types 
and the order in which their arrivals and departures occur,  is finite w.p.1.
The set of possible renewal cycles of a fixed length $m$, and then of all possible finite renewal cycles, is countable.
Therefore, w.p.1., the process gets into and stays within 
the set $\mathcal{S}$ of states that may be reached from the empty state within a finite renewal cycle; this set of states is countable.}

\edit{Starting at this point we will view the process as
a continuous-time Markov chain on a countable state space $\mathcal{S}$. Since this Markov chain is irreducible and positive recurrent (stable), it has a unique stationary distribution. The state at time $t$ \tcorange{will} be denoted by ${S}(t) = (F_i(\cdot;t), i \in \mathbb{N}_{+})$, where $F_i(x;t)$, $x \geq 0$, represents the number of $i$-items located entirely to the left of point $x$, i.e., completely within $[0, x)$.}
We also define ${S}(\infty) = (F_i(\cdot; \infty), i \in \mathbb{N}_{+})$ as a random element whose distribution is the stationary distribution of the process ${S}(t)$.

\indent Some further notation is now required. 
\tcr{Let $F_0(x)$ denote the  total empty space to the left of $x$.
For any $z>0$, denote $S_z =(F_i(x), 0 \le x \le z, i \in \mathbb Z_{+})$ as the projection
of state $S\in\mathcal{S}$, being the packing configuration in $[0,z]$. 
Since $S$ is a continuous-time Markov chain on a countable state space, it is easy to see that if for a finite $z$ a bounded real-valued function $h(S)$ only depends on $S_z$, then $h$ is within the domain of the generator; the fact that the total arrival rate of all items, as well as the total departure rate of the items completely or partially in $[0,z]$, are uniformly bounded above, is used here. Note that the value of $\mathcal{A}h(S)$ may depend not only on $S_z$, but on the packing configuration to the right of point $z$ as well — there is no contradiction here.
Furthermore, for any $h$ in the generator domain we have $\mathbb E[\mathcal{A}h(S(\infty))]=0$. Using this basic identity (which holds for 
a system with any $r$) to derive inequalities shall play an important role in our forthcoming analysis.} For any $i\in \mathbb{N}_{+}$, we denote
$\beta_i=\sum_{j=1}^i \alpha_j p_j$ and set $\beta_0=0$ by convention.
For every $\delta>0$ there exists a constant $A_{\delta}>0$ (depending only on $\delta$)  such that the set of sizes $\alpha_\ell \le A_{\delta}$ is finite and
$$
\sum_{\ell: \alpha_\ell>A_{\delta}} \alpha_\ell p_\ell<\delta.
$$

We now proceed to state Theorem \ref{thm-ff1} and Theorem \ref{thm-ff2}.
\begin{thm}\label{thm-ff1}
     For any $i\in \mathbb{N}_{+}$,  as $r \rightarrow \infty$,
    \begin{equation}
        \frac{1}{r}F_i\left( r \beta_i  ; \infty\right)\xrightarrow{P}p_i.
    \end{equation}
\end{thm}
\noindent  Theorem \ref{thm-ff1} \tcorange{is equivalent to the property that} for  any $y\in[\beta_{i-1},\beta_i)$ and any item type $k$,
  $$
    \frac{1}{r}\Big(F_k(ry;\infty)-F_k(r\beta_{i-1};\infty)\Big)\xrightarrow{P}
    \begin{cases}
      \dfrac{y-\beta_{i-1}}{\alpha_i}, & k=i,\\[6pt]
      0, & k\neq i.
    \end{cases}
$$

Theorem \ref{thm-ff1} states that in the steady state, the storage is almost ordered by item size from left to right: for large $r$, almost all size-$\alpha_i$ items concentrate within $[r\beta_{i-1},\, r\beta_i)$. Moreover, this configuration is “almost optimal” in the sense that the total empty space in $[0,\, rM)$ is $o(r)$. For brevity, we shall refer in the sequel to any system with this property as being ``asymptotically optimally packed."  
\begin{remark}
The statement of Theorem \ref{thm-ff1} and its proof in the paper can also be generalized to the case of different service rates across types, specifically to the case where items of different item types have service times that are exponentially distributed with type-dependent parameters that are uniformly bounded away from zero and infinity.
\end{remark}
\medskip

To prove Theorem \ref{thm-ff1}, we need to establish the result for all $i$. 
A natural proof approach is induction. The following Theorem~\ref{thm-ff2} performs this induction, and therefore implies Theorem \ref{thm-ff1}.

\begin{thm}\label{thm-ff2}
Let $i \in \mathbb{N}_{+}$. Assume that as $r \to \infty$,  $$ \frac{1}{r} F_{j}\left(r\beta_j; \infty\right) \xrightarrow{P} p_j $$ holds for all $j < i$. Furthermore, suppose that for some $y\in [\beta_{i-1},\beta_i)$,
$$
\frac{1}{r} F_i\left(r y ; \infty\right) \xrightarrow{P}\frac{y-\beta_{i-1}}{\alpha_i} .
$$
Then there exists $\epsilon>0$ (which may depend on y and $i$), such that $$\frac{1}{r} F_i\left(r\left(y+\epsilon\right) ; \infty\right) \xrightarrow{P}\frac{y+\epsilon-\beta_{i-1}} {\alpha_i}.$$ \end{thm}

\noindent The reason that Theorem \ref{thm-ff2} implies Theorem \ref{thm-ff1} is that Theorem \ref{thm-ff2} establishes a two–layer induction scheme. For each fixed $i$, define
$$
\mathcal{Y}_i:=\left\{y \in\left[\beta_{i-1}, \beta_i\right]: \frac{1}{r} F_i(r y ; \infty) \xrightarrow{P} \frac{y-\beta_{i-1}}{\alpha_i}\right\}.
$$
We first verify that $\beta_{i-1}\in\mathcal{Y}_i$. \tcblue{ Since
$\beta_j\leq\beta_{i-1}$ for every $j<i$ and no point of the
half-axis is occupied by more than one item,
\begin{equation}\label{eq:Yi-left}
0\leq
\frac{\alpha_i}{r}F_i(r\beta_{i-1};\infty)
\leq
\beta_{i-1}
-\sum_{j=1}^{i-1}\frac{\alpha_j}{r}F_j(r\beta_j;\infty)
\xrightarrow{P}
\beta_{i-1}-\sum_{j=1}^{i-1}\alpha_jp_j
=0.
\end{equation}
Thus $\beta_{i-1}\in\mathcal{Y}_i$}. If $y \in \mathcal{Y}_i \cap\left[\beta_{i-1}, \beta_i\right)$, then by Theorem \ref{thm-ff2} there exists $\epsilon>0$ such that $y+\epsilon \in \mathcal{Y}_i$. \tcb{Moreover, for any $z\in[y,y+\epsilon]$, monotonicity gives
$$
F_i(ry)\le F_i(rz)\le F_i(r(y+\epsilon)).
$$
Since each type-$i$ item has length $\alpha_i$, 
$$
F_i(rz)-F_i(ry)\le \frac{r(z-y)}{\alpha_i}+1,
\quad
F_i(r(y+\epsilon))-F_i(rz)\le \frac{r(y+\epsilon-z)}{\alpha_i}+1.
$$
Given convergence at $y$ and $y+\epsilon$, we have that
$$
\frac1rF_i(rz;\infty)\xrightarrow{P}\frac{z-\beta_{i-1}}{\alpha_i}.
$$
Hence every $z\in[y,y+\epsilon]$ belongs to $\mathcal{Y}_i$.} 
\indent We proceed to prove by contradiction that  $\mathcal{Y}_i=[\beta_{i-1},\beta_i]$. Let $y^*=\sup \mathcal{Y}_i$. \tcblue{Choose $y_n\in\mathcal{Y}_i$ such that $y_n\uparrow y^*$. Since each type-$i$ item has length $\alpha_i$, the fact that at most one type-$i$ item can intersect $ry_n$ gives
$$
\begin{aligned}
&\left|
\frac{1}{r}F_i(ry^*;\infty)
-\frac{y^*-\beta_{i-1}}{\alpha_i}
\right|\\
&\qquad\leq
\left|
\frac{1}{r}F_i(ry_n;\infty)
-\frac{y_n-\beta_{i-1}}{\alpha_i}
\right|
+\frac{y^*-y_n}{\alpha_i}
+\frac{1}{r}.
\end{aligned}
$$
For each fixed $n$, let $r\rightarrow\infty$, and then let $n\rightarrow\infty$. It follows that
$$
\frac{1}{r}F_i(ry^*;\infty)
\xrightarrow{P}
\frac{y^*-\beta_{i-1}}{\alpha_i},
$$
and hence $y^*\in\mathcal{Y}_i$.} Suppose, for the sake of contradiction, that $y^*<\beta_i$. Since $y^*\in\mathcal{Y}_i$, Theorem~\ref{thm-ff2} implies that there exists $\epsilon>0$ such that $y^*+\epsilon\in\mathcal{Y}_i$. This contradicts the definition of $y^*$ as $\sup\mathcal{Y}_i$. Thus, $y^*=\beta_i$, and so $\mathcal{Y}_i=[\beta_{i-1},\beta_i]$.

This completes the inductive step from $i-1$ to $i$. Then induction proceeds by considering type $i+1$ and $y \in [\beta_i,\beta_{i+1})$, and thus, by outer induction over $i=1,2\ldots$, Theorem \ref{thm-ff1} follows. In the sequel, we shall work under the assumptions of Theorem \ref{thm-ff2}. Theorem \ref{thm-ff2} will then be proven in Section \ref{sec55}.
\begin{remark}
    For the sake of brevity: for $y\in [\beta_{i-1},\beta_i)$ we shall use the phrase ``under the inductive hypothesis at $y$" to refer to the following inductive assumptions as $r\to\infty$
$$
\frac{1}{r} F_j\left(r \beta_j ; \infty\right) \xrightarrow{P} p_j \,\, \text { for all } j<i, \ \text{and} \,\,\,\, \frac{1}{r} F_i\left(ry ; \infty\right) \xrightarrow{P}\frac{y-\beta_{i-1}}{\alpha_i}.
$$
\end{remark}

\subsection{\tcgr{Roadmap and some key notation for the proof of Theorem~\ref{thm-ff2}}}\label{sectionroadmap}
\tcgr{In this section, we frontload the necessary notation and provide a roadmap  of the main ideas in the proof of Theorem~\ref{thm-ff2}. For fixed $y\in[\beta_{i-1},\beta_i)$ and $0<\epsilon<(\beta_i-y)/2$, the proof concerns packing in the interval $[r\beta_{i-1},\,r(y+\epsilon))$ under the inductive hypothesis at $y$. The purpose of  the proof is to show that, on the hydrodynamic scale, every mechanism that prevents this interval from being packed by size-$\alpha_i$ items (\tcbl{i.e. increases fragmentation with respect to packing by size-$\alpha_i$ items with no gaps}) contributes only $o(r)$ to the relevant drift estimates as $r\to\infty$.}

\tcgr{For convenience, we shall enumerate (from left to right) all items lying entirely in $[r\beta_{i-1},\,r\beta_i)$ and write them as $[u_1,v_1), [u_2,v_2), \ldots$, so that $u_1<u_2<\cdots$ and $v_j\le u_{j+1}$ for all $j$. Let $e_j$ be the length of item $j$ and let $g_j$ be the length of the hole between items $j$ and $j+1$; that is
\begin{equation}\label{eq:eg-def}
e_j:=v_j-u_j,\quad g_j:=u_{j+1}-v_j.
\end{equation}
}\tcp{For any nonempty collection of items, we shall refer to the ``first item" in the collection as the item with the smallest index in that collection. Similarly, we shall refer to the ``last item" in the collection as the item with the largest index in the collection.} \tcgr{For $y\in[\beta_{i-1},\beta_i)$, we define
\begin{equation}\label{eq:jtilde-def}
j_{\min}(y):=\min\{j:[u_j,v_j)\subset [r\beta_{i-1},\,ry)\},
\quad
j_{\max}(y):=\max\{j:[u_j,v_j)\subset [r\beta_{i-1},\,ry)\}.
\end{equation}
Here, $j_{\min}(y)$ is the index of the first item fully contained in $[r\beta_{i-1},\,ry)$, and $j_{\max}(y)$ is the index of the last item fully contained in $[r\beta_{i-1},\,ry)$. Throughout Section~\ref{sec55}, we shall interpret all sets, item lengths, and hole lengths with respect to this pre-event indexing.}

\tcgr{
The proof of Theorem~\ref{thm-ff2} follows the same overall strategy as the proof of Theorem~\ref{conj1} in Section~\ref{sec-proof-main}, although the Lyapunov functionals used in the general model are more involved. In both settings, we identify packing configurations that prevent the interval under consideration from being asymptotically optimally packed, and we introduce Lyapunov functionals which detect these ``bad" patterns. Under the inductive hypothesis,
every mechanism that can increase the chosen functional contributes only $o(r)$ to its drift. On the other hand, if one of the bad patterns were present in quantity of order $r$, then the first-fit dynamics would force the functional to decrease at a rate of at least $cr$ for some constant $c>0$.  Since the steady-state drift is zero, this implies that such ``bad" patterns can \tcorange{occupy} only \tcorange{$o(r)$ space}. In Sections~\ref{sec34}\textendash\ref{sec35}, this strategy is carried out with simpler functionals, such as $Z+D$ and the odd-hole counts. In Section ~\ref{sec55}, the same strategy is carried out using the Lyapunov functional $\mathrm{TF}$ together with additional counts of certain classes of holes and items.
}

\tcbl{We proceed to \tcbl{informally discuss} a few key objects to be used in the proof. The quantity $H_i(x)$ measures how much fragmentation an interval of length $x$ can contribute, and the quantity $h_i(x)$ is the remainder when $x$ is divided by $\alpha_i$. Their precise definitions are given in \eqref{H_i} and \eqref{def:hi}, respectively, and $H_i(x) \ge h_i(x)$ always holds. The Lyapunov functional $\mathrm{TF}$ measures the total fragmentation in $[r\beta_{i-1},\,ry)$. The other necessary objects are defined later in the proof. 
}

\tcgr{The proof proceeds in four steps. Firstly, Propositions~\ref{prop642}--\ref{prop643} and Corollary~\ref{corr644} \tcb{show that, for item types $k>i$ satisfying $H_i(\alpha_k)>h_i(\alpha_k)$, where $H_i$ is defined in \eqref{H_i} and $h_i$ is defined in \tcorange{\eqref{def:hi}}, the hydrodynamic-scaled numbers of such items in $[r\beta_{i-1},\,r(y+\epsilon))$ \tcorange{converge to zero}.} Secondly, Propositions~\ref{prop645}--\ref{prop646} \tcb{show that, for every fixed $x>0$, the hydrodynamic-scaled number of holes with length at least $x$ in $[r\beta_{i-1},\,r(y+\epsilon))$ \tcorange{converges to zero}. Corollary~\ref{prop647} then shows that, for item types $k>i$ satisfying $H_i(\alpha_k)=h_i(\alpha_k)>0$, the hydrodynamic-scaled numbers of such items also \tcorange{converge to zero}.} Proposition~\ref{prop648} then shows that the total hole length in $[r\beta_{i-1},\,r(y+\epsilon))$ is $o(r)$. Finally, Proposition~\ref{prop649} \tcb{proves that, for item types $k>i$ satisfying $H_i(\alpha_k)=0$, the hydrodynamic-scaled numbers of such items \tcorange{converge to zero}.}}

\tcr{In the sequel, we will use the following convention to avoid unnecessarily clogging expressions. If $B$ is some function of the process state $S$, 
then whether the symbol $B$ means {\em $B$ as a function of a given state} or {\em $B(\infty)$ as its random value in steady-state} will be determined by the context. Specifically, in expressions containing expectation $\E$ or probability $\mathbb P$, $B$ means $B(\infty)$; otherwise, it means a function of a state.
We shall sometimes explicitly write $B(\infty)$ to emphasize the meaning.}

\subsection{Preliminaries}\label{sectionprecoro}
This section contains some key preliminaries needed to prove Theorem \ref{thm-ff2}. Following the notation in Section \ref{sec35}, we let $G_j$ denote the total number of $j$-items that can potentially fit completely into the empty space in the interval $[r\beta_{j-1},\, r\beta_j)$. \tcgr{Further, let $\widetilde{G}_i$ denote the total number of size-$\alpha_i$ items that can potentially fit completely into the empty space in $[0,\,r\beta_{i-1}+\alpha_i)$.}

Proposition \ref{prop31} below shows that, under the induction hypothesis at $\beta_{i-1}$, for every $j<i$, as $r\to\infty$, the steady-state probability that $j$-items can arrive in the interval $[r\beta_{i-1},\, \infty)$ vanishes. \tcgr{ It also shows that, as $r\to\infty$, the steady-state probability that $\widetilde{G}_i>0$ vanishes. Equivalently, the probability that a size-$\alpha_i$ item can fit completely into the empty space in $[0,\,r\beta_{i-1}+\alpha_i)$ tends towards $0$.}
\begin{prop}\label{prop31}
Fix $i \in \mathbb{N}_{+}$ and assume that, as $r \to \infty$,  $$ \frac{1}{r} F_j\left(r\beta_j; \infty\right) \xrightarrow{P} p_j $$ holds for all $j < i$. Then for any $j<i$,
\begin{equation}\label{eqgj}
    \lim _{r\rightarrow\infty} \mathbb{P}\left(G_j=0\right)=0
\end{equation}
\tcgr{
and 
\begin{equation}\label{tildeg}
\lim _{r \rightarrow \infty} \mathbb{P}\left(\widetilde{G}_i>0\right)=0.
\end{equation}
}
\end{prop}
\renewcommand{\proofname}{Proof}
\begin{proof}[Proof of Proposition \ref{prop31}]
We first prove \eqref{eqgj}.
We consider the dynamics of $G_j$, specifically $\mathcal{A} G_j$. Note that $G_j$ increases when a size-$\alpha_j$ item \tcblue{completely contained} in the interval
$
\left[r \beta_{j-1}, r \beta_j\right)
$ exits the system. \tcblue{Among the items counted by
$F_j(r\beta_j)-F_j(r\beta_{j-1})$, at most one intersects
$r\beta_{j-1}$.} Hence, the rate of increase of $G_j$ is at least $
F_j\left(r \beta_j\right)-F_j\left(r \beta_{j-1}\right) \tcblue{-1}.
$

The rate of decrease of $G_j$ has two components. Firstly, $G_j$ may decrease due to the arrival of a size-$\alpha_j$ item in the interval
$
\left[r \beta_{j-1}, r \beta_j\right),
$ which occurs with rate at most 
$$p_j r \cdot I\left(G_j>0\right).$$ 
Secondly, a size-$\alpha_\ell$ item with $\ell\neq j$ may arrive in the same interval. Each such arrival can reduce $G_j$ by at most $\bigl\lfloor \alpha_\ell/\alpha_j \bigr\rfloor + 1$, since a larger item may occupy extra space that could otherwise accommodate  size-$\alpha_j$ items. In the steady state, the expectation of the rate at which size-$\alpha_\ell$ items are placed \tcblue{with their right endpoints in
$(r\beta_{j-1},\,r\beta_j]$} equals the expectation of the rate \tcblue{at which such items depart}, namely
$$
\mathbb{E}\left(F_\ell\!\left(r\beta_j;\,\infty\right) - F_\ell\!\left(r\beta_{j-1};\,\infty\right)\right).
$$
Therefore, for any $\ell \neq j$, the expectation of the rate of decrease of $G_j$ due to the arrival of \tcblue{such} size-$\alpha_\ell$ items can be upper bounded by 
$$
\left(\left\lfloor\frac{\alpha_{\ell}}{\alpha_j}\right\rfloor+1\right)\mathbb{E}\left(F_{\ell}\left(r \beta_j ; \infty\right)-F_{\ell}\left(r \beta_{j-1} ; \infty\right)\right) .
$$
\tcblue{A size-$\alpha_\ell$ arrival not counted by
$F_\ell(r\beta_j)-F_\ell(r\beta_{j-1})$ can reduce $G_j$ only if
the arriving item intersects $r\beta_j$. Since at most one item
can intersect this point and each item departs at rate one, in the
steady state the expectation of the aggregate rate at which such
items are placed is at most one. Hence, for items of size at most
$A_\delta$, the expectation of the rate of decrease of $G_j$ due
to such arrivals can be upper bounded by}
$$
\tcblue{
\left\lfloor\frac{A_\delta}{\alpha_j}\right\rfloor+1.
}
$$
\tcblue{The contribution from items of size greater than
$A_\delta$ is included in the last term in \eqref{eq:ag}.}
Therefore,
\begin{equation}\label{eq:ag}
    \begin{aligned}
0=\mathbb{E}\left(\mathcal{A} G_j\right)\geq & \mathbb{E}\left(F_j\left(r \beta_j ; \infty\right)-F_j\left(r \beta_{j-1} ; \infty\right)\right)\tcblue{-1}-p_j r \cdot \mathbb{P}\left(G_j>0\right)  \\
&- \sum_{\ell: \alpha_j\neq \alpha_\ell\leq A_{\delta}} \left(\left\lfloor\frac{ A_{\delta}}{\alpha_j}\right\rfloor+1\right)\mathbb{E}\left(F_{\ell}\left(r \beta_j ; \infty\right)-F_{\ell}\left(r \beta_{j-1} ; \infty\right)\right)\\
&\tcblue{-\left(\left\lfloor\frac{A_\delta}{\alpha_j}\right\rfloor+1\right)}\\
&- \mathbb{E} \left(\sum_{\ell: \alpha_\ell> A_{\delta}} \left(\left\lfloor\frac{ \alpha_{\ell}}{\alpha_j}\right\rfloor+1\right)F_{\ell}\left(\infty ; \infty\right)\right).
\end{aligned}
\end{equation}
\tcblue{For any $j<i$, the items counted by $F_j(r\beta_{j-1};\infty)$ and by $F_q(r\beta_q;\infty)$, $1\leq q<j$, all lie in $[0,r\beta_{j-1})$, so their total length is at most $r\beta_{j-1}$. Since $0\leq \alpha_qF_q(r\beta_q;\infty)/r\leq\beta_q$ for every $q\leq j$, our inductive hypothesis gives
$$
\begin{aligned}
    0
&\leq\limsup_r\frac{1}{r}\mathbb{E}\left(F_j(r\beta_{j-1};\infty)\right)\\
&\leq\limsup_r\frac{1}{\alpha_j}\left(\beta_{j-1}-\sum_{q=1}^{j-1}\frac{\alpha_q}{r}\mathbb{E}\left(F_q(r\beta_q;\infty)\right)\right)=0,\\
\end{aligned}
$$
hence,
\begin{equation}\label{eq:lj}
\begin{aligned}
\lim_r\frac{1}{p_jr}\mathbb{E}\left(F_j(r\beta_j;\infty)-F_j(r\beta_{j-1};\infty)\right)&=1.
\end{aligned}
\end{equation}
For $\ell<j$, since $\beta_\ell\leq\beta_{j-1}$, our inductive hypothesis and the fact
$\mathbb{E}\left(F_\ell(\infty;\infty)\right)=p_\ell r$ give
\begin{equation}\label{eq:lj1}
\begin{aligned}
0
&\leq\limsup_r
\frac{1}{r}\mathbb{E}\left(
F_\ell(r\beta_j;\infty)-F_\ell(r\beta_{j-1};\infty)
\right)\\
&\leq\limsup_r
\frac{1}{r}\mathbb{E}\left(
F_\ell(\infty;\infty)-F_\ell(r\beta_\ell;\infty)
\right)
= p_\ell-p_\ell=0.
\end{aligned}
\end{equation}
For $\ell>j$, since $\beta_q\leq\beta_j$ for every $q\leq j$, our inductive hypothesis gives
\begin{equation}\label{eq:lj2}
\begin{aligned}
0
&\leq
 \limsup_r \frac{1}{r}\mathbb{E}\left(
F_\ell(r\beta_j;\infty)-F_\ell(r\beta_{j-1};\infty)
\right)\\
&\leq \limsup_r
\frac{1}{r}\mathbb{E}\left(F_\ell(r\beta_j;\infty)\right) \leq\limsup_r
\frac{1}{\alpha_\ell}
\left(
\beta_j-
\sum_{q=1}^{j}
\frac{\alpha_q}{r}
\mathbb{E}\left(F_q(r\beta_q;\infty)\right)
\right)
= 0.
\end{aligned}
\end{equation}}
\tcblue{For fixed $\delta$, the two constant terms in \eqref{eq:ag}, divided by $p_jr$, vanish as $r\rightarrow\infty$.}
Combining \tcblue{\eqref{eq:ag}--\eqref{eq:lj2} } yields
$$
\begin{aligned}
        \liminf_{r}\mathbb{P}\left(G_j>0\right) \geq & \lim_{r}\frac{1}{p_jr}\mathbb{E}\left( F_j\left(r \beta_j ; \infty\right)-F_j\left(r \beta_{j-1} ; \infty\right)\right)\\
        &- \limsup_{r}\sum_{\ell: \alpha_j\neq\alpha_\ell\leq A_{\delta}}  \frac{1}{p_jr} \left(\left\lfloor\frac{ A_{\delta}}{\alpha_j}\right\rfloor+1\right)\mathbb{E}\left(F_{\ell}\left(r \beta_j ; \infty\right)-F_{\ell}\left(r \beta_{j-1} ; \infty\right)\right)\\
        &- \limsup_{r}\frac{2}{\alpha_1p_jr}\mathbb{E}\left(\sum_{\ell: \alpha_\ell> A_{\delta}}\alpha_{\ell}F_{\ell}\left(\infty ; \infty\right)\right)\geq 1-\frac{\tcgr{2}\delta}{\alpha_1p_j}.
\end{aligned}
$$
Thus, taking $\delta\rightarrow 0$, we obtain $\lim _r \mathbb{P}\left(G_j=0\right)=1-\lim _r \mathbb{P}\left(G_j>0\right)=0$.   

\tcgr{
We next prove (\ref{tildeg}). \tcblue{Taking expectations in the inequalities in \eqref{eq:Yi-left} gives}
$$
\lim_{r\to\infty}\frac{1}{r}\mathbb{E}\left(F_i\left(r\beta_{i-1};\infty\right)\right)=0.
$$
Since an interval of length $\alpha_i$ can contain at most one size-$\alpha_i$ item,
$$
0\le F_i\left(r\beta_{i-1}+\alpha_i\right)-F_i\left(r\beta_{i-1}\right)\le 1.
$$
Therefore,
$$
\lim_{r\to\infty}\frac{1}{r}\mathbb{E}\left(F_i\left(r\beta_{i-1}+\alpha_i;\infty\right)\right)=0.
$$
We now consider the dynamics of $F_i\left(r\beta_{i-1}+\alpha_i\right)$. Whenever $\widetilde{G}_i>0$, there exists an interval of length $\alpha_i$ that is completely contained in the empty space in $[0,r\beta_{i-1}+\alpha_i)$.  Hence, by the first-fit allocation rule, whenever $\widetilde{G}_i>0$, an arriving size-$\alpha_i$ item increases $F_i\left(r\beta_{i-1}+\alpha_i\right)$ by $1$. Thus, the rate of increase of $F_i\left(r\beta_{i-1}\tcblue{+\alpha_i}\right)$ is at least
$$
p_ir \cdot I\left(\widetilde{G}_i>0\right).
$$
On the other hand, the rate of decrease of $F_i\left(r\beta_{i-1}+\alpha_i;\infty\right)$ is exactly
$
F_i\left(r\beta_{i-1}+\alpha_i\right).
$
Therefore,
$$
0=\mathbb{E}\left(\mathcal{A}F_i\left(r\beta_{i-1}+\alpha_i;\infty\right)\right)
\geq p_ir\cdot \mathbb{P}\left(\widetilde{G}_i>0\right)-\mathbb{E}\left(F_i\left(r\beta_{i-1}+\alpha_i;\infty\right)\right).
$$
It then follows that
$$
\limsup_r\mathbb{P}\left(\widetilde{G}_i>0\right)
\leq \limsup_r
\frac{1}{p_ir}\mathbb{E}\left(F_i\left(r\beta_{i-1}+\alpha_i;\infty\right)\right)=0.
$$
}
\end{proof}
\noindent Note that Proposition \ref{prop31} implies the following: under the inductive hypothesis at $\beta_{i-1}$, and in the steady state, as $r\to\infty$, for every $j<i$, the probability that there exists at least one hole capable of accommodating a size-$\alpha_j$ item in the interval $[r\beta_{j-1},\, r\beta_j)$ tends to $1$. Consequently, any arriving size-$\alpha_j$ item lands \tcgr{completely}  in $[0,\, r\beta_{i-1})$ with probability $1-o(1)$ as $r\to\infty$. Thus, the rate at which size-$\alpha_j$ items arrive in the interval $[r\beta_{i-1},\,\infty)$ is $o(r)$.

\tcgr{Likewise, by the proof of (\ref{tildeg}), the aggregate rate at which size-$\alpha_i$ items are placed either entirely in $[0,\,r\beta_{i-1})$ or in intervals $[x,\,x+\alpha_i)$ where $$x<r\beta_{i-1}<x+\alpha_i,$$ is $o(r)$. \tcb{In particular, the rate at which size-$\alpha_i$ items with left endpoints in $[r\beta_{i-1},\,\infty)$ are placed is $p_ir-o(r)$.}}

\noindent For any $y$ and $\epsilon$, we denote by $D_i$ the total number of $i$-items that the interval $\left[0, r(y+\epsilon)\right)$ can potentially (completely) fit into its empty space. Then we have the following proposition\tcorange{, generalizing Lemma~\ref{lem1}}.
\begin{prop}\label{npr1}
   For any $y\in[\beta_{i-1},\beta_i)$, there exists  $\epsilon>0$ such that, under the inductive hypothesis at $y$,
   $$
       \lim_{r\rightarrow\infty}\frac{\mathbb{E}(D_i)}{r}=0.
   $$
 Moreover, for any \tcgr{$\tcorange{0<}\epsilon<(\beta_i-y)/2$}, $N\tcorange{\in\mathbb N_+}$ and $\delta>0$, \tcp{ letting}
 \tcgr{
 \begin{equation}\label{eq:B-di}
B_{\delta,i}:=\left\lceil\frac{A_\delta}{\alpha_i}\right\rceil\tcblue{+1},
\end{equation}
and
\begin{equation}\label{eq:z-ye}
z_{y,\epsilon}:=\frac{2\epsilon+y-\beta_{i-1}}{\alpha_i}<p_i,
\end{equation}}
     $$
      \limsup_{r\rightarrow\infty}\mathbb{P}(D_i>N)\leq  \frac{p_i+z_{y,\epsilon}}{p_i-z_{y,\epsilon}}\frac{\tcp{B_{\delta,i}+1}}{N}+\frac{6\delta}{\alpha_1(p_i-z_{y,\epsilon})} .
     $$

\end{prop}

\begin{proof}
     When $D_i>0$, its downward drift due to size-$\alpha_i$ arrivals is at least $p_ir$.
\tcblue{Suppose that an item of size $\alpha_\ell$ and its two adjacent holes are completely contained in $[0,r(y+\epsilon))$, and let \tcorange{$s,t\geq0$} be the lengths of these holes\tcorange{, taking the corresponding length to be zero when there is no adjacent vacant interval on that side}. Let $\Delta D_i$ denote the change in $D_i$ caused by the departure of this item. Then we have}
\begin{equation}\label{eq:Di-jump}
\tcblue{\begin{aligned}
\Delta D_i
&=
\left\lfloor\frac{s+\alpha_\ell+t}{\alpha_i}\right\rfloor
-
\left\lfloor\frac{s}{\alpha_i}\right\rfloor
-
\left\lfloor\frac{t}{\alpha_i}\right\rfloor\\
&\leq
\frac{
\alpha_\ell
+s-\alpha_i\left\lfloor s/\alpha_i\right\rfloor
+t-\alpha_i\left\lfloor t/\alpha_i\right\rfloor
}{\alpha_i}
.
\end{aligned}}
\end{equation}
\tcblue{Since the two remainder terms are each less than $\alpha_i$ and $\Delta D_i$ is an integer, \eqref{eq:Di-jump} also gives $\Delta D_i \leq \left\lceil\alpha_{\ell} / \alpha_i\right\rceil+1$. When the upper bound in \eqref{eq:Di-jump} is summed over departures, $\alpha_\ell$ is counted once for each departing item, while each hole can appear as $s$ or $t$ for at most two departures.}
The rate of increase of $D_i$ \tcblue{caused by departures for which the departing item and its two adjacent holes are completely contained in $[0,r\beta_{i-1})$, $[r\beta_{i-1},ry)$, or $[ry,r(y+\epsilon))$} can be decomposed into three parts.  The first part corresponds to the increase of available space for size‐$\alpha_i$ items in the interval $[0,r\beta_{i-1})$. \tcblue{If $i=1$, this interval is empty, and we set $\Lambda_1^{(D_1)}:=0$. Suppose that $i\geq2$. The quantity
$$
r\beta_{i-1}-\sum_{\ell=1}^{i-1}\alpha_\ell F_\ell(r\beta_{i-1})
$$
is an upper bound for the sum of the total hole length and the total length occupied by items of size at least $\alpha_i$ in $[0,r\beta_{i-1})$. Summing the upper bound in \eqref{eq:Di-jump} over departures of items of size at least $\alpha_i$ in the first part shows that their contribution to the rate of increase of $D_i$ is bounded above by}
$$
\tcblue{\frac{2}{\alpha_i}
\left(
r\beta_{i-1}-\sum_{\ell=1}^{i-1}\alpha_\ell F_\ell(r\beta_{i-1})
\right).}
$$
\tcblue{For departures of items of size less than $\alpha_i$ in the first part, we have $\alpha_\ell\leq\alpha_{i-1}$. Whenever $\Delta D_i>0$, the upper bound in \eqref{eq:Di-jump} gives}
\begin{equation}\label{eq:Di-small-jump}
\tcblue{
\Delta D_i
\leq
\frac{
 s-\alpha_i\left\lfloor s/\alpha_i\right\rfloor
+t-\alpha_i\left\lfloor t/\alpha_i\right\rfloor
}{\alpha_i-\alpha_{i-1}}\leq  \frac{
 s
+t
}{\alpha_i-\alpha_{i-1}}.
}
\end{equation}
\tcblue{Since each hole is adjacent to at most two items, summing \eqref{eq:Di-small-jump} over departures of items of size less than $\alpha_i$ in the first part shows that their contribution to the rate of increase of $D_i$ is bounded above by}
$$
\tcblue{\frac{2}{\alpha_i-\alpha_{i-1}}
\left(
r\beta_{i-1}-\sum_{\ell=1}^{i-1}\alpha_\ell F_\ell(r\beta_{i-1})
\right).}
$$
\tcblue{Thus, the }\tcp{rate of increase of }$D_i$ \tcp{contributed by this first }\tcblue{part is bounded above by}
$$
\Lambda_1^{\left(D_i\right)}:=\parens{\frac{2}{\alpha_i}+\frac{\tcblue{2}}{\alpha_i-\alpha_{i-1}}}\biggr(r\beta_{i-1} - \sum_{\ell=1}^{i-1} \alpha_{\ell}F_{\ell}(r\beta_{i-1})\biggr).
$$
The second \tcblue{part} corresponds to the increase in available space for size-$\alpha_i$ items  in the interval $[ry,r(y+\epsilon))$. \tcblue{Summing the upper bound in \eqref{eq:Di-jump} over departures in the second part shows that their contribution is bounded above by the total occupied length plus twice the total hole length in $[ry,r(y+\epsilon))$, divided by $\alpha_i$. Since the sum of the total occupied length and the total hole length in this interval is at most $r\epsilon$, this contribution is at most $2r\epsilon/\alpha_i$, and hence the }\tcp{rate of increase of }$D_i$ \tcp{contributed by this second }\tcblue{part} \tcp{ can be bounded above by}
$$
\Lambda_2^{\left(D_i\right)}:=2\biggr(\frac{r\epsilon}{\alpha_i}+1\biggr) .
$$
The third \tcblue{part} corresponds to the increase in available space for size-$\alpha_i$ items in the interval $[r\beta_{i-1}, ry)$. \tcblue{Since at most one size-$\alpha_i$ item counted by $F_i(r y)-F_i\left(r \beta_{i-1}\right)$ intersects $r \beta_{i-1}$, the total length occupied by size-$\alpha_i$ items completely contained in $\left[r \beta_{i-1}, r y\right)$ is at least $\alpha_i\left(F_i(r y)-F_i\left(r \beta_{i-1}\right)-1\right)$.}  \tcblue{Hence, the sum of the total hole length and the total length occupied by items whose sizes differ from $\alpha_i$ is bounded above by}
$$
\tcblue{r(y-\beta_{i-1})
-\alpha_i\bigl(F_i(ry)-F_i(r\beta_{i-1})\bigr)
+\alpha_i.}
$$
\tcblue{For a departing item of size $\alpha_i$, the upper bound in \eqref{eq:Di-jump} gives}
$$
\tcblue{\Delta D_i
\leq
1+
\frac{
s-\alpha_i\left\lfloor s/\alpha_i\right\rfloor
+t-\alpha_i\left\lfloor t/\alpha_i\right\rfloor
}{\alpha_i}\leq 1+
\frac{
s
+t
}{\alpha_i}.}
$$
\tcblue{After summing over departures of size-$\alpha_i$ items in the third part,  the constant $1$ in this bound contributes at most $r(y-\beta_{i-1})/\alpha_i$, and the terms involving $s$ and $t$ contribute at most twice the total hole length in $[r\beta_{i-1},ry)$, divided by $\alpha_i$. For departures of items whose sizes differ from $\alpha_i$ in the third part, the upper bound in \eqref{eq:Di-jump} gives the total length occupied by these items plus twice the total hole length in this interval, divided by $\alpha_i$. Hence, the }\tcp{rate of increase of }$D_i$ \tcp{contributed by this third }\tcblue{part} \tcp{ can be upper bounded by}
$$
   \Lambda_3^{\left(D_i\right)}:= \frac{\tcblue{4}}{\alpha_i}\biggr( r \left(y-\beta_{i-1}\right)-\alpha_i \big(F_i(ry) - F_i(r\beta_{i-1})\big)\tcblue{+\alpha_i}\biggr)+\frac{r \left(y-\beta_{i-1}\right)}{\alpha_i}.
$$    
\tcblue{It remains to consider departures not covered by the three upper bounds $\Lambda_1^{(D_i)}$, $\Lambda_2^{(D_i)}$, and $\Lambda_3^{(D_i)}$. Such a departure can occur only if the departing item or one of its two adjacent holes intersects $r\beta_{i-1}$, $ry$, or $r(y+\epsilon)$. For each of these three points, if the point is occupied by an item, only the departure of that item is not covered by the three upper bounds. If the point lies in a hole, only departures of items adjacent to that hole are not covered, and there are at most two such items. Thus, at most six departures need to be considered. If the departing item or one of its adjacent holes intersects $r(y+\epsilon)$, only the portions lying in $[0,r(y+\epsilon))$ can increase $D_i$. Applying \tcorange{the estimate analogous to} \eqref{eq:Di-jump} to these portions gives $\Delta D_i\leq\left\lceil\alpha_\ell/\alpha_i\right\rceil+1$. Therefore, since each of these at most six items departs at rate
one, the rate of increase of $D_i$ caused by departures of those
items having size at most $A_\delta$ is bounded above by
$6B_{\delta,i}$. The contribution from departures of items of size greater than $A_\delta$ is included in the term $\sum_{\ell:\alpha_\ell>A_\delta}\left(2+\frac{\alpha_\ell}{\alpha_i}\right)F_\ell(\infty)$.}

We denote $\Lambda^{(D_i)}:=\Lambda_1^{\left(D_i\right)}+\Lambda_2^{\left(D_i\right)}+\Lambda_3^{\left(D_i\right)}\tcblue{+6B_{\delta,i}}$ and for any \tcblue{$m\in\mathbb{N}_+$}, we take the Lyapunov function $\left(D_i-m\right)^{+}$. \tcgr{Recall the constant }$B_{\tcblue{\delta},i}$ \tcgr{ from \eqref{eq:B-di}}\tcblue{. Since \eqref{eq:Di-jump} gives $\Delta D_i\leq\left\lceil\alpha_\ell/\alpha_i\right\rceil+1$, }when $D_i<\tcp{m-B_{\delta,i}}$, increases of $\left(D_i-m\right)^{+}$ can only be caused by departures of items with size exceeding $A_{\delta}$. Hence the rate of increase when $D_i<\tcp{m-B_{\delta,i}}$ is upper bounded by $\sum_{\ell: \alpha_{\ell}>A_{\delta}}\left(2+\frac{\alpha_{\ell}}{\alpha_i}\right) F_{\ell}(\infty)$. Therefore,
\begin{equation*}
    \mathcal{A}\left(D_i-m\right)^{+} \leq-p_i r\cdot I\left(D_i>m\right)+\Lambda^{\left(D_i\right)}\cdot I\left(D_i\geq \tcp{m-B_{\delta,i}}\right)+\sum_{\ell: \alpha_{\ell}>A_{\delta}}\left(2+\frac{\alpha_{\ell}}{\alpha_i}\right) F_{\ell}(\infty).
\end{equation*}
Since $\mathbb{E}\left[\mathcal{A}\left(D_i-m\right)^{+}\right]=0$, it follows that
\begin{equation}\label{eq:dm}
    0 \leq -p_i\cdot \mathbb{P}\left(D_i>m\right)+\frac{1}{r} \mathbb{E}\left[\Lambda^{\left(D_i\right)} \cdot I\left(D_i \geq \tcp{m-B_{\delta,i}}\right)\right]+\frac{3 \delta}{\alpha_1}.
\end{equation}
Consider an arbitrary value of $\epsilon$ less than $(\beta_i-y)/2$, and recall the constant 
$
z_{y,\epsilon}
$
from \eqref{eq:z-ye}. We then have
\begin{equation}\label{eq:dm2}
\begin{aligned}
        \frac{1}{r} \mathbb{E}\left[\Lambda^{\left(D_i\right)} \cdot I\left(D_i \geq \tcp{m-B_{\delta,i}}\right)\right]\leq & \frac{1}{r} \mathbb{E}\left[\Lambda^{\left(D_i\right)} \cdot I\left(\frac{\Lambda^{(D_i)}}{r} > \frac{1}{2}(p_i+z_{y,\epsilon})\right)\right]\\&+ \frac{1}{2}(z_{y,\epsilon}+p_i)\cdot \mathbb{P}\left(D_i\geq \tcp{m-B_{\delta,i}}\right).
\end{aligned}
\end{equation}
Combining inequalities \eqref{eq:dm} and \eqref{eq:dm2} yields
\begin{equation*}
\begin{aligned}
        \frac{1}{2}(p_i-z_{y,\epsilon})\cdot \mathbb{P}\left(D_i>m\right)\leq& \frac{1}{2}(z_{y,\epsilon}+p_i)\cdot \mathbb{P}\left(\tcp{m-B_{\delta,i}} \leq D_i\leq m\right)\\&+\frac{1}{r} \mathbb{E}\left[\Lambda^{\left(D_i\right)} \cdot I\left(\frac{\Lambda^{\left(D_i\right)}}{r} > \frac{1}{2}(p_i+z_{y,\epsilon})\right)\right]+\frac{3\delta}{\alpha_1}.
\end{aligned}
\end{equation*}
For any $\tcp{N>B_{\delta,i}}$, there exists an integer $m$ with $1 \leq m \leq N$ such that
$$
\mathbb{P}\left(\tcp{m-B_{\delta,i}} \leq D_i \leq m\right) \leq \frac{\tcp{B_{\delta,i}}+1}{N} .
$$
Consequently, we have that for any $r,\delta$ and $\tcp{N>B_{\delta,i}}$, 
\begin{equation}\label{dmn}
    \frac{1}{2}(p_i-z_{y,\epsilon})\cdot \mathbb{P}\left(D_i>N\right)\leq \frac{1}{2}(z_{y,\epsilon}+p_i)\cdot\frac{\tcp{B_{\delta,i}}+1}{N}+\frac{1}{r} \mathbb{E}\left[\Lambda^{\left(D_i\right)} \cdot I\left(\frac{\Lambda^{(D_i)}}{r} > \frac{1}{2}(p_i+z_{y,\epsilon})\right)\right]+\frac{3\delta}{\alpha_1}.
\end{equation}
Hence for any $\eta>0$, \tcblue{taking $N=\lfloor\eta r\rfloor$ in \eqref{dmn} gives, for all sufficiently large $r$,} 
\begin{equation}\label{dmn2}
    \frac{1}{2}(p_i-z_{y,\epsilon})\cdot \mathbb{P}\left(\frac{D_i}{r}>\eta\right)\leq \frac{1}{2}(z_{y,\epsilon}+p_i)\cdot\frac{\tcp{B_{\delta,i}}+1}{\tcblue{\lfloor\eta r\rfloor}}+\frac{1}{r} \mathbb{E}\left[\Lambda^{\left(D_i\right)} \cdot I\left(\frac{\Lambda^{(D_i)}}{r}> \frac{1}{2}(p_i+z_{y,\epsilon})\right)\right]+\frac{3\delta}{\alpha_1}.
\end{equation}
\tcblue{For $i=1$, $\Lambda_1^{(D_1)}=0$. For $i\geq2$, the inductive hypothesis at $y$ gives} 
$$
 \Lambda_1^{\left(D_i\right)}/r=(\frac{2}{\alpha_i}+\frac{\tcblue{2}}{\alpha_i-\alpha_{i-1}})\biggr(\beta_{i-1}- \frac{1}{r} \sum_{\ell=1}^{i-1} \alpha_{\ell}F_{\ell}(r\beta_{i-1};\infty)\biggr)\xrightarrow{P} 0, \ \text{as} \ r \rightarrow \infty,
$$
and
$$
 \Lambda_3^{\left(D_i\right)}/r= \frac{\tcblue{4}}{\alpha_i}\biggr(  \left(y-\beta_{i-1}\right)- \frac{\alpha_i}{r}\big(F_i(ry;\infty) - F_i(r\beta_{i-1};\infty)\big)\tcblue{+\frac{\alpha_i}{r}}\biggr)+\frac{\left(y-\beta_{i-1}\right)}{\alpha_i}\xrightarrow{P}  \frac{\left(y-\beta_{i-1}\right)}{\alpha_i},  \ \text{as} \ r \rightarrow \infty.
$$
Therefore,
$$
 \Lambda^{\left(D_i\right)}/r= \Lambda_1^{\left(D_i\right)}/r + \Lambda_2^{\left(D_i\right)}/r+\Lambda_3^{\left(D_i\right)}/r\tcblue{+\frac{6B_{\delta,i}}{r}}\xrightarrow{P}  z_{y,\epsilon}  \ \text{as} \ r \rightarrow \infty.
$$
Thus by bounded convergence, 
\begin{equation}\label{eqlambdad}
    \lim_{r\rightarrow\infty} \frac{1}{r} \mathbb{E}\left[\Lambda^{\left(D_i\right)} \cdot I\left(\frac{\Lambda^{(D_i)}}{r}
 > \frac{1}{2}(p_i+z_{y,\epsilon})\right)\right]=0.
\end{equation}
Combining \eqref{eqlambdad} with inequality \eqref{dmn2}, \tcblue{first letting $r\rightarrow\infty$ with $\delta$ fixed and then letting $\delta\downarrow0$, we obtain}
$$
 \frac{D_i}{r}\xrightarrow{P}  0,  \ \text{as} \ r \rightarrow \infty.
$$
Applying bounded convergence once again, we obtain
\begin{equation*}
\lim_{r\rightarrow\infty}\frac{\mathbb{E}\left(D_i\right)}{r}=0.
\end{equation*}
\tcblue{Then, }combining \eqref{dmn} and \eqref{eqlambdad} yields
     $$\tcp{
      \limsup_{r\rightarrow\infty}\mathbb{P}(D_i>N)\leq  \frac{p_i+z_{y,\epsilon}}{p_i-z_{y,\epsilon}}\frac{B_{\delta,i}+1}{N}+\frac{6\delta}{\alpha_1(p_i-z_{y,\epsilon})} .
     }$$
\tcorange{For $1\leq N\leq B_{\delta,i}$, the bound is immediate
because its right-hand side is at least $1$.}
     
\end{proof}

Proposition \ref{npr1} implies that, for any $y \in\left[\beta_{i-1},\beta_i\right)$ and $\epsilon$ sufficiently small, under the inductive hypothesis at $y$,  in the interval $
\left[r \beta_{i-1}, r(y+\epsilon)\right)$, the hydrodynamic-scaled number of such ``large" holes (holes of size greater \tcp{than or equal to} $\alpha_i$) converges to zero as $r \rightarrow \infty$. However, in a system with countably many item types, with sizes that are not necessarily mutually rational, the arrivals and departures of different items induce fragmentation that produces a variety of ``small” holes (holes with size less than $\alpha_i$) with arbitrarily small positive sizes. In order to study how many holes of size less than $\alpha_i$ can appear, in the next subsection, we introduce a measure of fragmentation, which is an upper bound of all possible fragmentation patterns.

\tcblue{The next corollary is proved in the Appendix using a drift argument for a bounded functional measuring the distance from $r\beta_{i-1}$ to the left-most item left endpoint in $[r\beta_{i-1},\,r(y+\epsilon))$, with a convention for the case when no such endpoint exists.}
\begin{cor}\label{cor:gmin}
Fix $y\in[\beta_{i-1},\beta_i)$ and $0<\epsilon<(\beta_i-y)/2$ and assume the inductive hypothesis holds at $y$. Recall from \eqref{eq:jtilde-def} that $j_{\min}(y+\epsilon)$ denotes the index of the first item fully contained in
$
[r\beta_{i-1},\,r(y+\epsilon)),
$ and that $u_{j_{\min}(y+\epsilon)}$ is the left endpoint of the first item fully contained in
$
 [r\beta_{i-1},\,r(y+\epsilon)).
$
\tcbl{If there exists at least one item fully
contained in $[r\beta_{i-1},\,r(y+\epsilon))$}, let $G^{\min}_{i,y,\epsilon}$ be the total number of size-$\alpha_i$ items that can fit completely into the available empty space in
$$
[r\beta_{i-1},\,u_{j_{\min}(y+\epsilon)}).
$$
\tcbl{If no item is fully
contained in $[r\beta_{i-1},\,r(y+\epsilon))$, set}
$$
G^{\min}_{i,y,\epsilon}:=0.
$$
Then
\begin{equation}\label{eq:gmin}
\lim_{r\to\infty}\mathbb P\left(G^{\min}_{i,y,\epsilon}>0\right)=0.
\end{equation}
\end{cor}
\begin{proof}
See Appendix.
\end{proof}

Combining \eqref{tildeg} and \eqref{eq:gmin} implies that the aggregate rate with which a size-$\alpha_i$ arrival is placed entirely in  
$[0,\,u_{j_{\min}(y+\epsilon)})$  is $o(r)$.

\subsection{A characterization of fragmentation \tcorange{with respect to size-$\alpha_i$ items}}\label{sectiontf}
Given the above analysis, we know that the asymptotic dynamics on the half‐line $\left[r\beta_{i-1}, \infty\right)$ mainly involve the arrivals and departures of items of types $i, i+1, ...$ To study these dynamics we will introduce a ``measure of potential maximal total fragmentation'' of a subinterval of length $x$. 
In turn, to define a fragmentation measure, we must first understand how many distinct maximal ``admissible packing combinations" exist, which consist of items of types $j \ge i$ that can simultaneously fit within an interval of length $x$.

Let $\mathcal{J}_i(x)$ be the set of all maximal admissible packing combinations (using only items of types $j \geq i$ ) that can fit into an interval of length $x$. The formal definition of $\mathcal{J}_i(x)$ is given by Definition \ref{def1} below.

 \begin{definition}\label{def1}
     For any $x>\alpha_i$, let $\mathcal{J}_i(x)$ denote the collection of all finite index‐sequences satisfying
 \begin{equation}
  \mathcal{J}_i(x) = \left\{\, \mathbf{j} = (j_1, j_2, \dots, j_m) \in \{i, i+1, i+2, \ldots\}^m \ \bigg| \ m \in \mathbb{N}_{+},\ \forall \ 1\leq \ell \leq m, \ \alpha_{j_{\ell}}<x, \ x - \alpha_i < \sum_{\ell=1}^m \alpha_{j_\ell} \leq x \,\right\}.
 \end{equation}
  \end{definition}
  
Based on the construction of $\mathcal{J}_i(x)$, we define $H_i(x)$ for an interval of length $x$ as the maximum (over all admissible packing combinations) of the unused space length plus the sum of the recursive terms $H_i\left(\alpha_{j_\ell}\right)$ contributed by the packed items. We define $H_i(x)$ using the following recursive definition.
\begin{definition}
    For any $\tcp{x\geq 0}$,
    \begin{equation}
\begin{aligned}\label{H_i}
&H_i(x)= \begin{cases}x, & \tcp{0\leq x<\alpha_i} \\ 0, & x=\alpha_i \\ \max _{\mathbf{j} \in \tcblue{\mathcal{J}_i(x)}}\left[\left(x-\sum_{\ell=1}^m \alpha_{j_\ell}\right)+\sum_{\ell=1}^m H_i\left(\alpha_{j_\ell}\right)\right] & x>\alpha_i\end{cases}\\
\end{aligned}
\end{equation}
\end{definition}
We now pause to make a few remarks about $H_i$.
\begin{itemize}
  \item If $\tcp{0\leq x < \alpha_i}$, no item of size $\alpha_i$ or larger can fit, so the entire interval of length $x$ remains as holes. In this case,
  $
    H_i(x) = x.\,
  $\tcp{Note that $H_i(0)=0$.}
  
  \item If $x = \alpha_i$, we can place exactly one item of size $\alpha_i$ and leave zero residual. In this case,
  $
    H_i(\alpha_i) = 0.
  $
  
  \item For $x > \alpha_i$, consider every admissible sequence $\mathbf{j}\in \mathcal{J}_i(x)$. Placing items of sizes $\{\alpha_{j_\ell}\}$ consumes
  $
    \sum_{\ell=1}^m \alpha_{j_\ell},
  $
  leaving a leftover of $x - \sum_{\ell=1}^m \alpha_{j_\ell}$, which lies in $\tcblue{[0,\alpha_i)}$. Additionally, each \tcorange{interval occupied by an} item of size $\alpha_{j_\ell}$ may itself generate smaller holes of total length $H_i(\alpha_{j_\ell})$. Hence the total residual under pattern $\mathbf{j}$ is
  $$
    \parens{x - \sum_{\ell=1}^m \alpha_{j_\ell}}
    + \sum_{\ell=1}^m H_i\bigl(\alpha_{j_\ell}\bigr).
  $$
  Therefore, $H_i(x)$ is defined as the maximum of this quantity over all $\mathbf{j}\in \mathcal{J}_i(x)$.
\end{itemize}

Proposition \ref{hyx} below illustrates a fundamental property of $H_i$.
\begin{prop}\label{hyx}
    For any $0<x<y$, \begin{equation}\label{eq:hyx}
        H_i(y)-H_i(x)\leq y-x.
    \end{equation}
\end{prop}

\begin{proof} 
 Without loss of generality, assume $0\le x<y$ and $0<y - x<\alpha_i$. If $y<\alpha_i$, then by the definition of $H_i$ on $[0,\alpha_i)$ we have $H_i(y)-H_i(x)=y-x$. \tcblue{If $y=\alpha_i$, \eqref{eq:hyx} follows directly from the definition of $H_i$.} 
We therefore focus on the case $y>\alpha_i$. Let $\mathbf{j}=(j_1,\dots,j_m)\in \mathcal{J}_i(y)$ be any sequence attaining the maximum in the definition of $H_i(y)$, so that
$$
H_i(y)
=\Bigl(y-\sum_{\ell=1}^m\alpha_{j_\ell}\Bigr)
+\sum_{\ell=1}^m H_i(\alpha_{j_\ell}).
$$
\tcblue{Since the order of the indices does not affect the two sums above, if some $j_\ell=i$, we reorder the indices so that $j_m=i$.}
If $\sum_{\ell=1}^m\alpha_{j_\ell}\leq x$ \tcblue{and $m=1$ with $\alpha_{j_1}=x$}, then \tcblue{$H_i(y)=y-x+H_i(x)$, and \eqref{eq:hyx} follows immediately.}
\tcblue{Otherwise, if $\sum_{\ell=1}^m\alpha_{j_\ell}\leq x$, every $\alpha_{j_\ell}$ is strictly less than $x$, and}
$$
\tcblue{x-\alpha_i<y-\alpha_i<\sum_{\ell=1}^m\alpha_{j_\ell}\leq x.}
$$
\tcblue{Thus, the same sequence belongs to $\mathcal{J}_i(x)$, and the definition of $H_i(x)$ gives}
$$
\begin{aligned}
H_i(x)
\ge\sum_{\ell=1}^m H_i(\alpha_{j_\ell})+x-\sum_{\ell=1}^m\alpha_{j_\ell}
=\Bigl(y-\sum_{\ell=1}^m\alpha_{j_\ell}\Bigr)
+\sum_{\ell=1}^m H_i(\alpha_{j_\ell})-(y-x)
=H_i(y)-(y-x).
\end{aligned}
$$
If $\sum_{\ell=1}^m\alpha_{j_\ell}> x$ and $j_m=i$, then
$$
\tcblue{x-\alpha_i<\sum_{\ell=1}^{m-1}\alpha_{j_\ell}\leq y-\alpha_i<x.}
$$
\tcblue{If $m>1$, the sequence $(j_1,\dots,j_{m-1})$ belongs to $\mathcal{J}_i(x)$. If $m=1$, then $x<\alpha_i$, and the first inequality below follows from $H_i(x)=x$. Thus,}
$$
\begin{aligned}
H_i(x)
\ge\sum_{\ell=1}^{m-1}H_i(\alpha_{j_\ell})+x-\sum_{\ell=1}^{m-1}\alpha_{j_\ell}
\geq \Bigl(y-\sum_{\ell=1}^m\alpha_{j_\ell}\Bigr)
+\sum_{\ell=1}^m H_i(\alpha_{j_\ell})-(y-x)
=H_i(y)-(y-x).
\end{aligned}
$$
If $\sum_{\ell=1}^m\alpha_{j_\ell}> x$, $j_m\neq i$, \tcblue{and $m>1$}, then \tcblue{$\alpha_{j_m}>\alpha_i$ and}
$$
\tcblue{\sum_{\ell=1}^{m-1}\alpha_{j_\ell}\leq y-\alpha_{j_m}<y-\alpha_i<x.}
$$
\tcblue{Therefore, }by \tcblue{the} definitions of $H_i(x)$ and $H_i\Bigr(y-\sum_{\ell=1}^{m-1}\alpha_{j_\ell}\Bigr)$, we have
$$
H_i(x)\ge H_i\Bigl(x-\sum_{\ell=1}^{m-1}\alpha_{j_\ell}\Bigr)
+\sum_{\ell=1}^{m-1}H_i(\alpha_{j_\ell}),
$$
and, since $y-\sum_{\ell=1}^{m-1}\alpha_{j_\ell}-\alpha_{j_m}<\alpha_i$,
$$
H_i\Bigr(y-\sum_{\ell=1}^{m-1}\alpha_{j_\ell}\Bigr)\ge \Bigl(y-\sum_{\ell=1}^{m}\alpha_{j_\ell}\Bigr)
+H_i(\alpha_{j_m})=H_i(y)-\sum_{\ell=1}^{m-1}H_i(\alpha_{j_\ell}).
$$
Hence, to prove $H_i(x) \geq H_i(y)-(y-x)$,
it suffices to show
$$
H_i\Bigl(x-\sum_{\ell=1}^{m-1}\alpha_{j_\ell}\Bigr)
\ge H_i\Bigl(y-\sum_{\ell=1}^{m-1}\alpha_{j_\ell}\Bigr)-(y-x).
$$
Setting
$$
x^{\prime}=x-\sum_{\ell=1}^{m-1} \alpha_{j_\ell}, \quad y^{\prime}=y-\sum_{\ell=1}^{m-1} \alpha_{j_\ell},
$$
we still have $0 < y' - x' < \alpha_i$.

Finally, if $\sum_{\ell=1}^m\alpha_{j_\ell}>x$, $j_m\neq i$ and $m=1$, then
$H_i(y)=y-\alpha_{j_1}+H_i(\alpha_{j_1})$, so it suffices to show
$$H_i(x)\ge H_i(\alpha_{j_1})-(\alpha_{j_1}-x).$$ 
In \tcblue{this} case, setting
$$
x^{\prime}=x, \quad y^{\prime}=\alpha_{j_1},
$$
we also have $0 < y' - x' < \alpha_i$.  

In the last two cases, we therefore have a pair $(x',y')$ with \tcblue{$0<x'<y'$ and }$0<y'-x'<\alpha_i$.  \tcblue{If \eqref{eq:hyx} for $(x',y')$ does not follow from the definition of $H_i$ or from one of the first two cases above, choose a sequence $(j_1,\dots,j_m)\in\mathcal{J}_i(y')$ attaining the maximum in the definition of $H_i(y')$, reorder it as above, and apply the same reduction.}
More precisely, when $m>1$, \tcblue{the same argument as above gives $\sum_{\ell=1}^{m-1}\alpha_{j_\ell}<x'$, and} we replace $(x',y')$ by
$$
\left(x'-\sum_{\ell=1}^{m-1}\alpha_{j_\ell},\,
      y'-\sum_{\ell=1}^{m-1}\alpha_{j_\ell}\right),
$$
and when $m=1$, \tcblue{we have $x'<\alpha_{j_1}<y'$, and} we replace $(x',y')$ by
$$
(x',\alpha_{j_1}).
$$
In both cases, the new pair, denoted by $(x'',y'')$, still satisfies
$$
\tcblue{0<x''<y'',\qquad}0<y''-x''<\alpha_i.
$$
\tcblue{At each subsequent step, if \eqref{eq:hyx} for the new pair does not follow from the definition of $H_i$ or from one of the first two cases above, we relabel $(x'',y'')$ as $(x',y')$, choose and reorder a sequence attaining the maximum in the definition of $H_i(y')$, and repeat. If $m>1$, $y'$ decreases by at least $\alpha_i$, and if $m=1$, $y'$ is replaced by the strictly smaller item size $\alpha_{j_1}$. Since $y'$ remains positive and only finitely many item sizes are smaller than $y$, the procedure terminates after finitely many steps.}
Applying the corresponding bound \tcblue{at the final step} and then substituting back through the previous reductions yields
$$
H_i(x)\ge H_i(y)-(y-x),
$$
as required.
\end{proof}
A useful corollary to Proposition \ref{hyx} is as follows.
\begin{cor}\label{hyx2}
    For any $j^*\in \{i,i+1,i+2  \dots\}$ and any $x\in (0,\alpha_i)$, 
    \begin{equation*}
        H_i(\alpha_{j^*}+x)=H_i(\alpha_{j^*})+x.
    \end{equation*}  
    \end{cor}
\begin{proof}[Proof of Corollary \ref{hyx2}]
By the definition of $H_i$, for any $0<x<\alpha_i$,
$$
H_i(\alpha_{j^*}+x)\;\ge\;H_i(\alpha_{j^*})\;+\;x.
$$
On the other hand, Proposition \ref{hyx} implies that
$$
H_i(\alpha_{j^*}+x)\;\le\;H_i(\alpha_{j^*})\;+\;x.
$$
Hence, the corollary immediately follows by combining these two inequalities.
\end{proof}
With the above preparations in hand,  we now introduce a functional $\mathrm{TF}$ of a system state $S$ to capture the total fragmentation of (a part of) the interval $\left[r \beta_{i-1}, r  \beta_{i}\right)$, with respect to size-$\alpha_i$ items.
\tcgr{Recall the notation $e_j$ and $g_j$ in \eqref{eq:eg-def}. For $y\in[\beta_{i-1},\beta_i)$, whenever there exists at least one item fully contained in $[r\beta_{i-1},\,ry)$, let $j_{\min}(y)$ and $j_{\max}(y)$ be the indices defined in \eqref{eq:jtilde-def}. Thus $u_{j_{\min}(y)}$ is the left endpoint of the first item fully contained in $[r\beta_{i-1},\,ry)$, and $u_{j_{\max}(y)}$ is the left endpoint of the last item fully contained in $[r\beta_{i-1},\,ry)$.}
\begin{definition}\label{def:tf}
For any $
y \in \left[\beta_{i-1}, \beta_{i}\right),
$ we define the total fragmentation $
\mathrm{TF}(ry)
$ as
\begin{equation}
\begin{aligned}\label{eq:tf}
        \mathrm{TF}(ry)=&\sum_{\substack{\tcgr{g_j<\alpha_i} \\ v_{j+1} \tcorange{\leq} ry}} H_i\left(\tcgr{e_j+g_j}\right)+\sum_{\substack{\tcgr{g_j \geq \alpha_i} \\ v_{j+1} \tcorange{\leq} ry}}\left[H_i\left(\tcgr{e_j}\right)+H_i\left(\tcgr{g_j}\right)\right] \\
        &+\tcp{\left(u_{j_{\min}(y)}-r\beta_{i-1}+ry-u_{j_{\max}(y)}\right).}
\end{aligned}
\end{equation}
In particular, if there is no item fully contained in the interval $\left[r \beta_{i-1}, r y\right)$, then $\mathrm{TF}(r y)=r y-r \beta_{i-1}$. 
\end{definition}

\tcbl{Note that, via $y\in[\beta_{i-1},\beta_i)$, $\mathrm{TF}(ry)$ implicitly depends on $i$ as a parameter. Some general remarks about \eqref{eq:tf} are now in order. In \eqref{eq:tf}, the two sums treat the cases $g_j < \alpha_i$ and the opposite, differently. Note that the first sum considers the fragmentation of the entire interval containing the item plus the hole whereas the second sum considers the fragmentation of the item interval plus the fragmentation of the hole interval. We shall refer to the third term in \eqref{eq:tf} as the ``boundary" term. This term accounts for the ``worst case" fragmentation contributions adjacent to the endpoints $r\beta_{i-1}$ and $ry$}: we add $u_{\tcp{j_{\min}(y)}}-r\beta_{i-1}$, the distance from $r\beta_{i-1}$ to the left endpoint of the first item fully contained in $[r\beta_{i-1},ry)$, and $ry-u_{\tcp{j_{\max}(y)}}$, the distance from $ry$ to the left endpoint of the last item fully contained in $[r\beta_{i-1},ry)$. 

With the definition of $\mathrm{TF}(ry)$ in place, $\mathrm{TF}(ry)/r$ tends to zero in probability as $r\to\infty$ when the system is asymptotically optimally packed on $[r\beta_{i-1},\, ry)$.  \tcbl{When considering the interval $[r\beta_{i-1},\,r(y+\epsilon))$, under the inductive assumption at $y$, we typically see items of size at least $\alpha_i$ and holes of length less than $\alpha_i$. This is the same picture as in Sections~2 and~3: in the two-size case, the typical packing configuration treated in Subsection~3.5 for the interval $[p_1r,\,\tcorange{r(y+\epsilon)})$ consists of 2-items separated by holes of length $1$. In this more general case, Proposition~\ref{prop31} shows that arrivals of items of types $j<i$ into $[r\beta_{i-1},\,r(y+\epsilon))$ have total rate $o(r)$ and Proposition~\ref{npr1} shows that the hydrodynamic-scaled number of holes of length at least $\alpha_i$ in $[r\beta_{i-1},\,r(y+\epsilon))$ converges to zero. Therefore, arrivals of items of types $j<i$ into $[r\beta_{i-1},\,r(y+\epsilon))$ and departures of items with an adjacent hole of length at least $\alpha_i$ are not included among the ``typical" events considered below, since the total rate of such events tends to $0$ on the hydrodynamic scale. Thus, Events~1 and~2 in Proposition \ref{prop632} below are the ``typical" arrival and departure events we shall consider in the interval $[r\beta_{i-1},\,r(y+\epsilon))$.}

\tcbl{The proposition below states that the total fragmentation $\mathrm{TF}$ remains non‐increasing during each of the system’s ``typical" events.} 
\begin{prop}\label{prop632}
    For any $y\in\left[\beta_{i-1},\beta_i\right)$,
$\mathrm{TF}(ry)$ remains non‐increasing on the following events.
\begin{itemize}
  \item \textbf{Event 1}: For any $k\geq i$, a size-$\alpha_k$ item arrives in the interval 
  $\bigl[r\beta_{i-1},\;ry \bigr)$.
  
  \item \textbf{Event 2}: For any $k\in \mathbb{N}_{+}$, a size-$\alpha_k$ item completely lying in $\bigl[r\beta_{i-1}, ry \bigr)$, but not the left-most or the right-most such item, departs from the system,  and immediately to its left and to its right there are no holes of length greater \tcp{than or equal to} $\alpha_i$. In  addition, the item immediately to its left (possibly separated by a hole) has size at least $\alpha_i$.
\end{itemize}
\end{prop}
\begin{proof}[Proof of Proposition \ref{prop632}]
\tcgr{Consider Event 1. \tcblue{If no item is fully contained in $[r\beta_{i-1},ry)$ before the arrival, then after the arrival there is at most one item fully contained in this interval. By Definition~\ref{def:tf}, we have $\mathrm{TF}(ry)=r(y-\beta_{i-1})$ both before and after the arrival, so $\mathrm{TF}(ry)$ does not increase.} Suppose first that a size-$\alpha_k$ item with $k\geq i$ arrives and occupies the hole $[v_j,u_{j+1})$, where both adjacent items $[u_j,v_j)$ and $[u_{j+1},v_{j+1})$ are fully contained in $[r\beta_{i-1},ry)$. Then $g_j\geq \alpha_k$. Since $g_j\geq \alpha_i$, before the arrival, the term indexed by $j$ in the second sum of Definition~\ref{def:tf} is $H_i(e_j)+H_i(g_j)$. After the arrival, the item $[u_j,v_j)$ is unchanged and the arriving size-$\alpha_k$ item contributes $H_i(\alpha_k)$. The remaining hole has length $g_j-\alpha_k$. If $g_j-\alpha_k\geq \alpha_i$, then the new term is $H_i(e_j)+H_i(\alpha_k)+H_i(g_j-\alpha_k)$. If $g_j-\alpha_k<\alpha_i$, Corollary~\ref{hyx2} yields}
$$
\tcgr{H_i(g_j)=H_i(\alpha_k)+H_i(g_j-\alpha_k).}
$$
\tcgr{Thus the post-arrival expression can be written as $H_i(e_j)+H_i(\alpha_k)+H_i(g_j-\alpha_k)$. Therefore,}
$$
\tcgr{\Delta \mathrm{TF}=-H_i(g_j)+H_i(g_j-\alpha_k)+H_i(\alpha_k).}
$$
\tcgr{\tcbl{We now show that $\Delta \mathrm{TF}$ is non-positive.} Let $x:=g_j-\alpha_k$. In Definition~\ref{def1}, one admissible packing of an interval of length $x+\alpha_k$ is obtained by placing one size-$\alpha_k$ item and then using an admissible packing for the remaining interval of length $x$. Therefore}
$$
\tcgr{H_i(g_j)=H_i(x+\alpha_k)\geq H_i(\alpha_k)+H_i(x)=H_i(\alpha_k)+H_i(g_j-\alpha_k),}
$$
\tcgr{and hence $\Delta \mathrm{TF}\leq 0$.}

\tcorange{Suppose instead that the arriving item is placed to the
left of the first item fully contained in $[r\beta_{i-1},ry)$
before the arrival. The interval
$[r\beta_{i-1},u_{j_{\min}(y)})$ is then divided into the new
left boundary interval and the interval consisting of the
arriving item and the hole on its right. By
Definition~\ref{def:tf} and $H_i(z)\leq z$ for $z\geq0$,
the sum of their contributions is at most
$u_{j_{\min}(y)}-r\beta_{i-1}$, so $\mathrm{TF}(ry)$
does not increase.}

\tcorange{If the arriving item is placed to the right of the
last item fully contained in $[r\beta_{i-1},ry)$ before the
arrival, the first-fit rule places its left endpoint at
$v_{j_{\max}(y)}$. The contribution $ry-u_{j_{\max}(y)}$
is then replaced by
$H_i(e_{j_{\max}(y)})+ry-v_{j_{\max}(y)}$.
Hence
$\Delta\mathrm{TF}
=H_i(e_{j_{\max}(y)})-e_{j_{\max}(y)}\leq0$.}

\tcgr{Now consider Event~2. Before the departure, we have $g_{j-1}<\alpha_i$, $g_j<\alpha_i$, and $e_{j-1}\geq \alpha_i$. Therefore the two terms in Definition~\ref{def:tf} that involve the item $[u_j,v_j)$ are $H_i(e_{j-1}+g_{j-1})$ and $H_i(e_j+g_j)$. After the departure, the two adjacent holes and the interval $[u_j,v_j)$ merge into a single hole of length $g_{j-1}+e_j+g_j$.} \tcbl{The corresponding contribution to $\mathrm{TF}$ from the item of length $e_{j-1}$ and this merged hole is}
$$
\tcbl{H_i(e_{j-1})+H_i(g_{j-1}+e_j+g_j).}
$$
\tcbl{We see this by considering two cases. If $g_{j-1}+e_j+g_j\ge \alpha_i$, then this contribution is counted by the second sum in Definition \ref{def:tf}. If $g_{j-1}+e_j+g_j<\alpha_i$, then Corollary~\ref{hyx2} gives}
$$
\tcbl{H_i(e_{j-1}+g_{j-1}+e_j+g_j)=H_i(e_{j-1})+H_i(g_{j-1}+e_j+g_j).}
$$
\tcgr{Hence}
$$
\tcgr{\begin{aligned}
\Delta \mathrm{TF}
&=H_i(e_{j-1})+H_i(g_{j-1}+e_j+g_j)-H_i(e_{j-1}+g_{j-1})-H_i(e_j+g_j)\\
&=H_i(g_{j-1}+e_j+g_j)-g_{j-1}-H_i(e_j+g_j)\leq 0,
\end{aligned}}
$$
\tcgr{where the second line uses Corollary~\ref{hyx2}, and the last inequality follows from Proposition~\ref{hyx} applied to $x=e_j+g_j$ and $y=g_{j-1}+e_j+g_j$. This completes the proof.}
\end{proof}

In the next proposition, we show that under the inductive hypothesis at $y\in[\beta_{i-1},\beta_i)$ and for any $0<\epsilon< (\beta_i-y)/2$, the expectation of the aggregate contributions of the scenarios which may increase $
\mathrm{TF}(r(y+\epsilon))$ is $o(r)$.
\begin{prop}\label{prop636}
    For $y\in[\beta_{i-1},\beta_i)$ and $0<\epsilon<(\beta_i-y)/2$, let $\Lambda^{\mathrm{TF}}$ denote the instantaneous aggregate rate at which $\mathrm{TF}(r(y+\epsilon))$ may increase. Specifically, $\Lambda^{\mathrm{TF}}$ is the sum of the terms in the generator $\mathcal{A}  \mathrm{TF}(r(y+\epsilon))$ corresponding to transitions that increase $\mathrm{TF}(r(y+\epsilon))$. Then, under the inductive hypothesis at $y$, 
    \begin{equation}\label{tf}
\lim_{r\rightarrow\infty}\frac{\mathbb{E}\left(\Lambda^{\mathrm{TF}}(\infty)\right)}{r}=0.
    \end{equation}
\end{prop}
\begin{proof}[Proof of Proposition \ref{prop636}]
By Proposition \ref{prop632}, $\mathrm{TF}(r(y+\epsilon))$ can only increase in three scenarios: (i) when an item of size less than $\alpha_i$ arrives into $[r\beta_{i-1},r(y+\epsilon))$; (ii) when an item (but not the left-most or the right-most such item) lying in $[r\beta_{i-1},\,r(y+\epsilon))$
departs and either is adjacent to a hole of length at least $\alpha_i$ or its nearest left neighbor has size less than $\alpha_i$; (iii) when the first or the last item fully lying in $\left[r \beta_{i-1}, r (y+\epsilon)\right)$ departs. 

In the first case, the rate of arrival is at most $r\cdot\sum_{j=1}^{i-1}I(G_j=0)$. Each such arrival increases $\mathrm{TF}(r(y+\epsilon))$ by at most $\alpha_i$, and so the contribution to the rate of increase is bounded by $r\alpha_i\cdot\sum_{j=1}^{i-1}I(G_j=0)$.

\tcblue{For a departure in the second case, by Proposition~\ref{hyx} and Corollary~\ref{hyx2}, the positive change of $\mathrm{TF}(r(y+\epsilon))$ is at most}
$$
\tcblue{
2\alpha_i\left(
\left\lfloor\frac{g_{j-1}}{\alpha_i}\right\rfloor+
\left\lfloor\frac{g_j}{\alpha_i}\right\rfloor
\right)
+\alpha_i I(e_{j-1}<\alpha_i).
}
$$
\tcblue{Each hole can appear as $g_{j-1}$ or $g_j$ for at most two departures, and each item can appear as $e_{j-1}$ for at most one departure. }   Hence, the contribution to the rate of increase in the second case is bounded \tcblue{above} by
\begin{equation}\label{disum}
    4\alpha_i\left(D_i+\sum_{\ell<i}\left(F_{\ell}(r(y+\epsilon))-F_{\ell}\left(r \beta_{i-1}\right)\right)\right). 
\end{equation}
For the third case, \tcp{$\mathrm{TF}(r (y+\epsilon))$} increases by at most the distance between the left endpoints of the first and the second items, and the distance between the left endpoints of the last and the second-to-last items in \tcp{$\left[r \beta_{i-1}, r (y+\epsilon)\right)$}. Therefore, for any $\delta>0$, the contribution to the rate at which \tcp{$\mathrm{TF}(r(y+\epsilon))$ may increase} is upper bounded by
$$
4\alpha_i(D_i+1)+2(A_{\delta}+\sum_{\ell: \alpha_{\ell}>A_\delta}\alpha_{\ell} F_{\ell}(\infty ; \infty)),
$$
\tcp{where $4\alpha_i(D_i+1)$ bounds the contribution from the two boundary holes (i.e., the holes between the first and second items, and between the second-to-last and last items) in $[r\beta_{i-1},r(y+\epsilon))$, while $2\big(A_\delta+\sum_{\ell:\alpha_\ell>A_\delta}\alpha_\ell F_\ell(\infty;\infty)\big)$ bounds the contribution from the lengths of the two boundary items.}

The above analysis allows us to write \begin{equation*} 
\begin{aligned}
  \limsup_{r\rightarrow\infty}\frac{\mathbb{E}\left(\Lambda^{\mathrm{TF}}(\infty)\right)}{r}\leq & \alpha_i \lim_{r\rightarrow\infty} \sum_{j=1}^{i-1}\mathbb{P}\left(G_j=0\right)+ 8\alpha_i\lim_{r\rightarrow\infty} \frac{\mathbb{E}\left(D_i+A_{\delta}+1\right)}{r}+\limsup _{r \rightarrow \infty} \mathbb{E}\left(\frac{2}{r} \sum_{\ell: \alpha_{\ell}>A_\delta}\alpha_{\ell} F_{\ell}(\infty ; \infty)\right)\\
  &+\limsup _{r \rightarrow \infty} \mathbb{E}\left(\frac{4\alpha_i }{r} \sum_{\ell<i}\left(F_{\ell}(r(y+\epsilon))-F_{\ell}\left(r \beta_{i-1}\right)\right) \right) <2\delta.  
\end{aligned}
 \end{equation*}
Hence, \begin{equation*} \lim_{r\rightarrow\infty}\frac{\mathbb{E}\left(\Lambda^{\mathrm{TF}}(\infty)\right)}{r}=0.
\end{equation*}
\end{proof}
Equality \eqref{tf} states that the hydrodynamic-scaled average rate at which $\mathrm{TF}(r(y+\epsilon))$ grows will vanish. Intuitively, this observation suggests that  under hydrodynamic scaling, and in the limit, the system’s total fragmentation is non-increasing. Consequently, since all mechanisms that may increase $\mathrm{TF}(r(y+\epsilon))$ contribute only $o(r)$ to the overall rate of change of $\mathrm{TF}(r(y+\epsilon))$, in the steady state the contribution of any mechanism that could decrease $\mathrm{TF}(r(y+\epsilon))$ must also be $o(r)$. We will  prove asymptotic optimality by showing that if all such potential decreases are $o(r)$, then the system must be ``asymptotically optimally packed" (recall this was defined after Theorem \ref{thm-ff1}).  We proceed to prove Theorem \ref{thm-ff2} in Section \ref{sec55}.

\section{Proof of Theorem \ref{thm-ff2}}\label{sec55}

This section builds on the analysis in Section \ref{sectionprecoro}, uses the notation introduced in \tcgr{Section~\ref{sectionroadmap}}, and uses the properties of $\mathrm{TF}(ry)$ established above to prove Theorem \ref{thm-ff2}. For any $y\in[\beta_{i-1},\beta_i)$ and sufficiently small $\epsilon>0$, we shall use $\mathrm{TF}(r(y+\epsilon))$ as a Lyapunov functional. For any $x\tcorange{\geq}0$, define
\begin{equation}\label{def:hi}
\tcbl{h_i(x)=x-\left\lfloor\frac{x}{\alpha_i}\right\rfloor \alpha_i.}
\end{equation}
From \eqref{def:hi}, for  $m \in \mathbb Z_{+}$, we have
\begin{equation}\label{eq:Hhi}
    \tcp{h_i(x)=H_i \parens{x-\left\lfloor\frac{x}{\alpha_i}\right\rfloor \alpha_i}}=\min _{m \geq 0, \ m \alpha_i\leq x} H_i\parens{x-m \alpha_i} \leq H_i(x).
\end{equation}
Therefore, any hole or item of length $x$ falls into one of two categories: (i) $H_i(x)>h_i(x)$; and (ii) $H_i(x)=h_i(x)$. 

\tcp{For convenience, in the sequel, for a  quantity $X_r$ depending on $r$, we shall write
$$
X_r=\widetilde o(r)
$$
to mean that
$$
\mathbb E[X_r]=o(r), \qquad r\to\infty.
$$
Additionally, for fixed $y \in\left[\beta_{i-1}, \beta_i\right)$ and $0<\epsilon<\left(\beta_i-y\right) / 2$, we define the positive constant
\begin{equation}\label{ciye}
    c_{i, y, \epsilon}:=p_i e^{-\left((y+\epsilon) / \alpha_1+2\right)} .
\end{equation}
\tcbl{Recall from Proposition~\ref{prop31} and Corollary~\ref{cor:gmin} that 
$\widetilde G_i$ and $G^{\min}_{i,y,\epsilon}$ are the numbers of size-$\alpha_i$ items that can fit, respectively, in $[0,r\beta_{i-1}+\alpha_i)$ and in the empty space before the first item fully contained in $[r\beta_{i-1},r(y+\epsilon))$}. 
}\tcp{We also define the event}
$$
\tcp{
B^{\min}_{i,y,\epsilon}
:=
\left\{G^{\min}_{i,y,\epsilon}=0\right\}
\cap
\left\{\widetilde G_i=0\right\}.
}
$$
\tcp{By \eqref{tildeg} and \eqref{eq:gmin},}
\begin{equation}\label{eq:bmin}
    \tcp{
\mathbb P\left(B^{\min}_{i,y,\epsilon}\right)
\ge
1-\mathbb P\left(\widetilde G_i>0\right)-\mathbb P\left(G^{\min}_{i,y,\epsilon}>0\right)
=
1-o(1).
}
\end{equation}
\tcgr{Before proceeding, we recall and expand on the roadmap from Section~\ref{sectionroadmap}, giving a more detailed outline of the proof of Theorem~\ref{thm-ff2}.} \tcbl{At a high level, the proof considers the following four \tcbl{classes} of items and holes: items of types $k>i$ with $H_i(\alpha_k)>h_i(\alpha_k)$, items of types $k>i$ with $H_i(\alpha_k)=h_i(\alpha_k)>0$, holes with length at least a fixed positive number, and items of types $k>i$ with $H_i(\alpha_k)=0$. For each \tcbl{class}, we use a Lyapunov functional and prove that if items or holes of that \tcbl{class} are present in quantity of order $r$, then the drift of the corresponding Lyapunov functional is at most $-cr$ for some constant $c>0$, while all terms that increase that functional are $o(r)$ under the inductive hypothesis at $y$.}

\tcgr{
More precisely, under the inductive hypothesis at $y$, we use $\operatorname{TF}(r(y+\epsilon))$ as the basis for the Lyapunov functionals. For each class of packing configurations $B$, we consider a Lyapunov functional built from $\operatorname{TF}(r(y+\epsilon))$ together with a count such as $\left|\widetilde{\mathcal{R}}_i^\delta\right|$, $\left|L_x^{\delta,n}\right|$, or $Q_i$. 
We then show that all terms that increase \tcgr{the chosen functional} are $o(r)$, while if $|B|$ is of order $r$, the first-fit dynamics decrease the functional at rate at least $cr$ for some $c>0$. Since the steady-state drift is zero, this implies that $|B|=o(r)$.} \tcp{To carry out this strategy, we divide the item types $k>i$ into three classes distinguished by the relationship between $H_i(\alpha_k)$ and $h_i(\alpha_k)$:
$$
H_i(\alpha_k)>h_i(\alpha_k), \quad
H_i(\alpha_k)=h_i(\alpha_k)>0, \quad
H_i(\alpha_k)=0.
$$
The proof proceeds by showing that the hydrodynamic-scaled numbers of items in each of these three classes, as well as the hydrodynamic-scaled total length of the holes within $[r\beta_{i-1},r(y+\epsilon))$, converge to zero as $r\to\infty$.
\begin{enumerate}
    \item  Firstly, Propositions~\ref{prop642}--\ref{prop643} and Corollary~\ref{corr644} show that, for item types in the first class, the hydrodynamic-scaled numbers of such items in $[r\beta_{i-1},r(y+\epsilon))$ converge to zero as $r\to\infty$.
   \item Secondly, Propositions~\ref{prop645}--\ref{prop646} show that, for every fixed $x>0$, the hydrodynamic-scaled number of holes of length at least $x$ in $[r\beta_{i-1},r(y+\epsilon))$ vanishes as $r\to\infty$. This strengthens the conclusion of Proposition~\ref{npr1}, which controls holes of length at least $\alpha_i$, and in turn allows us to treat the second class. Specifically, Corollary~\ref{prop647} shows that, for item types in the second class, the hydrodynamic-scaled numbers of such items in $[r\beta_{i-1},r(y+\epsilon))$ also vanish asymptotically.
   \item Furthermore, Proposition~\ref{prop648} shows that the hydrodynamically scaled total length of the holes in $[r\beta_{i-1},r(y+\epsilon))$ converges to zero.
   \item Finally, Proposition~\ref{prop649} shows that, for item types in the third class, the hydrodynamic-scaled numbers of such items in $[r\beta_{i-1},r(y+\epsilon))$ converge to zero as $r\to\infty$.  Combining these results, we conclude that under the inductive hypothesis at $y$,  the interval $[r\beta_{i-1}, r(y+\epsilon))$ is asymptotically optimally packed.
\end{enumerate}
}
\tcgr{We recall the notation $e_j$ and $g_j$ from \eqref{eq:eg-def}. Here, $[v_j,u_{j+1})$ is the hole between the consecutive items $[u_j,v_j)$ and $[u_{j+1},v_{j+1})$, and its length is $g_j$ (possibly $0$). Throughout Section~\ref{sec55}, all sets, item lengths, and hole lengths are interpreted with respect to the pre-event indexing. When an item departs, any merged hole is described using its pre-event endpoints; for example, we write $[v_j,u_{j+2})$ without relabeling indices during the argument.}

Some further notation is now necessary. When there exists at least one hole of length at least $\alpha_i$ whose two adjacent items both lie completely in $[r\beta_{i-1},\,r(y+\epsilon))$, we shall set $j_i^*$ to be the smallest index such that both 1. and 2. below hold:
\begin{enumerate}
    \item $\tcgr{g_j}\ge \alpha_i$.
    \item The items $[u_j,v_j)$ and $[u_{j+1},v_{j+1})$ are contained completely in $[r\beta_{i-1},\,r(y+\epsilon))$. This means that $[v_{j_i^*},u_{j_i^*+1})$ is the first hole of length at least $\alpha_i$ with adjacent items lying entirely in $[r\beta_{i-1},\,r(y+\epsilon))$.
\end{enumerate}\tcp{
We now introduce the notation
$FH_i := \tcgr{g_{j_i^*}}$ to denote the length of such a hole. If no such hole exists, we set $j_i^*=\infty$ and $FH_i=0$.}  \tcgr{The quantities $j_i^*$ and $FH_i$ are determined by the current state. Thus, whenever we write $j_i^*$ or $FH_i$ at a given time, these quantities are understood in the state at that time.}

\tcp{Suppose that $FH_i$ satisfies $H_i\left(FH_i\right)>h_i\left(FH_i\right)$. Then, when there are consecutive arrivals of size-$\alpha_i$ items into the hole $\left[v_{j_i^*},u_{j_i^*+1}\right)$ and no departures of other items in $[r\beta_{i-1}, r(y+\epsilon))$, successively filling the hole with size-$\alpha_i$ items until its residual length is strictly less than $\alpha_i$ will reduce $\mathrm{TF}(r(y+\epsilon))$ by $H_i\left(FH_i\right)-h_i\left(FH_i\right)$. This key observation provides the essential intuition for Proposition \ref{prop642}.
}
\begin{prop}\label{prop642}
     For any $\delta >0$, any $y\in[\beta_{i-1},\beta_i)$ and $0<\epsilon< (\beta_{i}-y)/2 $, under the inductive hypothesis at $y$, 
\begin{equation}\tcp{
   \lim _{r \rightarrow \infty} \mathbb{P}\left(H_i(FH_i)>h_i(FH_i), FH_i\leq A_{\delta}+2\alpha_i\right)=0 .
}\end{equation}
\end{prop}
\begin{proof}[Proof of Proposition \ref{prop642}]
\tcp{Firstly, for any $\delta>0$, note that if
$$
S_i^\delta
:=
\left\{
x\in(\alpha_i,A_\delta+2\alpha_i]:
H_i(x)-H_i(x-\alpha_i)>0
\right\},
$$
is empty, then Proposition~\ref{prop642} is immediate. Indeed, recall that, by the definition of $H_i$,
$m\mapsto H_i(x-m\alpha_i)$ is nonincreasing and that
$H_i(\alpha_i)=h_i(\alpha_i)=0$. Hence, if $S_i^\delta=\varnothing$, then
$H_i(x)=h_i(x)$ for every $x\in[\alpha_i,A_\delta+2\alpha_i]$, and therefore for every $r$,
$$
\mathbb{P}\left(H_i(FH_i)>h_i(FH_i),\ FH_i\le A_\delta+2\alpha_i\right)=0.
$$
Thus, in the remainder of the proof, we may assume that
$S_i^\delta\neq\varnothing$, and define
\begin{equation}\label{delta_di}
    \Delta_i^\delta
:=
\inf_{x\in S_i^\delta}
\left(H_i(x)-H_i(x-\alpha_i)\right)>0
\end{equation}}to represent the minimal positive increment of $H_i$ under a shift by $\alpha_i$ for $x\leq A_{\delta}+2\alpha_i$.  Note that on the compact range $x \leq A_\delta+2 \alpha_i$, the following relation holds for the increment $ H_i(x)- H_i\left(x-\alpha_i\right)$
$$
\begin{aligned}
    H_i(x)- H_i\left(x-\alpha_i\right)\in \Bigl\{\,
&\alpha_i-\Bigl(\sum_{\ell=1}^{m}\alpha_{j_\ell}-\sum_{\ell=1}^{m'}\alpha_{k_\ell}\Bigr)
+\Bigl(\sum_{\ell=1}^{m} H_i(\alpha_{j_\ell})-\sum_{\ell=1}^{m'} H_i(\alpha_{k_\ell})\Bigr)
\ :\
j_\ell,k_\ell\ge i,\\&
 \sum_{\ell=1}^{m}\alpha_{j_\ell}\le A_\delta+2\alpha_i,\ \ 
  \sum_{\ell=1}^{m'}\alpha_{k_\ell}\le A_\delta+\alpha_i
\Bigr\}.
\end{aligned}
$$
\tcgr{Because the sums $\sum_{\ell=1}^m \alpha_{j_{\ell}}$ and $\sum_{\ell=1}^{m^{\prime}} \alpha_{k_{\ell}}$ are upper bounded, respectively, by $A_\delta+2 \alpha_i$ and $A_\delta+\alpha_i$, each admissible size can occur only finitely many times. In addition, since $x\le A_\delta+2\alpha_i$ and the size sequence $0<\alpha_1<\alpha_2<\cdots$ is either
finite or goes to infinity, only finitely many indices $j\ge i$ can satisfy $\alpha_j\le A_\delta+2\alpha_i$.  Hence there are only finitely many admissible pairs of finite index sequences $\left(j_1, \ldots, j_m\right)$ and $\left(k_1, \ldots, k_{m^{\prime}}\right)$, and therefore the set of possible values of $H_i(x)-H_i\left(x-\alpha_i\right)$ is finite.}  Thus
$\Delta_i^\delta>0$. Moreover, by \eqref{eq:Hhi}, if $H_i(x)>h_i(x)$, then there exists some $m\in\mathbb{N}_+$ such that
$$
H_i\big(x-(m-1)\alpha_i\big) \;>\; H_i\big(x-m\alpha_i\big).
$$
Therefore, for any $x\le A_\delta+2\alpha_i$, if $H_i(x)>h_i(x)$, then
\begin{equation}\label{eqn23}
    H_i(x)-h_i(x)=\sum_{m=1}^{\left\lfloor x / \alpha_i\right\rfloor}\left(H_i\left(x-(m-1) \alpha_i\right)-H_i\left(x-m \alpha_i\right)\right) \geq \Delta_i^{\delta}.
\end{equation}

Recall that $\mathrm{TF}(r(y+\epsilon))$ decreases when $\tcp{j_i^*<\infty}$ and the hole $[v_{j_i^*},u_{j_i^*+1})$ satisfies 1.-3. below:
\begin{enumerate}
\item $\tcp{\alpha_i< FH_i\leq A_{\delta}+2\alpha_i}$.
\item $\tcp{H_i(FH_i)-H_i\left(FH_i-\alpha_i\right)>0}$.
\item An item of size $\alpha_i$ arrives into the system and occupies the hole $[v_{j_i^*},u_{j_i^*+1})$.
\end{enumerate}
By \eqref{delta_di}, when a scenario satisfying 1.-3. occurs, $\mathrm{TF}(r(y+\epsilon))$ decreases by at least $\Delta_i^{\delta}$. Hence, the rate at which this scenario decreases $\mathrm{TF}(r(y+\epsilon))$ is at least
$$\tcp{
 \Delta^{\delta}_i p_i r\cdot I\left(H_i(FH_i)-H_i\left(FH_i-\alpha_i\right)>0,\, \alpha_i< FH_i\leq A_{\delta}+2\alpha_i\right)+\widetilde{o}(r).
}$$
By Proposition \ref{prop636}, $\mathrm{TF}(r(y+\epsilon))$ can increase with rate at most $\Lambda^{\mathrm{TF}}$. We thus have
\begin{equation}\tcp{\label{fh1}
    \mathcal{A}(\mathrm{TF}(r(y+\epsilon))) \leq \Lambda^{\mathrm{TF}}-  \Delta^{\delta}_i p_i r\cdot I\left(H_i(FH_i)-H_i\left(FH_i-\alpha_i\right)>0,\, \alpha_i< FH_i\leq A_{\delta}+2\alpha_i\right)+\widetilde{o}(r).
}\end{equation}
Since $\lim_{r\rightarrow\infty}\mathbb{E}\left(\Lambda^{\mathrm{TF}}(\infty)\right)/{r}=0$, inequality \eqref{fh1} implies that
\begin{equation}\tcp{\label{fh2}
    \limsup _{r \rightarrow \infty} \mathbb{P}\left(H_i(FH_i)-H_i\left(FH_i-\alpha_i\right)>0, \, \alpha_i< FH_i\leq A_{\delta}+2\alpha_i\right)=0.
}\end{equation}
{Since type-$i$ item arrivals are Poisson with rate $p_i r$ and item service times are independent and exponential with unit mean,  and because the number of items completely contained in $[0,r(y+\epsilon))$ is at most $r(y+\epsilon)/\alpha_1$, the probability that no such item departs during a window of length $1/r$ is at least $e^{- (y+\epsilon)/\alpha_1}$. 
By the Poisson law, \tcorange{for a window of length $1/r$,} }$$
\mathbb{P}\left(\tcorange{\text{exactly one arrival in total, of type }i}\right)
=(p_i r)(1/r)e^{-p_i r(1/r)}\cdot e^{-(r-p_i r)(1/r)} = p_i e^{-1}.   
$$
\tcorange{Conditional on the state at time $0$, consider the event
that no item initially completely contained in
$[0,r(y+\epsilon))$ departs during $(0,1/r]$, that exactly one
arrival occurs in total during this interval and is of type $i$,
and that this item remains in the system at time $1/r$.
Given its arrival time, its probability of remaining until time
$1/r$ is at least $e^{-1/r}$. By independence of the arrival
process and the item lifetimes, for $r\geq1$, the conditional
probability of this event is at least}
\begin{equation}\label{constantc}\tcp{
p_i e^{-1}\cdot e^{(- (y+\epsilon)/\alpha_1)}\tcorange{\cdot e^{-1/r}}  \geq p_i e^{-\left( (y+\epsilon) / \alpha_1+2\right)}=c_{i,y,\epsilon}.}
\end{equation}
\tcorange{In all subsequent arguments over a time interval
of length $1/r$ with exactly one arrival, we also require
the arriving item to remain in the system until the end
of the interval.}
 \tcp{Recall the event $B^{\min}_{i,y,\epsilon}$ in \eqref{eq:bmin}}. For any $n\geq 2$, whenever the event
$$
\tcp{B^{\min}_{i,y,\epsilon}}\cap\left\{H_i\left(FH_i-(n-1) \alpha_i\right)>H_i\left(FH_i-n \alpha_i\right), \, n \alpha_i<FH_i \leq A_{\delta}+2\alpha_i\right\}
$$\tcp{
occurs at time $0$, then at time $1/r$, with probability at least $c_{i,y,\epsilon}$, the event
}$$
 \tcp{B^{\min}_{i,y,\epsilon} }\cap\left\{H_i\left(FH_i-(n-2) \alpha_i\right)>H_i\left(FH_i-(n-1) \alpha_i\right), \, (n-1)\alpha_i< FH_i \leq A_{\delta}+2\alpha_i\right\}
$$
will occur. Thus we have

\begin{equation*}\tcp{
\begin{aligned}
     &\mathbb{P}\biggr(H_i\left(FH_i-(n-2) \alpha_i\right)>H_i\left(FH_i-(n-1) \alpha_i\right), \, (n-1) \alpha_i<FH_i \leq A_{\delta}+2\alpha_i\biggr)  \\ \geq & c_{i,y,\epsilon}\cdot\mathbb{P}\biggr(H_i\left(FH_i-(n-1) \alpha_i\right)>H_i\left(FH_i-n \alpha_i\right), \, n \alpha_i<FH_i \leq A_{\delta}+2\alpha_i\biggr)-o(1).
\end{aligned}
}\end{equation*}

This implies that for any $n\in \mathbb{N}_+$,
\begin{equation*}\tcp{
\begin{aligned}
        &\limsup _{r \rightarrow \infty} \mathbb{P}\biggr(H_i\left(FH_i-(n-1) \alpha_i\right)>H_i\left(FH_i-n \alpha_i\right), \, n \alpha_i<FH_i \leq A_{\delta}+2\alpha_i\biggr)\\
        \leq & \left(\frac{1}{c_{i,y,\epsilon}}\right)^{(n-1)} \limsup _{r \rightarrow \infty} \mathbb{P}\biggr(H_i(FH_i)-H_i\left(FH_i-\alpha_i\right)>0, \, \alpha_i< FH_i\leq A_{\delta}+2\alpha_i\biggr).
\end{aligned}
}\end{equation*}
Since
$$
\begin{aligned}
    &\{H_i(FH_i)>h_i(FH_i),  FH_i\leq A_{\delta}+2\alpha_i\}\\\subseteq &\bigcup_{n=1}^{\left\lfloor{A_{\delta}}/{\alpha_i}\right\rfloor+2}\left\{H_i\left(FH_i-(n-1) \alpha_i\right)>H_i\left(FH_i-n \alpha_i\right), \ n \alpha_i<FH_i \leq A_{\delta}+2\alpha_i\right\},
\end{aligned}
$$
invoking  \eqref{fh2} yields
\begin{equation}\tcp{\label{fhn}
\begin{aligned}
      & \limsup _{r \rightarrow \infty} \mathbb{P}\biggr(H_i(FH_i)>h_i(FH_i),\, FH_i\leq A_{\delta}+2\alpha_i\biggr)\\
      \leq & \limsup _{r \rightarrow \infty} \sum_{n=1}^{\left\lfloor\frac{A_{\delta}}{\alpha_i}\right\rfloor+2}\mathbb{P}\biggr(H_i\left(FH_i-(n-1) \alpha_i\right)>H_i\left(FH_i-n \alpha_i\right), \, n \alpha_i<FH_i \leq A_{\delta}+2\alpha_i\biggr)\\
      \leq &  \left(\left\lfloor\frac{A_{\delta}}{\alpha_i}\right\rfloor+2\right) \left(\frac{1}{c_{i,y,\epsilon}}\right)^{\left(\left\lfloor\frac{A_{\delta}}{\alpha_i}\right\rfloor+2\right)} \limsup _{r \rightarrow \infty} \mathbb{P}\biggr(H_i(FH_i)-H_i\left(FH_i-\alpha_i\right)>0, \, \alpha_i< FH_i\leq A_{\delta}+2\alpha_i\biggr)\\
      =&0.
\end{aligned}
}\end{equation}
\end{proof}
Recall that any hole falls into one of two categories: (i) holes of length $x$ with $H_i(x)>h_i(x)$; and (ii) holes of length $x$ with $H_i(x)=h_i(x)$. Proposition \ref{prop642} implies that, within $[r\beta_{i-1},\, r(y+\epsilon))$, the  probability that the hole $[v_{j_i^*},u_{j_i^*+1})$ is of category (i) vanishes as $r\to\infty$. \tcp{Consequently, invoking Proposition \ref{prop642}, we can control the number of such category (i) packing configurations, which we do via the index set $\mathcal{R}_i$ defined below.}

\tcp{\tcgr{Note that for any index $j$, if $H_i(e_j)>h_i(e_j)$ and the hole $[v_j,u_{j+1})$ on the right has length $g_j<\alpha_i$, then
$$
e_j+g_j=(v_j-u_j)+(u_{j+1}-v_j)=u_{j+1}-u_j.
$$
Since $0\le g_j<\alpha_i$, Corollary~\ref{hyx2} gives
$$
H_i(e_j+g_j)=H_i(e_j)+g_j,
$$
and the definition of $h_i$ \tcgr{in \eqref{def:hi}} gives
$$
h_i(e_j+g_j)\le h_i(e_j)+g_j.
$$
Therefore $H_i(e_j+g_j)>h_i(e_j+g_j)$.}}
Moreover, by Proposition~\ref{npr1}, in the hydrodynamic scaling, the number of holes of length at least $\alpha_i$ in $[r\beta_{i-1},\,r(y+\epsilon))$ vanishes as $r\to\infty$. 
Hence, in order to study the number of type-$k$ items with $H_i(\alpha_k)>h_i(\alpha_k)$ in $[r\beta_{i-1},\,r(y+\epsilon))$, it will suffice to count the indices $j$ with $H_i\left(e_j+g_j\right)>h_i\left(e_j+g_j\right)$ whose left hole $[v_{j-1},u_j)$ has length less than $\alpha_i$. This key observation motivates the following definition:
\begin{equation}
\begin{aligned}
        \mathcal{R}_i =& \{\, j : \tcorange{j>j_{\min}(y+\epsilon)},\ v_{j+1}\tcorange{\leq} r(y+\epsilon),\ g_{j-1}<\alpha_i,\\
        & e_j \geq \alpha_i,\; g_j < \alpha_i,\ H_i\left(e_j+g_j\right) > h_i\left(e_j+g_j\right) \}.
\end{aligned}
\end{equation}
The above index set collects the indices $j$ for which the item in $[u_j,v_j)$ is of size at least $\alpha_i$, both adjacent holes have length less than $\alpha_i$, and the interval $[u_j,u_{j+1})$ satisfies $H_i\left(e_j+g_j\right) > h_i\left(e_j+g_j\right)$. Let $|\mathcal{R}_i|$ be the cardinality of the index set $\mathcal{R}_i$, that is, the number of indices $j$ satisfying the above conditions. This definition and the preceding observation lead to Proposition \ref{prop643} below.

\begin{prop}\label{prop643}
     For any $y\in[\beta_{i-1},\beta_i)$ and $0<\epsilon< (\beta_{i}-y)/2 $, under the inductive hypothesis at $y$, the random variable $|\mathcal{R}_i|$  satisfies 
    \begin{equation*}
        \limsup_{r\rightarrow\infty}\frac{\mathbb{E}(|\mathcal{R}_i|)}{r}=0.
    \end{equation*}
\end{prop}
\begin{proof}[Proof of Proposition \ref{prop643}]
For any $\delta>0$, we define two auxiliary index sets related to $\mathcal{R}_i$
$$
\begin{aligned}
       \mathcal{R}^{\delta}_i = & \{\, j : \tcorange{j>j_{\min}(y+\epsilon)},\ v_{j+1}\tcorange{\leq} r(y+\epsilon),\ g_{j-1}<\alpha_i,\\
       & \alpha_i\leq e_j \leq A_{\delta},\; g_j < \alpha_i,\ H_i\left(e_j+g_j\right) > h_i\left(e_j+g_j\right) \}, 
\end{aligned}
$$
and
$$
\begin{aligned}
    \widetilde{\mathcal{R}}^{\delta}_i =& \{\, j : \; \;  v_{j+1}\tcorange{\leq} r(y+\epsilon), \ \alpha_i \leq g_j<A_{\delta}+2\alpha_i,\ H_i(g_j) > h_i(g_j) \}.
\end{aligned}
$$
The index set $\widetilde{\mathcal{R}}^{\delta}_i$ collects the indices of holes whose length is at least $\alpha_i$ and less than $A_{\delta}+2\alpha_i$ and for which the corresponding value of $H_i$ is strictly greater than $h_i$.

\tcp{As in the proof of Proposition~\ref{prop642}, 
 if $S_i^\delta=\varnothing$, then $H_i(x)=h_i(x)$ for all
$x\in[\alpha_i,A_\delta+2\alpha_i]$, which implies
$$
\mathcal R_i^\delta=\widetilde{\mathcal R}_i^\delta=\varnothing.
$$
Hence, we may assume that
$S_i^\delta\neq\varnothing$ so that $\Delta_i^\delta$ is well-defined.}

Here, we recall that $j_i^*$ denotes the index of the first hole of length at least $\alpha_i$ with adjacent items lying entirely in $[r\beta_{i-1},\,r(y+\epsilon))$. We also recall that $FH_i$ is the length of this hole. We have
$$
\{j_i^* \in \widetilde{\mathcal{R}}^{\delta}_i\}\subseteq \{H_i(FH_i)>h_i(FH_i),\;  FH_i \leq A_{\delta}+2\alpha_i\}.
$$
By Proposition~\ref{prop642},
\begin{equation}\label{jri}
  \limsup_{r \to \infty}\,\mathbb{P}\!\left(j_i^* \in \widetilde{\mathcal{R}}^{\delta}_i\right)
\leq \limsup_{r \to \infty}\,\mathbb{P}\!\left(H_i(FH_i)>h_i(FH_i),\; FH_i \leq A_{\delta}+2\alpha_i\right)
= 0.  
\end{equation}
Note that event $\{j_i^* \in \widetilde{\mathcal{R}}^{\delta}_i\}$ \tcp{contains} the event 
$$
\left\{|\widetilde{\mathcal{R}}_i^\delta|>0, \text { and there is no hole of size at least } \alpha_i \text { in }\left[0, v_{\min (\widetilde{\mathcal{R}}_i^\delta)}\right)\right\},
$$ where $\min (\widetilde{\mathcal{R}}_i^\delta)$ denotes the smallest index in $\widetilde{\mathcal{R}}_i^\delta$ and $v_{\min (\widetilde{\mathcal{R}}_i^\delta)}$ denotes the right endpoint of the item with index $\min (\widetilde{\mathcal{R}}_i^\delta)$.  \tcp{Indeed, on $\{|\widetilde{\mathcal{R}}_i^\delta|>0\}$ the index $\min(\widetilde{\mathcal{R}}_i^\delta)$ identifies the leftmost hole in $\widetilde{\mathcal{R}}_i^\delta$; the additional requirement that there is no hole of size at least $\alpha_i$ in $[0, v_{\min(\widetilde{\mathcal{R}}_i^\delta)})$ forces this hole to be the first hole of length at least $\alpha_i$, and hence its index \tcgr{is} $j_i^*$.}

 \tcp{We proceed to define $\widetilde D^{\mathcal{R}}_{i,\delta}$ as the total number of $i$-items that can potentially (completely) fit into the available empty space within $[0,\,v_{\min(\widetilde{\mathcal{R}}_i^\delta)})$. If $\left|\widetilde{\mathcal{R}}_i^\delta\right|=0$, we set $\widetilde{D}_{i, \delta}^{\mathcal{R}}=D_i$. On the event $\{|\widetilde{\mathcal{R}}_i^\delta|>0\}$, the condition $\widetilde D^{\mathcal{R}}_{i,\delta}=0$ is the requirement that there is no hole of size at least $\alpha_i$ in $[0, v_{\min(\widetilde{\mathcal{R}}_i^\delta)})$, and thus the event $\left\{|\widetilde{\mathcal{R}}_i^\delta|>0,\,\widetilde D^{\mathcal{R}}_{i,\delta}=0\right\}$ is a subset of $\{j_i^* \in \widetilde{\mathcal{R}}^{\delta}_i\}$.} Therefore, by \eqref{jri},
\begin{equation}\label{jri2}
    \tcp{\lim_{r \to \infty}\mathbb{P}\bigl(|\widetilde{\mathcal{R}}_i^\delta|>0,\,\widetilde D^{\mathcal{R}}_{i,\delta}=0\bigr) \leq \lim_{r \to \infty}\mathbb{P} \!\left(j_i^* \in \widetilde{\mathcal{R}}^{\delta}_i\right)=0.}
\end{equation}
\tcp{Recall (as in the proof of Proposition~\ref{prop642}) that, conditional on the state at time $0$, arrivals in $(0,1/r]$ are independent of the remaining service times of the items present at time $0$. In particular, over a time window of length $1/r$, the event that there is no departure of items lying entirely in $[0,r(y+\epsilon))$ and that exactly one arrival occurs in total and its type is $i$ has probability at least $c_{i,y,\epsilon}$ as defined in \eqref{constantc}. On the above event and on $\{|\widetilde{\mathcal{R}}_i^\delta|>0,\,\tcblue{1\leq}\widetilde D^{\mathcal{R}}_{i,\delta}\tcblue{\leq m}\}$}, the single \tcp{ type-$i$ arrival is placed within $[0,\,v_{\min(\widetilde{\mathcal{R}}_i^\delta)})$ \tcblue{and reduces the total number of $i$-items that can potentially (completely) fit into the available empty space within this interval by one. A hole further to the left of the hole with index $\min(\widetilde{\mathcal{R}}_i^\delta)$ may become the leftmost hole in $\widetilde{\mathcal{R}}_i^\delta$ at time $1/r$; hence $\widetilde D^{\mathcal{R}}_{i,\delta}$ may decrease by more than one, and in particular $\widetilde D^{\mathcal{R}}_{i,\delta}\leq m-1$.}} Hence, for any $m\in \mathbb{N}_+$,
\begin{equation}\label{eq:ad1}
\tcblue{\begin{aligned}
&\mathbb{P}\bigl(|\widetilde{\mathcal{R}}_i^\delta|>0,\,\widetilde D^{\mathcal{R}}_{i,\delta}\leq m\bigr)\\
&\qquad\leq \mathbb{P}\bigl(|\widetilde{\mathcal{R}}_i^\delta|>0,\,\widetilde D^{\mathcal{R}}_{i,\delta}=0\bigr)
+\frac{1}{c_{i,y,\epsilon}}\mathbb{P}\bigl(|\widetilde{\mathcal{R}}_i^\delta|>0,\,\widetilde D^{\mathcal{R}}_{i,\delta}\leq m-1\bigr).
\end{aligned}}
\end{equation}
\tcblue{Applying \eqref{eq:ad1} successively for $m=1,\ldots,N$ and using} \eqref{jri2}, for any fixed $N\in\mathbb{Z}_{+}$, $$ \lim_{r\to\infty}\mathbb{P}\bigl(|\widetilde{\mathcal{R}}_i^\delta|>0,\ \widetilde D^{\mathcal{R}}_{i,\delta}\le N\bigr)=0. $$ 
Hence, by Proposition \ref{npr1} and by the fact that $\tcp{ \tcor{\widetilde D^{\mathcal{R}}_{i,\delta}}\leq D_i}$, for any $N$ and $\widetilde\delta>0$, \tcgr{with $B_{\widetilde\delta,i}$ and $z_{y,\epsilon}$ defined in \eqref{eq:B-di} and \eqref{eq:z-ye}, respectively, we have} \begin{equation}\tcp{\label{eq:34}
    \begin{aligned} \limsup_{r\to\infty}\mathbb{P}\bigl(|\widetilde{\mathcal{R}}_i^\delta|>0\bigr) \leq & \limsup_{r\to\infty}\mathbb{P}\bigl(|\widetilde{\mathcal{R}}_i^\delta|>0,\ \widetilde D^{\mathcal{R}}_{i,\delta}\le N\bigr) +\limsup_{r\to\infty}\mathbb{P}\bigl(\widetilde D^{\mathcal{R}}_{i,\delta}> N\bigr)\\
    \leq & \ 0+ \limsup_{r\to\infty}\mathbb{P}\bigl( D_i\tcblue{>}N\bigr) \leq \frac{p_i+z_{y,\epsilon}}{p_i-z_{y,\epsilon}}\frac{B_{\widetilde\delta,i}+1}{N}+\frac{6\widetilde\delta}{\alpha_1(p_i-z_{y,\epsilon})}, \end{aligned}
}\end{equation} 
yielding
\begin{equation}\label{tr} \lim_{r\to\infty}\mathbb{P}\bigl(|\widetilde{\mathcal{R}}_i^\delta|>0\bigr)=0.
\end{equation}

\tcp{Equation \eqref{tr} indicates that the probability of  
$|\widetilde{\mathcal{R}}^{\delta}_i| > 0$ 
vanishes as $r \to \infty$. This will help us  control \tcp{the rate of increase of }the Lyapunov functional  
$
\mathrm{TF}(r(y+\epsilon)) - |\widetilde{\mathcal{R}}^{\delta}_i|.$}

We next examine the mechanism by which $\mathrm{TF}(r(y+\epsilon))-\big|\widetilde{\mathcal{R}}^{\delta}_i\big|$ decreases. For a given index $j$, we shall study whether or not the departure of the item $[u_j,v_j)$ causes a decrease of $\mathrm{TF}(r(y+\epsilon))-\big|\widetilde{\mathcal{R}}^{\delta}_i\big|$. To do so, we distinguish between two cases: (i) $H_i(\tcgr{g_{j-1}+e_{j}+g_{j}})=h_i(\tcgr{g_{j-1}+e_{j}+g_{j}})$; and (ii) $H_i(\tcgr{g_{j-1}+e_{j}+g_{j}})>h_i(\tcgr{g_{j-1}+e_{j}+g_{j}})$.  Recall that $\tcp{j_{\min}(y+\epsilon)}$ is the index of the first item fully contained in $[r \beta_{i-1},r(y+\epsilon))$. For any $j\in\mathcal{R}^{\delta}_i$ with $j\neq \tcp{j_{\min}(y+\epsilon)}$, if the item occupying $[u_{j-1},v_{j-1})$ has size at least $\alpha_i$, then when the item occupying the interval $[u_j,v_j)$ departs from the system the corresponding change in $\mathrm{TF}(r(y+\epsilon))$ is
\begin{equation}\label{eq24}
\begin{aligned}
    \Delta\mathrm{TF}=&H_i\left(\tcgr{g_{j-1}+e_{j}+g_{j}}\right)-\left(\tcgr{g_{j-1}}+H_i\left(\tcgr{e_j+g_j}\right)\right)\\
    \leq& I\left(H_i\left(\tcgr{g_{j-1}+e_{j}+g_{j}}\right)=h_i\left(\tcgr{g_{j-1}+e_{j}+g_{j}}\right)\right) \cdot\left(h_i\left(\tcgr{g_{j-1}+e_{j}+g_{j}}\right)-\left(\tcgr{g_{j-1}}+H_i\left(\tcgr{e_j+g_j}\right)\right)\right)\\
    \leq & I\left(H_i\left(\tcgr{g_{j-1}+e_{j}+g_{j}}\right)=h_i\left(\tcgr{g_{j-1}+e_{j}+g_{j}}\right)\right) \cdot\left(h_i\left(\tcgr{e_j+g_j}\right)-H_i\left(\tcgr{e_j+g_j}\right)\right)\\
    \leq & -\Delta_i^{\delta}\cdot I\left(H_i\left(\tcgr{g_{j-1}+e_{j}+g_{j}}\right)=h_i\left(\tcgr{g_{j-1}+e_{j}+g_{j}}\right)\right).
\end{aligned}
\end{equation}
\tcp{Some remarks about the above derivation in \eqref{eq24} are now in order. In the first inequality, if
$H_i(\tcgr{g_{j-1}+e_{j}+g_{j}})=h_i(\tcgr{g_{j-1}+e_{j}+g_{j}})$, then the inequality becomes an equality. If instead
$H_i(\tcgr{g_{j-1}+e_{j}+g_{j}})>h_i(\tcgr{g_{j-1}+e_{j}+g_{j}})$, then the indicator vanishes and, by Proposition~\ref{prop632},
$\Delta \mathrm{TF}\le 0$ . In the second inequality,
\tcgr{by the definition of $h_i$ in \eqref{def:hi}}, for any $a\ge0$ and 
\tcor{$b\geq 0$}, $
h_i(a+b)-b\le h_i(a).$ Finally, the third inequality follows from the definition of $\mathcal{R}_i^\delta$. Since $j \in \mathcal{R}_i^\delta$, we have that $H_i\left(\tcgr{e_j+g_j}\right) > h_i\left(\tcgr{e_j+g_j}\right)$, which, by \eqref{eqn23}, guarantees that $H_i\left(\tcgr{e_j+g_j}\right) - h_i\left(\tcgr{e_j+g_j}\right) \ge \Delta_i^\delta$.} 

We now fix $j\in\mathcal{R}^{\delta}_i$ with $j\neq\tcp{j_{\min}(y+\epsilon)}$ and suppose that the item occupying $[u_{j-1},v_{j-1})$ has length at least $\alpha_i$. If
$$
H_i\big(\tcgr{g_{j-1}+e_{j}+g_{j}}\big) > h_i\big(\tcgr{g_{j-1}+e_{j}+g_{j}}\big),
$$
when the item $[u_j, v_j)$ departs, the index $j-1$ will belong to $\widetilde{\mathcal{R}}^{\delta}_i$ and thus $\big|\widetilde{\mathcal{R}}^{\delta}_i\big|$ will increase by one.  
\tcp{Recall that, by Proposition~\ref{prop632}, $\mathrm{TF}$ remains non-increasing upon such a departure}. Therefore, in this case, the Lyapunov functional 
$$
\mathrm{TF}(r(y+\epsilon)) - \big|\widetilde{\mathcal{R}}^{\delta}_i\big|
$$
decreases by at least $1$.
On the other hand, if 
$$
H_i\big(\tcgr{g_{j-1}+e_{j}+g_{j}}\big) = h_i\big(\tcgr{g_{j-1}+e_{j}+g_{j}}\big),
$$
then, as established by \eqref{eq24}, $\mathrm{TF}(r(y+\epsilon))$ decreases by at least $\Delta^{\delta}_i$. Thus, the Lyapunov functional will decrease by at least $\Delta^{\delta}_i$.  
The indicator appearing in inequality \eqref{eq24} represents the event 
$H_i(\tcgr{g_{j-1}+e_{j}+g_{j}}) = h_i(\tcgr{g_{j-1}+e_{j}+g_{j}})$,  
which implies that either $\mathrm{TF}(r(y+\epsilon))$ decreases or $\big|\widetilde{\mathcal{R}}^{\delta}_i\big|$ increases.  
Thus, the combined functional $\mathrm{TF}(r(y+\epsilon)) - \big|\widetilde{\mathcal{R}}^{\delta}_i\big|$ decreases by at least $\min\!\left(\Delta^{\delta}_i, 1\right)$ upon the departure of the item $[u_j,v_j)$.

\tcgr{The number of indices $j\in \mathcal R_i^\delta$ such that either
$j=j_{\min}(y+\epsilon)$ or the left neighboring item $[u_{j-1},v_{j-1})$
has size less than $\alpha_i$ is at most}
$$
\tcp{
1+\sum_{\ell<i}\left(F_{\ell}(r(y+\epsilon))-F_{\ell}(r\beta_{i-1})\right)
=\widetilde{o}(r).
}
$$
\tcgr{Hence, the number of indices $j\in \mathcal R_i^\delta$ such that
$j\neq j_{\min}(y+\epsilon)$ and the left neighboring item $[u_{j-1},v_{j-1})$
has size at least $\alpha_i$ is}
$
\tcp{
|\mathcal R_i^\delta|+\widetilde{o}(r).
}
$
Therefore, \tcp{the rate of decrease of
$\mathrm{TF}(r(y+\epsilon))-|\widetilde{\mathcal R}_i^\delta|$
is at least}
$$
\tcp{
\min\left(\Delta_i^\delta,1\right)\cdot|\mathcal R_i^\delta|
+\widetilde{o}(r).
}
$$

On the other hand, as noted in Proposition \ref{prop636}, \tcp{the rate of increase of }$\mathrm{TF}(r(y+\epsilon))$ can be upper bounded by $\Lambda^{\mathrm{TF}}.$ There are two mechanisms by which $|\widetilde{\mathcal{R}}^{\delta}_i|$ can decrease. The first occurs when a new item arrives to occupy the hole in a packing configuration belonging to $\widetilde{\mathcal{R}}^{\delta}_i$. The rate of \tcp{the decrease of $|\widetilde{\mathcal{R}}^{\delta}_i|$} caused by such arrivals can be upper bounded by $r\cdot I\left(|\widetilde{\mathcal{R}}^{\delta}_i|>0\right)$. Secondly, \tcp{a decrease in $|\widetilde{\mathcal{R}}^{\delta}_i|$} may occur when one of the items adjacent to the hole in a packing configuration from $\widetilde{\mathcal{R}}^{\delta}_i$ departs from the system, and the rate of \tcp{the decrease of $|\widetilde{\mathcal{R}}^{\delta}_i|$} caused by such departures can be upper bounded by $2 |\widetilde{\mathcal{R}}^{\delta}_i|$. Thus, the rate of decrease of $|\widetilde{\mathcal{R}}^{\delta}_i|$ can be upper bounded by $$r \cdot I(|\widetilde{\mathcal{R}}^{\delta}_i|>0)+2 |\widetilde{\mathcal{R}}^{\delta}_i|.$$ Combining all of the above dynamics of the Lyapunov functional $\mathrm{TF}(r(y+\epsilon))- |\widetilde{\mathcal{R}}^{\delta}_i|$, we write
\begin{equation}\label{ad2}
\begin{aligned}
    \mathcal{A}\left(\mathrm{TF}(r(y+\epsilon))-| \widetilde{\mathcal{R}}^{\delta}_i|\right) \leq &-\min \left(\Delta^{\delta}_i, 1\right) \cdot  | \mathcal{R}^{\delta}_i|+\Lambda^{\mathrm{TF}}+2|\widetilde{\mathcal{R}}^{\delta}_i|+r \cdot I(|\widetilde{\mathcal{R}}^{\delta}_i| > 0)+\widetilde{o}(r).
\end{aligned}
\end{equation}
Note that $|\widetilde{\mathcal{R}}^{\delta}_i|\leq D_i$. Combining \eqref{ad2} with \eqref{tr} yields
       \begin{equation*}
        \limsup_{r\rightarrow\infty}\frac{\mathbb{E}(|{\mathcal{R}}^{\delta}_i|)}{r}=0.
    \end{equation*}
    \tcgr{It now remains to compare $\mathcal R_i^\delta$ and $ \mathcal R_i$, and then use the estimate for $|\mathcal R_i^\delta|$ to obtain an upper bound on $| \mathcal R_i|$.}
     By the definitions of $\mathcal{R}_i$ and $\mathcal{R}^{\delta}_i$, the set $\mathcal{R}^{\delta}_i$ is obtained from $\mathcal{R}_i$ by removing those indices $j$ for which the item $[u_j,v_j)$ has size greater than $A_{\delta}$. The number of such indices is upper bounded by the total number of items in the system with size exceeding $A_{\delta}$, i.e., by $\sum_{\ell:\,\alpha_\ell>A_{\delta}} F_{\ell}(\infty)$. Therefore, 
    \begin{equation*}
   \limsup_{r\rightarrow\infty}\left(\frac{\mathbb{E}(|\mathcal{R}_i|)}{r}-\frac{\mathbb{E}(|\mathcal{R}^{\delta}_i|)}{r}\right)\leq \limsup_{r\rightarrow\infty} \mathbb{E}\left( \frac{1}{r}\sum_{\ell: \alpha_\ell> A_{\delta}}F_{\ell}\left(\infty ; \infty\right)\right)\leq \frac{\delta}{\alpha_1},
    \end{equation*}
 and thus
           \begin{equation*}
        \limsup_{r\rightarrow\infty}\frac{\mathbb{E}(|\mathcal{R}_i|)}{r}=0.
    \end{equation*}
\end{proof}
\vspace{-3mm}
With Proposition \ref{prop643} in hand, and based on the preceding analysis of the relationship between type-$k$ items with $H_i(\alpha_k)>h_i(\alpha_k)$ in $[r\beta_{i-1},\,r(y+\epsilon))$ and the indices counted by $|\mathcal{R}_i|$, it is straightforward to obtain Corollary \ref{corr644} below. Corollary \ref{corr644} states that under the inductive hypothesis at $y$, the hydrodynamic-scaled  number of any type-$k$ item with $H_i(\alpha_k)>h_i(\alpha_k)$ in $[r\beta_{i-1},r(y+\epsilon))$ vanishes.
\begin{corollary}\label{corr644}
 For any $y\in[\beta_{i-1},\beta_i)$ and $0<\epsilon< (\beta_{i}-y)/2 $, under the inductive hypothesis at $y$,
for any $k>i$ satisfying
\begin{equation*}
    H_i\left(\alpha_k\right)>h_i(\alpha_k),
\end{equation*}
\begin{equation*}
    \lim _{r \rightarrow \infty} \frac{1}{r} \mathbb{E}\left(F_k(r(y+\epsilon) ; \infty)-F_k\left(r \beta_{i-1} ; \infty\right)\right)=0.
\end{equation*}
\end{corollary}
\begin{proof}\tcp{
    Fix $k>i$ such that $H_i(\alpha_k)>h_i(\alpha_k)$. If a type-$k$ item $[u_j,v_j)$ completely in $[r\beta_{i-1},r(y+\epsilon))$ is not one of the at most two boundary items and both of its adjacent holes have length strictly less than $\alpha_i$,  Corollary~\ref{hyx2} yields
$$
H_i\left(e_j+g_j\right)=H_i(\alpha_k)+g_j>h_i(\alpha_k)+g_j\geq h_i\left(e_j+g_j\right).
$$
Therefore, by the definition of $\mathcal{R}_i$, the index $j$ belongs to $\mathcal R_i$.  Combining this with the fact that the number of items that are either boundary items or adjacent to holes with length at least $\alpha_i$ is upper bounded by $2+2D_i=\widetilde{o}(r)$ yields $$F_k(r(y+\epsilon))-F_k(r\beta_{i-1}) \leq |\mathcal{R}_i| + \widetilde{o}(r).$$
Hence,
$$\lim_{r \rightarrow \infty} \frac{1}{r} \mathbb{E}\big[F_k(r(y+\epsilon) ; \infty)-F_k(r \beta_{i-1} ; \infty)\big] \leq \lim _{r \rightarrow \infty} \frac{\mathbb{E}\big(|\mathcal{R}_i|\big) }{r} = 0.$$}
\end{proof}

We are now in the position to analyze and bound the number of holes of length $x$ with $H_i(x)=h_i(x)$, as well as the number of type-$k$ items with $k>i$ satisfying $H_i(\alpha_k)=h_i(\alpha_k)$. Suppose that two holes of length less than $\alpha_i$ are separated by an item. \tcgr{If that separating item departs and, before either adjacent item departs, the first item subsequently placed in the merged interval created by that departure has the same size as the departed item, then the first-fit rule places it flush with the left endpoint of that merged interval.} \tcp{Consequently, the remaining empty interval to the right of that item is a single hole whose length equals the sum of the lengths of the two holes \tcgr{present immediately before that departure.}} Iterating this mechanism causes holes of lengths less than $\alpha_i$ to merge into \tcgr{fewer, larger} holes. \tcgr{Starting from a hole with length in $[x,\alpha_i)$, repeated departures of separating items and subsequent first-fit arrivals may merge it with neighboring holes of lengths less than $\alpha_i$; in this way, a hole of length at least $\alpha_i$ may eventually be created. A subsequent size-$\alpha_i$ arrival may then be placed in that hole, thereby possibly decreasing $\operatorname{TF}(r(y+\epsilon))$.} Motivated by this mechanism, for any $\delta>0$, $0<x<\alpha_i$, and $m\in\mathbb{N}_{+}$, we define the following two index sets.
\begin{equation}\label{eq:L-sets}
\begin{aligned}
L_x^{\delta,m} = &\big\{\, j \ \big|\ x \le \tcgr{g_j} < \alpha_i,\ \forall\, 0 < \ell \le m,\ \alpha_i\le \tcgr{e_{j+\ell}} \le A_\delta,\\
& \tcgr{g_{j+\ell}} < \alpha_i,\ v_{j+m+1}\tcorange{\leq} r(y+\epsilon),\ H_i\left(\tcgr{e_{j+\ell}+g_{j+\ell}}\right)=h_i\left(\tcgr{e_{j+\ell}+g_{j+\ell}}\right),\\
& \tcgr{g_j}+\sum_{\ell=1}^{m-1} H_i\left(\tcgr{e_{j+\ell}+g_{j+\ell}}\right) < \alpha_i,\ \tcgr{g_j}+\sum_{\ell=1}^{m} H_i\left(\tcgr{e_{j+\ell}+g_{j+\ell}}\right) \ge \alpha_i \big\}, \\[3pt]
\end{aligned}
\end{equation}
\begin{equation}\label{eq:tL-sets}
\begin{aligned}
\widetilde{L}_x^{\delta,m} = &\big\{\, j \ \big|\ x+\alpha_i \le \tcgr{g_j} < A_\delta+2\alpha_i,\ H_i\left(\tcgr{g_j}\right)=h_i\left(\tcgr{g_j}\right)\ge x,\ \forall\, 0 < \ell \le m,\\
& \alpha_i\le \tcgr{e_{j+\ell}} \le A_\delta,\ \tcgr{g_{j+\ell}} < \alpha_i,\ H_i\left(\tcgr{e_{j+\ell}+g_{j+\ell}}\right)=h_i\left(\tcgr{e_{j+\ell}+g_{j+\ell}}\right),\\
& v_{j+m+1}\tcorange{\leq} r(y+\epsilon),\ h_i\left(\tcgr{g_j}\right)+\sum_{\ell=1}^{m-1} H_i\left(\tcgr{e_{j+\ell}+g_{j+\ell}}\right) < \alpha_i,\\
& h_i\left(\tcgr{g_j}\right)+\sum_{\ell=1}^{m} H_i\left(\tcgr{e_{j+\ell}+g_{j+\ell}}\right) \ge \alpha_i \big\}.
\end{aligned}
\end{equation}
\noindent We pause to explain the definitions of the above index sets $L_x^{\delta,m}$ and $\widetilde{L}_x^{\delta,m}$. The set $L_x^{\delta,m}$ defined in \eqref{eq:L-sets} consists of indices $j$ for which the packing configuration of holes and items beginning at $j$ has the following properties:
\begin{itemize}
    \item The first hole $[v_j,u_{j+1})$ has length in $[x,\alpha_i)$.
    \item For each $1\le \ell\le m$, the item $[u_{j+\ell},v_{j+\ell})$ has size in $[\alpha_i,A_\delta]$.
    \item For each $1\le \ell\le m$, the hole $[v_{j+\ell},u_{j+\ell+1})$ has length less than $\alpha_i$.
    \item In addition, for each $1\le \ell\le m$, we require
$$
H_i\left(\tcgr{e_{j+\ell}+g_{j+\ell}}\right)=h_i\left(\tcgr{e_{j+\ell}+g_{j+\ell}}\right),
$$
together with
$$
\tcgr{g_j}+\sum_{\ell=1}^{m-1} H_i\left(\tcgr{e_{j+\ell}+g_{j+\ell}}\right)<\alpha_i,
$$
and
$$
\tcgr{g_j}+\sum_{\ell=1}^{m} H_i\left(\tcgr{e_{j+\ell}+g_{j+\ell}}\right)\ge \alpha_i.
$$
\end{itemize}
\tcor{A schematic example of a packing configuration counted by $L_x^{\delta,m}$ is shown in Figure~\ref{fig:L-config-revised}.}

With this definition, \tcp{for any $m\ge 2$}, each time the leftmost item in a packing configuration from $L_x^{\delta,m}$ departs, the resulting \tcor{packing configuration} belongs to $\widetilde{L}_x^{\delta,m-1}$. Then, when successive size-$\alpha_i$ item arrivals occupy the leftmost hole of a packing configuration in $\widetilde{L}_x^{\delta,m-1}$, the system transitions to a packing configuration in $L_x^{\delta,m-1}$. \tcor{Figure~\ref{fig:L-transition-revised} illustrates the two steps that take a configuration counted by $L_x^{\delta,m}$ first to one counted by $\widetilde{L}_x^{\delta,m-1}$ and then to one counted by $L_x^{\delta,m-1}$.} Repeatedly applying these two steps reduces the parameter $m$ by $1$ at each stage, until one reaches a packing configuration in $L_x^{\delta,1}$.

\begin{figure}[!htbp]
\centering

\begin{tikzpicture}[xscale=1.15, yscale=1.15]

\def\baseY{0}
\def\barH{0.8}

\fill[gray!30] (0,\baseY) rectangle (1.2, \baseY+\barH);
\draw (0,\baseY) rectangle (1.2, \baseY+\barH);
\node at (0.6, \baseY+\barH/2) {\small $e_j$};

\draw[dashed] (1.2,\baseY) -- (1.2, \baseY+\barH);
\draw[dashed] (2.0,\baseY) -- (2.0, \baseY+\barH);
\draw (1.2,\baseY) -- (2.0,\baseY);
\draw[decorate, decoration={brace, mirror, amplitude=5pt}] (1.2, \baseY-0.15) -- (2.0, \baseY-0.15)
  node[midway, below=6pt] {\small $g_j$};

\fill[blue!20] (2.0,\baseY) rectangle (3.6, \baseY+\barH);
\draw (2.0,\baseY) rectangle (3.6, \baseY+\barH);
\node at (2.8, \baseY+\barH/2) {\small $e_{j+1}$};

\draw[dashed] (3.6,\baseY) -- (3.6, \baseY+\barH);
\draw[dashed] (4.0,\baseY) -- (4.0, \baseY+\barH);
\draw (3.6,\baseY) -- (4.0,\baseY);
\draw[decorate, decoration={brace, mirror, amplitude=5pt}] (3.6, \baseY-0.15) -- (4.0, \baseY-0.15)
  node[midway, below=6pt] {\small $g_{j+1}$};

\fill[blue!20] (4.0,\baseY) rectangle (5.4, \baseY+\barH);
\draw (4.0,\baseY) rectangle (5.4, \baseY+\barH);
\node at (4.7, \baseY+\barH/2) {\small $e_{j+2}$};

\draw[dashed] (5.4,\baseY) -- (5.4, \baseY+\barH);
\draw[dashed] (5.7,\baseY) -- (5.7, \baseY+\barH);
\draw (5.4,\baseY) -- (5.7,\baseY);
\draw[decorate, decoration={brace, mirror, amplitude=5pt}] (5.4, \baseY-0.15) -- (5.7, \baseY-0.15)
  node[midway, below=6pt] {\small $g_{j+2}$};

\fill[blue!20] (5.7,\baseY) rectangle (7.3, \baseY+\barH);
\draw (5.7,\baseY) rectangle (7.3, \baseY+\barH);
\node at (6.5, \baseY+\barH/2) {\small $e_{j+3}$};

\draw[dashed] (7.3,\baseY) -- (7.3, \baseY+\barH);
\draw[dashed] (7.8,\baseY) -- (7.8, \baseY+\barH);
\draw (7.3,\baseY) -- (7.8,\baseY);
\draw[decorate, decoration={brace, mirror, amplitude=5pt}] (7.3, \baseY-0.15) -- (7.8, \baseY-0.15)
  node[midway, below=6pt] {\small $g_{j+3}$};

\draw[<->] (1.2, \baseY+\barH+0.25) -- (2.0, \baseY+\barH+0.25);
\node[above] at (1.6, \baseY+\barH+0.25) {\small $g_j\in[x,\alpha_i)$};

\draw[decorate, decoration={brace, amplitude=6pt}] (2.0, \baseY+\barH+0.65) -- (7.3, \baseY+\barH+0.65)
  node[midway, above=7pt] {\small for $1\le \ell\le m$: $\alpha_i\le e_{j+\ell}\le A_\delta$ and $g_{j+\ell}<\alpha_i$};

\end{tikzpicture}
\caption{\tcgr{Schematic representation of a packing configuration counted by $L_x^{\delta,m}$ in \eqref{eq:L-sets}, shown here for $m=3$. The leftmost hole satisfies $g_j\in[x,\alpha_i)$. For each $1\le \ell\le m$, the corresponding item-hole pair satisfies $\alpha_i\le e_{j+\ell}\le A_\delta$, $g_{j+\ell}<\alpha_i$, and $H_i(e_{j+\ell}+g_{j+\ell})=h_i(e_{j+\ell}+g_{j+\ell})$. In addition, the packing configuration shown also satisfies
$g_j+\sum_{\ell=1}^{m-1}H_i(e_{j+\ell}+g_{j+\ell})<\alpha_i$
and
$g_j+\sum_{\ell=1}^{m}H_i(e_{j+\ell}+g_{j+\ell})\ge \alpha_i$.
The figure is schematic and not to scale.}}
\label{fig:L-config-revised}

\vspace{1em}      %

\begin{tikzpicture}[xscale=1.1, yscale=1.25]

\def\rowA{5.5}
\def\barH{0.7}

\node[anchor=east] at (-0.3, \rowA+\barH/2) {\small $L_x^{\delta,3}$};

\fill[gray!30] (0,\rowA) rectangle (1.0, \rowA+\barH);
\draw (0,\rowA) rectangle (1.0, \rowA+\barH);
\node at (0.5, \rowA+\barH/2) {\tiny $e_j$};

\draw (1.0,\rowA) -- (1.7,\rowA);
\draw[dashed] (1.0,\rowA) -- (1.0,\rowA+\barH);
\draw[dashed] (1.7,\rowA) -- (1.7,\rowA+\barH);
\node[below] at (1.35, \rowA-0.05) {\tiny $g_j$};

\fill[red!25] (1.7,\rowA) rectangle (3.1, \rowA+\barH);
\draw (1.7,\rowA) rectangle (3.1, \rowA+\barH);
\node at (2.4, \rowA+\barH/2) {\tiny $e_{j+1}$};
\draw[->, thick, red!70!black] (2.4, \rowA+\barH+0.15) -- (2.4, \rowA+\barH+0.55)
  node[above] {\tiny\color{red!70!black} departs};

\draw (3.1,\rowA) -- (3.4,\rowA);
\draw[dashed] (3.1,\rowA) -- (3.1,\rowA+\barH);
\draw[dashed] (3.4,\rowA) -- (3.4,\rowA+\barH);

\fill[blue!20] (3.4,\rowA) rectangle (4.7, \rowA+\barH);
\draw (3.4,\rowA) rectangle (4.7, \rowA+\barH);
\node at (4.05, \rowA+\barH/2) {\tiny $e_{j+2}$};

\draw (4.7,\rowA) -- (4.95,\rowA);
\draw[dashed] (4.7,\rowA) -- (4.7,\rowA+\barH);
\draw[dashed] (4.95,\rowA) -- (4.95,\rowA+\barH);

\fill[blue!20] (4.95,\rowA) rectangle (6.35, \rowA+\barH);
\draw (4.95,\rowA) rectangle (6.35, \rowA+\barH);
\node at (5.65, \rowA+\barH/2) {\tiny $e_{j+3}$};

\draw (6.35,\rowA) -- (6.7,\rowA);
\draw[dashed] (6.35,\rowA) -- (6.35,\rowA+\barH);
\draw[dashed] (6.7,\rowA) -- (6.7,\rowA+\barH);

\draw[->, very thick] (4.5, \rowA - 0.3) -- (4.5, \rowA - 0.9)
  node[midway, right=3pt, align=left] {\small Step 1: departure of $[u_{j+1},v_{j+1})$};

\def\rowB{2.8}

\node[anchor=east] at (-0.3, \rowB+\barH/2) {\small $\widetilde{L}_x^{\delta,2}$};

\fill[gray!30] (0,\rowB) rectangle (1.0, \rowB+\barH);
\draw (0,\rowB) rectangle (1.0, \rowB+\barH);
\node at (0.5, \rowB+\barH/2) {\tiny $e_j$};

\draw (1.0,\rowB) -- (3.4,\rowB);
\draw[dashed] (1.0,\rowB) -- (1.0,\rowB+\barH);
\draw[dashed] (3.4,\rowB) -- (3.4,\rowB+\barH);
\draw[decorate, decoration={brace, mirror, amplitude=5pt}] (1.0, \rowB-0.15) -- (3.4, \rowB-0.15)
  node[midway, below=6pt] {\tiny $L:=g_j+e_{j+1}+g_{j+1}$};
\node[above, align=center] at (2.2, \rowB+\barH+0.08) {\tiny $\alpha_i+x\le L<2\alpha_i$};

\fill[blue!20] (3.4,\rowB) rectangle (4.7, \rowB+\barH);
\draw (3.4,\rowB) rectangle (4.7, \rowB+\barH);
\node at (4.05, \rowB+\barH/2) {\tiny $e_{j+2}$};

\draw (4.7,\rowB) -- (4.95,\rowB);
\draw[dashed] (4.7,\rowB) -- (4.7,\rowB+\barH);
\draw[dashed] (4.95,\rowB) -- (4.95,\rowB+\barH);

\fill[blue!20] (4.95,\rowB) rectangle (6.35, \rowB+\barH);
\draw (4.95,\rowB) rectangle (6.35, \rowB+\barH);
\node at (5.65, \rowB+\barH/2) {\tiny $e_{j+3}$};

\draw (6.35,\rowB) -- (6.7,\rowB);
\draw[dashed] (6.35,\rowB) -- (6.35,\rowB+\barH);
\draw[dashed] (6.7,\rowB) -- (6.7,\rowB+\barH);

\draw[->, very thick] (4.5, \rowB - 0.3) -- (4.5, \rowB - 0.9)
  node[midway, right=3pt, align=left] {\small Step 2: one size-$\alpha_i$ arrival};

\def\rowC{0}

\node[anchor=east] at (-0.3, \rowC+\barH/2) {\small $L_x^{\delta,2}$};

\fill[gray!30] (0,\rowC) rectangle (1.0, \rowC+\barH);
\draw (0,\rowC) rectangle (1.0, \rowC+\barH);
\node at (0.5, \rowC+\barH/2) {\tiny $e_j$};

\fill[green!25] (1.0,\rowC) rectangle (2.3, \rowC+\barH);
\draw (1.0,\rowC) rectangle (2.3, \rowC+\barH);
\node at (1.65, \rowC+\barH/2) {\tiny $\alpha_i$};
\draw[->, thick, green!50!black] (1.65, \rowC+\barH+0.55) -- (1.65, \rowC+\barH+0.15);
\node[above] at (1.65, \rowC+\barH+0.55) {\tiny\color{green!50!black} arrives};

\draw (2.3,\rowC) -- (3.4,\rowC);
\draw[dashed] (2.3,\rowC) -- (2.3,\rowC+\barH);
\draw[dashed] (3.4,\rowC) -- (3.4,\rowC+\barH);
\draw[decorate, decoration={brace, mirror, amplitude=5pt}] (2.3, \rowC-0.15) -- (3.4, \rowC-0.15)
  node[midway, below=6pt] {\tiny residual $h_i(L)\in[x,\alpha_i)$};

\fill[blue!20] (3.4,\rowC) rectangle (4.7, \rowC+\barH);
\draw (3.4,\rowC) rectangle (4.7, \rowC+\barH);
\node at (4.05, \rowC+\barH/2) {\tiny $e_{j+2}$};

\draw (4.7,\rowC) -- (4.95,\rowC);
\draw[dashed] (4.7,\rowC) -- (4.7,\rowC+\barH);
\draw[dashed] (4.95,\rowC) -- (4.95,\rowC+\barH);

\fill[blue!20] (4.95,\rowC) rectangle (6.35, \rowC+\barH);
\draw (4.95,\rowC) rectangle (6.35, \rowC+\barH);
\node at (5.65, \rowC+\barH/2) {\tiny $e_{j+3}$};

\draw (6.35,\rowC) -- (6.7,\rowC);
\draw[dashed] (6.35,\rowC) -- (6.35,\rowC+\barH);
\draw[dashed] (6.7,\rowC) -- (6.7,\rowC+\barH);

\end{tikzpicture}
\caption{\tcgr{Schematic representation of the two updates used in the proof of Proposition~5.3, shown here for $m=3$. Start with a packing configuration in $L_x^{\delta,3}$. In Step~1, the leftmost item $[u_{j+1},v_{j+1})$ departs, and the two adjacent holes merge with that interval into a single hole of length $L=g_j+e_{j+1}+g_{j+1}$. The resulting packing configuration is in $\widetilde{L}_x^{\delta,2}$. In Step~2, we show the case $\alpha_i+x\le L<2\alpha_i$, so that one size-$\alpha_i$ arrival, placed at the left endpoint of this hole by the first-fit rule, leaves a hole of length $h_i(L)\in[x,\alpha_i)$. The resulting packing configuration is then in $L_x^{\delta,2}$. The figure is schematic and not to scale.}}
\label{fig:L-transition-revised}

\end{figure}

We proceed to define $L_x^\delta$ as the union of all $L_x^{\delta,m}$, that is,
\begin{equation}
\begin{aligned}
L_x^{\delta} =&\big\{\, j \ \big|\ x \le \tcgr{g_j} < \alpha_i,\ \exists\ m\in \mathbb{N}_{+},\ \text{s.t.}\ \forall\, 0 < \ell \le m,\ \alpha_i\le \tcgr{e_{j+\ell}} \le A_\delta,\\
& \tcgr{g_{j+\ell}} < \alpha_i,\ H_i\left(\tcgr{e_{j+\ell}+g_{j+\ell}}\right)=h_i\left(\tcgr{e_{j+\ell}+g_{j+\ell}}\right),\ v_{j+m+1}\tcorange{\leq} r(y+\epsilon),\\
& \tcgr{g_j}+\sum_{\ell=1}^{m-1} H_i\left(\tcgr{e_{j+\ell}+g_{j+\ell}}\right) < \alpha_i,\ \tcgr{g_j}+\sum_{\ell=1}^{m} H_i\left(\tcgr{e_{j+\ell}+g_{j+\ell}}\right) \ge \alpha_i \big\}.
\end{aligned}
\end{equation}

We now state Proposition \ref{prop645} below.

\begin{prop}\label{prop645}
 For any $y\in[\beta_{i-1},\beta_i)$ and $0<\epsilon< (\beta_{i}-y)/2 $, under the inductive hypothesis at $y$,
for any \tcblue{$\delta>0$ and }$0 < x < \alpha_i$, 
$$
\lim_{r \to \infty} \frac{\mathbb{E}({|L^\delta_x}|)}{r} = 0.
$$
\end{prop}
\begin{proof}[Proof of Proposition \ref{prop645}]
    We first analyze the asymptotic behavior of $| L^{\delta,1}_x|/r$ \tcp{ through the Lyapunov functional $\mathrm{TF}(r(y+\epsilon))-\big|\widetilde{\mathcal{R}}_i^\delta\big|$.}  Consider any $j\in L_x^{\delta,1}$ with $j\neq\tcp{j_{\min}(y+\epsilon)}$ where the item occupying $[u_j,v_j)$ has size at least $\alpha_i$.  \tcp{By the definition of $L_x^{\delta,1}$, the hole $[v_j,u_{j+1})$ has length in $[x,\alpha_i)$, the item $[u_{j+1},v_{j+1})$ has size in $[\alpha_i,A_\delta]$, the hole $[v_{j+1},u_{j+2})$ has length $<\alpha_i$, and $H_i(\tcgr{e_{j+1}+g_{j+1}})=h_i(\tcgr{e_{j+1}+g_{j+1}})$.}
    
    We now distinguish between two cases: Case 1,
$$
\tcgr{H_i(g_j+e_{j+1}+g_{j+1})>h_i(g_j+e_{j+1}+g_{j+1}),}
$$
and Case 2,
$$
\tcgr{H_i(g_j+e_{j+1}+g_{j+1})=h_i(g_j+e_{j+1}+g_{j+1}).}
$$
In Case 1, when the item located in $[u_{j+1},v_{j+1})$ departs, \tcp{the two holes $[v_j,u_{j+1})$ and $[v_{j+1},u_{j+2})$ merge into the single hole $[v_j,u_{j+2})$ of length $\tcgr{g_j+e_{j+1}+g_{j+1}}\in[\alpha_i,A_\delta+2\alpha_i)$ and Case 1 gives $H_i(\tcgr{g_j+e_{j+1}+g_{j+1}})>h_i(\tcgr{g_j+e_{j+1}+g_{j+1}})$. This is exactly the condition for $j\in\widetilde{\mathcal{R}}^{\delta}_i$} and thus the index $j$ enters $\widetilde{\mathcal{R}}^{\delta}_i$ and  $|\widetilde{\mathcal{R}}^{\delta}_i|$ increases by $1$. \tcp{As stated in Proposition \ref{prop632}, $\mathrm{TF}(r(y+\epsilon))$ remains non-increasing during such a departure event.} Hence
$ 
\mathrm{TF}(r(y+\epsilon)) - |\widetilde{\mathcal{R}}^{\delta}_i| 
$ 
decreases by at least $1$. In Case 2, \tcp{by the definition of $L_x^{\delta,1}$, we have $H_i(\tcgr{e_{j+1}+g_{j+1}})=h_i(\tcgr{e_{j+1}+g_{j+1}})$ and $\tcgr{g_j}+H_i(\tcgr{e_{j+1}+g_{j+1}})\geq \alpha_i$.} Therefore, when the item in $[u_{j+1},v_{j+1})$ departs, we have
$$
\tcgr{h_i(g_j+e_{j+1}+g_{j+1})=g_j+h_i(e_{j+1}+g_{j+1})-\alpha_i,}
$$
\tcgr{since $0\leq g_j<\alpha_i$, $0\leq h_i(e_{j+1}+g_{j+1})<\alpha_i$, and $g_j+h_i(e_{j+1}+g_{j+1})\geq \alpha_i$. Hence}
$$
\begin{aligned}
\Delta\mathrm{TF}
&= H_i(\tcgr{g_j+e_{j+1}+g_{j+1}})\;-\;\Big(H_i(\tcgr{e_{j+1}+g_{j+1}}) + \tcgr{g_j}\Big)\\
&= h_i(\tcgr{g_j+e_{j+1}+g_{j+1}})\;-\;\Big(h_i(\tcgr{e_{j+1}+g_{j+1}}) + \tcgr{g_j}\Big)\\
& = -\,\alpha_i.
\end{aligned}
$$
The equation for $\Delta\mathrm{TF}$ indicates that, for this case, $\mathrm{TF}(r(y+\epsilon))- |\widetilde{\mathcal{R}}^{\delta}_i|$ decreases by $\alpha_i$. We now combine the two cases. For any $j \in L_x^{\delta,1}$ with $j\neq \tcp{j_{\min}(y+\epsilon)}$ where the item occupying $[u_{j},v_{j})$ has size at least $\alpha_i$,  when the item located in $[u_{j+1},v_{j+1})$ departs, the quantity $\mathrm{TF}(r(y+\epsilon))- |\widetilde{\mathcal{R}}^{\delta}_i|$ decreases by at least $\min \left(\alpha_i, 1\right)$ (either $|\widetilde{\mathcal{R}}^\delta_i|$ increases by $1$ or $\mathrm{TF}(r(y+\epsilon))$ drops by $\alpha_i$). Hence, the rate of decrease of $\mathrm{TF}(r(y+\epsilon)) - |\widetilde{\mathcal{R}}^{\delta}_i|$ is at least 
   $$\min \left(\alpha_i, 1\right)\cdot|L_x^{\delta,1}|+\widetilde{o}(r).$$ As we derived in the proof of the last proposition, the rate of increase of $\mathrm{TF}(r(y+\epsilon)) - |\widetilde{\mathcal{R}}^{\delta}_i|$ can be upper bounded by $\Lambda^{\mathrm{TF}}+r\cdot I(|\widetilde{\mathcal{R}}^{\delta}_i| > 0)+2|\widetilde{\mathcal{R}}^{\delta}_i|$. Therefore, 
    \begin{equation*}
    \begin{aligned}
        \mathcal{A}\left(\mathrm{TF}(r(y+\epsilon))-|\widetilde{\mathcal{R}}^{\delta}_i|\right)\leq & -\min \left(\alpha_i, 1\right)\cdot|L_x^{\delta,1}|+\Lambda^{\mathrm{TF}}+ r\cdot I(|\widetilde{\mathcal{R}}^{\delta}_i| > 0)+2|\widetilde{\mathcal{R}}^{\delta}_i|+\widetilde{o}(r) ,
    \end{aligned}
    \end{equation*}
    which implies that
    \begin{equation}\label{eq:l1}
        \lim_{r \to \infty} \frac{\mathbb{E}({|L^{\delta,1}_x}|)}{r} = 0.
    \end{equation}
\tcp{We next derive recursive estimates for $|L_x^{\delta,n}|$ and $|\widetilde{L}_x^{\delta,n}|$.}  \tcgr{For convenience, for $j\in L_x^{\delta,n}$ or $j\in \widetilde{L}_x^{\delta,n}$, we shall refer to the items with indices $j+1,\ldots,j+n$ as the ``internal" items of the corresponding packing configuration. We shall also refer to the item with index $j+1$ as the ``leftmost item" and to the items with indices $j$ and $j+n+1$ (when present) as ``adjacent items."}

For any $n \in \mathbb{N}_{+}$, we analyze the dynamics of $| \widetilde{L}_x^{\delta,n}|$. \tcp{To lower-bound the positive drift of $|\widetilde{L}_x^{\delta,n}|$, it suffices to consider one scenario that increases $|\widetilde{L}_x^{\delta,n}|$, namely, the departure of the leftmost item in a packing configuration counted by $L_x^{\delta,n+1}$. Specifically, for \tcor{every $j\in L_x^{\delta, n+1}$}}, when the item $[u_{j+1}, v_{j+1})$ departs, $j$ becomes a new element of $\widetilde{L}_x^{\delta, n}$, increasing $|\widetilde{L}_x^{\delta, n}|$ by $1$. \tcp{Each packing configuration in $L_x^{\delta,n+1}$ contributes at most one such departure event, and the corresponding leftmost item $[u_{j+1},v_{j+1})$ departs at rate $1$.} Hence the rate of \tcp{the increase of $|\widetilde{L}_x^{\delta,n}|$} is \tcp{at least} $\left|L_x^{\delta, n+1}\right|$. 
\tcgr{Here there is no over-count in the rate term: if $j_1\neq j_2$ are indices in $L_x^{\delta,n+1}$, then the corresponding leftmost items $[u_{j_1+1},v_{j_1+1})$ and $[u_{j_2+1},v_{j_2+1})$ are different  items.} 

There are two mechanisms by which $| \widetilde{L}_x^{\delta,n}|$ can decrease. The first occurs when $|\widetilde{L}_x^{\delta,n}|> 0$ and a new item arrives and occupies one of the holes in a packing configuration from $\widetilde{L}_x^{\delta,n}$. The rate of \tcp{the decrease of $|\widetilde{L}_x^{\delta,n}|$} caused by such arrivals can be upper bounded by $r \cdot I(|\widetilde{L}_x^{\delta,n}| > 0)$. Secondly, \tcp{a decrease in $|\widetilde{L}_x^{\delta,n}|$} may occur when one of the internal or adjacent items in a packing configuration from $\widetilde{L}_x^{\delta,n}$ departs from the system. \tcp{\tcgr{A packing configuration counted by $\widetilde{L}_x^{\delta,n}$ involves at most $n$ internal items together with two adjacent boundary items. If any one of these items departs, then the resulting state may no longer satisfy one of the conditions defining $\widetilde{L}_x^{\delta,n}$.}} Therefore, the rate of \tcp{the decrease of }$|\widetilde{L}_x^{\delta,n}|$ caused by such departures can be upper bounded by $(n+2)\, |\widetilde{L}_x^{\delta,n}|$.  Combining these bounds, the generator applied to $| \widetilde{L}_x^{\delta,n}|$ satisfies
\begin{equation}\label{eqln1}
 \mathcal{A}\big(|\widetilde{L}_x^{\delta,n}|\big) 
   \;\geq\; | L_x^{\delta,n+1}|- r \cdot I\big(| \widetilde{L}_x^{\delta,n}| > 0\big) 
          -(n+2)\,|\widetilde{L}_x^{\delta,n}|.     
\end{equation}
We now examine the dynamics of $|L_x^{\delta,n}|$. It again suffices to consider one scenario that increases $|L_x^{\delta,n}|$: \tcp{When $j_i^* \in \widetilde{L}_x^{\delta,n}$ (which implies $H_i(\tcgr{FH_i})=h_i(\tcgr{FH_i})\ge x$) and $\tcgr{FH_i}<2\alpha_i$, placing an arriving size-$\alpha_i$ item flush with the left boundary of the hole $[v_{j_i^*},u_{j_i^*+1})$ leaves a remaining right hole of length $h_i(\tcgr{FH_i})\in[x,\alpha_i)$, thereby producing a new packing configuration in $L_x^{\delta,n}$ and increasing $\left|L_x^{\delta, n}\right|$ by $1$.} The rate of \tcp{the increase of $|L_x^{\delta,n}|$} caused by this event can be bounded below by 
$$
 p_i r\cdot I\left(j_i^* \in \widetilde{L}_x^{\delta,n},\, \tcgr{FH_i}<2 \alpha_i\right)+\widetilde{o}(r).
$$
On the other hand, \tcp{a decrease in $|L_x^{\delta,n}|$} \tcp{may occur through two mechanisms.} The first occurs when one of the internal or adjacent items of a packing configuration in $L_x^{\delta,n}$ departs from the system, which occurs with rate at most $(n+2)\,|L_x^{\delta,n}|$. The second occurs when an item of size less than $\alpha_i$ arrives in $[r\beta_{i-1},r(y+\epsilon))$ and occupies a hole of a packing configuration in $L_x^{\delta,n}$. \tcorange{An arrival can destroy at most $n+1$ configurations counted by $L_x^{\delta,n}$. Consequently, for fixed $n$, the contribution of arrivals of types $\ell<i$ to the rate of decrease is at most}
$$
\tcorange{(n+1)r\sum_{\ell<i}p_\ell I(G_\ell=0)=\widetilde{o}(r).}
$$
Combining these mechanisms yields
\begin{equation}\label{eqln2}
 \mathcal{A}\left(|L_x^{\delta,n}|\right) 
\geq  p_i r\cdot I\left(j_i^* \in \widetilde{L}_x^{\delta,n},\, \tcgr{FH_i}<2 \alpha_i\right)  
- (n+2)| L_x^{\delta,n}| +\widetilde{o}(r).
\end{equation}
\tcp{\tcgr{Conditional on the state at time $0$, the arrival process on $(0,1/r]$ is independent of the remaining service times of the items present at time $0$. Therefore, the event that no item lying entirely in $[0,r(y+\epsilon))$ departs during $(0,1/r]$, and that exactly one arrival occurs during $(0,1/r]$, namely a size-$\alpha_i$ arrival, has probability at least $c_{i,y,\epsilon}$, where $c_{i, y, \epsilon}$ is defined in \eqref{ciye}.}  \tcgr{Recalling $B^{\min}_{i,y,\epsilon}$ as defined in \eqref{eq:bmin}}, we fix an integer $q$ satisfying $2 \leq q \leq\left\lfloor{A_\delta}/{\alpha_i}\right\rfloor+2$ \tcgr{and suppose that, at time 0, the system is in a state satisfying}}\begin{equation}\label{eq:jistar1}
     \tcp{ B^{\min}_{i,y,\epsilon}\cap
\left\{j_i^* \in \widetilde{L}_x^{\delta, n}, q \alpha_i \leq \tcgr{FH_i}<(q+1) \alpha_i\right\}}.
\end{equation}
\tcbl{On the event that no item completely contained in $[0,r(y+\epsilon))$ departs during $(0,1/r]$ and that \tcorange{exactly one arrival occurs in total during $(0,1/r]$, of type $i$,} the first-fit rule places that arrival flush with the left boundary of the hole $[v_{j_i^*},u_{j_i^*+1})$. Therefore, at time $1/r$, the interval $ [v_{j_i^*}+\alpha_i,u_{j_i^*+1})$ 
is a hole of length
$ 
FH_i-\alpha_i.
$ 
Since \eqref{eq:jistar1} gives
$ 
q\alpha_i \le FH_i < (q+1)\alpha_i,
$ 
we have
$ 
(q-1)\alpha_i \le FH_i-\alpha_i < q\alpha_i.
$ 
Therefore, with probability at least $c_{i,y,\epsilon}$, at time $1/r$ the system is in a state satisfying} 
\begin{equation}\label{49a}
    \tcp{
B^{\min}_{i,y,\epsilon}\cap\left\{j_i^* \in \widetilde{L}_x^{\delta, n},(q-1) \alpha_i \leq FH_i<q \alpha_i\right\},
}
\end{equation}
\tcbl{
where all quantities in the event \eqref{49a}  are in the state at time $1/r$. In steady state, we have
$$
\begin{aligned}
&\mathbb{P}\!\left(B^{\min}_{i,y,\epsilon} \cap
\left\{j_i^*\in \widetilde{L}_x^{\delta,n},\ (q-1)\alpha_i \le FH_i < q\alpha_i\right\}\right)\\
&\qquad \ge c_{i,y,\epsilon}\,
\mathbb{P}\!\left(B^{\min}_{i,y,\epsilon} \cap
\left\{j_i^*\in \widetilde{L}_x^{\delta,n},\ q\alpha_i \le FH_i < (q+1)\alpha_i\right\}\right).
\end{aligned}
$$}
\tcbl{Using \eqref{eq:bmin}, we obtain}
\begin{equation}\label{ldn2}
    \begin{aligned}
   \mathbb{P}\left(j_i^* \in \widetilde{L}_x^{\delta,n}\right)=&\sum_{m=1}^{\left\lfloor\frac{A_{\delta}}{\alpha_i}\right\rfloor+2}\mathbb{P}\left(j_i^* \in \widetilde{L}_x^{\delta, n}, m \alpha_i \leq FH_i<(m+1) \alpha_i\right)\\\leq & \left(\left\lfloor\frac{A_{\delta}}{\alpha_i}\right\rfloor+2\right) \left(\frac{1}{c_{i,y,\epsilon}}\right)^{\left(\left\lfloor\frac{A_{\delta}}{\alpha_i}\right\rfloor+2\right)} \mathbb{P}\left(j_i^* \in \widetilde{L}_x^{\delta,n}, FH_i<2 \alpha_i\right)+\tcp{o(1)}.
\end{aligned}
\end{equation}

Similarly, we define $\widetilde{D}_{x, n, \delta}^L$ as the total number of $i$-items that can potentially (completely) fit into the available empty space within $\left[0, v_{\min \left(\widetilde{L}_x^{\delta, n}\right)}\right)$, where $\min \left(\widetilde{L}_x^{\delta, n}\right)$ denotes the smallest index in $\widetilde{L}_x^{\delta, n}$, and $v_{\min \left(\widetilde{L}_x^{\delta, n}\right)}$ denotes the right endpoint of the item with that index. \tcorange{If $\widetilde{L}_x^{\delta,n}=\varnothing$, set $\widetilde{D}_{x,n,\delta}^L=D_i$.} \tcp{\tcgr{Suppose that no item lying entirely in $[0,r(y+\epsilon))$ departs during $(0,1/r]$, and that \tcorange{exactly one arrival occurs in total during this interval, of type $i$.} Then the first-fit assignment mechanism places that arrival in $[0,v_{\min(\widetilde{L}_x^{\delta,n})})$ and so one unit of empty space that could have accommodated a size-$\alpha_i$ item before $v_{\min(\widetilde{L}_x^{\delta,n})}$ is now used.}} Hence, for any $m\in \mathbb{N}_+$, whenever the event
$$
\left\{|\widetilde{L}_x^{\delta,n}| > 0,\tcblue{1\leq}\widetilde{D}_{x, n, \delta}^L\tcblue{\leq m}\right\}
$$
\tcp{occurs} at time $0$ , then at time $1 / r$, with probability at least $c_{i,y,\epsilon}$, the event
$$
\left\{|\widetilde{L}_x^{\delta,n}| > 0,\widetilde{D}_{x, n, \delta}^L\tcblue{\leq}m-1\right\}
$$
occurs. Thus for any $N\in \mathbb{N}_+$,
\begin{equation}\label{ldn4}
\begin{aligned}
        &\mathbb{P}\left(|\widetilde{L}_x^{\delta,n}| > 0,\ D_i\leq N \right)\leq \mathbb{P}\left(|\widetilde{L}_x^{\delta,n}| > 0,\widetilde{D}_{x, n, \delta}^L\leq N \right)\\
        &\qquad\tcblue{\leq}\ \tcp{(N+1)}\left(\frac{1}{c_{i,y,\epsilon}}\right)^{N}\mathbb{P}\left(|\widetilde{L}_x^{\delta,n}| > 0,\widetilde{D}_{x, n, \delta}^L=0 \right) \tcp{\leq} \tcp{(N+1)}\left(\frac{1}{c_{i,y,\epsilon}}\right)^{N}\mathbb{P}\left(j_i^* \in \widetilde{L}_x^{\delta,n} \right),
\end{aligned}
\end{equation}
\tcp{where in the final inequality in \eqref{ldn4} we used the fact that on the event $\{\,|\widetilde{L}_x^{\delta,n}|>0,\ \widetilde{D}_{x,n,\delta}^L=0\,\}$ there is no hole of length at least $\alpha_i$ in $[0, v_{\min(\widetilde{L}_x^{\delta,n})})$. Thus the smallest index in $\widetilde{L}_x^{\delta,n}$ must \tcgr{be} $j_i^*$, i.e., $j_i^*\in \widetilde{L}_x^{\delta,n}$.}
Combining  \eqref{ldn2} with \eqref{ldn4} yields
\begin{equation}\label{ldn3}
    \mathbb{P}\left(|\widetilde{L}_x^{\delta,n}| > 0,\ D_i\leq N \right) \leq \tcp{(N+1)}\left(\left\lfloor\frac{A_{\delta}}{\alpha_i}\right\rfloor+2\right) \left(\frac{1}{c_{i,y,\epsilon}}\right)^{\left(\left\lfloor\frac{A_{\delta}}{\alpha_i}\right\rfloor+N+2\right)} \mathbb{P}\left(j_i^* \in \widetilde{L}_x^{\delta,n}, \tcgr{FH_i}<2 \alpha_i\right)+\tcp{o(1)}.
\end{equation}
Hence, by  Proposition \ref{npr1}, for any $N$ and $\widetilde{\delta}>0$,  \tcgr{with $B_{\widetilde\delta,i}$ and $z_{y,\epsilon}$ defined in \eqref{eq:B-di} and \eqref{eq:z-ye}, respectively, we have} 
\begin{equation}\tcp{\label{ldn}
    \begin{aligned}
    \limsup_{r\rightarrow\infty} \mathbb{P}\left(|\widetilde{L}_x^{\delta,n}| > 0\right)\leq &\limsup_{r\rightarrow\infty} \mathbb{P}\left(|\widetilde{L}_x^{\delta,n}| > 0,\ D_i\leq N \right)+ \limsup_{r\rightarrow\infty} \mathbb{P}\left( D_i> N \right)\\ \leq & \tcp{(N+1)}\left(\left\lfloor\frac{A_{\delta}}{\alpha_i}\right\rfloor+2\right) \left(\frac{1}{c_{i,y,\epsilon}}\right)^{\left(\left\lfloor\frac{A_{\delta}}{\alpha_i}\right\rfloor+N+2\right)} \limsup_{r\rightarrow\infty}\mathbb{P}\left(j_i^* \in \widetilde{L}_x^{\delta,n}, \tcgr{FH_i}<2 \alpha_i\right) \\ &+\frac{p_i+z_{y,\epsilon}}{p_i-z_{y,\epsilon}}\frac{B_{\widetilde\delta,i}+1}{N}+\frac{6\widetilde\delta}{\alpha_1(p_i-z_{y,\epsilon})}.
\end{aligned}
}\end{equation}
We also note that $| \widetilde{L}_x^{\delta,n}|\leq D_i$. Combining \eqref{ldn} with inequalities \eqref{eqln1} and \eqref{eqln2}, we obtain that for any $N$ and $\widetilde{\delta}>0$,
\begin{equation}\tcp{\label{eq36}
\begin{aligned}
\limsup_{r\rightarrow\infty}\frac{\mathbb{E}\left(|L_x^{\delta,n+1}|\right)}{r}  \leq &  \tcp{(N+1)}\left(\left\lfloor\frac{A_{\delta}}{\alpha_i}\right\rfloor+2\right) \left(\frac{1}{c_{i,y,\epsilon}}\right)^{\left(\left\lfloor\frac{A_{\delta}}{\alpha_i}\right\rfloor+N+2\right)}\frac{n+2}{p_i}\limsup_{r\rightarrow\infty}\frac{ \mathbb{E}\left(| L_x^{\delta,n}|\right)}{r}\\
&+  \left(\frac{p_i+z_{y,\epsilon}}{p_i-z_{y,\epsilon}}\frac{B_{\widetilde\delta,i}+1}{N}+\frac{6\widetilde\delta}{\alpha_1(p_i-z_{y,\epsilon})}\right). \\
 \\
\end{aligned}
}\end{equation}
This implies that, iterating over $n$, \tcp{starting from the base case $\lim_{r\to\infty}\mathbb{E}(|L_x^{\delta,1}|)/r=0$ established above and applying \eqref{eq36} recursively, \tcblue{first letting $N\rightarrow\infty$ for fixed $\widetilde{\delta}>0$ and then letting $\widetilde{\delta}\downarrow0$ at each step,} we obtain} 
\begin{equation}\label{eqln}
   \lim_{r\rightarrow\infty} \frac{\mathbb{E}\left(|L_x^{\delta,n+1}|\right)}{r}=0. 
\end{equation}

\tcp{\tcgr{We proceed to use the estimates for $|L_x^{\delta,n}|$ to obtain an upper bound on}
$|L_x^\delta|$. We fix an item index $j$ such that $[u_j,v_j)$ lies completely in
$[r\beta_{i-1},r(y+\epsilon))$ and consider all pairs $(\widetilde j,\widetilde m)$
with $\widetilde j\in L_x^{\delta,\widetilde m}$ and
$\widetilde j<j\le \widetilde j+\widetilde m$. These are the pairs corresponding to
packing configurations in $\bigcup_{m\ge1}L_x^{\delta,m}$ whose \tcgr{internal} items 
contain $[u_j,v_j)$.} \tcp{For any such pair, we write
$j=\widetilde j+k$ with $1\le k\le \widetilde m$. Since $k-1<\widetilde m$,
the conditions defining $L_x^{\delta,\widetilde m}$ \tcgr{in \eqref{eq:L-sets}} give
}
$$\tcp{
\tcgr{g_{\widetilde j}}
+\sum_{\ell=1}^{k-1}H_i(\tcgr{e_{\widetilde j+\ell}+g_{\widetilde j+\ell}})
<\alpha_i .
}$$
\tcp{Moreover, for each $1\le \ell\le k-1$, the \tcgr{conditions defining $L_x^{\delta,\widetilde m}$ \tcgr{in \eqref{eq:L-sets}}} imply
$\tcgr{g_{\widetilde j+\ell}}<\alpha_i$ and
$H_i(\tcgr{e_{\widetilde j+\ell}+g_{\widetilde j+\ell}})
=h_i(\tcgr{e_{\widetilde j+\ell}+g_{\widetilde j+\ell}})$.
Hence, by Corollary~\ref{hyx2}
}
$$\tcp{
H_i(\tcgr{e_{\widetilde j+\ell}+g_{\widetilde j+\ell}})
\ge \tcgr{g_{\widetilde j+\ell}}.
}$$
\tcp{Therefore,}
\begin{equation}\label{eq:Lx-hole-sum}
\tcgr{\sum_{\ell=0}^{k-1} g_{\widetilde j+\ell}<\alpha_i.}
\end{equation}
\tcp{\tcgr{Among the indices $\widetilde j,\widetilde j+1,\ldots,j-1$, there can be at most $\lfloor \alpha_i/x\rfloor+1$ indices $s$ satisfying $g_s\ge x$. Otherwise, the total length of those holes would be at least $\bigl(\lfloor \alpha_i/x\rfloor+2\bigr)x\ge \alpha_i$, contradicting  \eqref{eq:Lx-hole-sum}. Since every index $s\in L_x^\delta$ satisfies $g_s\ge x$, it follows that among $\widetilde j,\widetilde j+1,\ldots,j-1$ there are at most $\lfloor \alpha_i/x\rfloor+1$ indices belonging to $L_x^\delta$.}}

\tcp{\tcgr{Now let $\widetilde j_*$ be the smallest index among all pairs $\left(\widetilde j,\widetilde m\right)$ satisfying $\widetilde j\in L_x^{\delta,\widetilde m}$ and $\widetilde j<j\le \widetilde j+\widetilde m$. \tcbl{Any index $\widetilde{j}$ for which there exists $\widetilde{m}\in\mathbb{N}_+$ such that
$\widetilde{j}\in L_x^{\delta,\widetilde{m}}$ and
$\widetilde{j}<j\le \widetilde{j}+\widetilde{m}$
lies in $\{\widetilde{j}_*,\widetilde{j}_*+1,\ldots,j-1\}$,
and therefore there are at most $\lfloor \alpha_i/x \rfloor + 1$
possible values of $\widetilde{j}$.}  \tcgr{For each fixed index $\widetilde j$, there is at most one value of $\widetilde m$ such that $\widetilde j\in L_x^{\delta,\widetilde m}$ and $\widetilde j<j\le \widetilde j+\widetilde m$. This is because the quantities
$$
g_{\widetilde j}+\sum_{\ell=1}^{m}H_i(e_{\widetilde j+\ell}+g_{\widetilde j+\ell})
$$
are non-decreasing in $m$.} Therefore, there are at most $\left\lfloor{\alpha_i}/{x}\right\rfloor+1$ pairs $\left(\widetilde j,\widetilde m\right)$ with $\widetilde j\in L_x^{\delta,\widetilde m}$ and $\widetilde j<j\le \widetilde j+\widetilde m$. There are  thus at most $\left\lfloor{\alpha_i}/{x}\right\rfloor+1$ packing configurations in $\bigcup_{m\ge1}L_x^{\delta,m}$ for which $\widetilde j<j\le \widetilde j+\widetilde m$.}}  Based on the above analysis, we proceed to compute the total item length counted over all packing configurations in
$
\bigcup_{m \geq n+1} L_x^{\delta, m}
$ as
$$
\sum_{m=n+1}^{\infty} \sum_{j \in L_x^{\delta, m}} \sum_{k=1}^m\tcgr{e_{j+k}}.
$$
For any $n$, the
length of any single item $[u_j,v_j)$ entirely in $[r\beta_{i-1},r(y+\epsilon))$
is counted at most $\lfloor \alpha_i/x\rfloor+1$ times, and hence
$$
\sum_{m=n+1}^{\infty} \sum_{j \in L_x^{\delta, m}} \sum_{k=1}^m
\tcgr{e_{j+k}}
\le
r(y+\epsilon)\left(\left\lfloor\frac{\alpha_i}{x}\right\rfloor+1\right).
$$
Meanwhile,
$$
\sum_{m=n+1}^{\infty} \sum_{j \in L_x^{\delta, m}} \sum_{k=1}^m\tcgr{e_{j+k}} \geq \sum_{m=n+1}^{\infty} \sum_{j \in L_x^{\delta, m}} m \alpha_i \geq n \alpha_i \sum_{m=n+1}^{\infty}\left|L_x^{\delta, m}\right| .
$$
Combining the previous two displays yields
\begin{equation}\label{ad3}
    |L^{\delta}_x|-\sum_{m=1}^n | L_x^{\delta,m}| \leq \left(\left\lfloor\frac{\alpha_i}{x}\right\rfloor+1\right)\frac{r(y+\epsilon)}{n\alpha_i}.
\end{equation}
Now, combining \eqref{ad3} with \eqref{eqln}, we finally obtain
$$
\lim_{r \to \infty} \frac{\mathbb{E}({| L^{{\delta}}_x|})}{r} = 0.
$$   
\end{proof}
By Proposition \ref{prop645}, for any $0<x<\alpha_i$,  $|L_x^{\delta}|/r$ tends to zero in probability as $r\to\infty$. Moreover, each packing configuration \tcor{counted} by $L_x^{\delta}$ contains at most $\lfloor \alpha_i/x\rfloor + 1$ holes of length at least $x$. \tcp{Thus, up to the multiplicative factor $\left\lfloor \alpha_i/x\right\rfloor + 1$, $|L_x^{\delta}|$ controls the number of holes with length in $[x,\alpha_i)$ that arise within these packing configurations. To bound all holes with length in $[x,\alpha_i)$ inside $[r\beta_{i-1},r(y+\epsilon))$,}  we further define the index set of the holes in $[r\beta_{i-1},r(y+\epsilon))$ with a length of at least $x$ and less than $\alpha_i$ as
$$
\mathcal{G}_x:=\left\{j: v_{j+1}\tcorange{\leq} r(y+\epsilon), \ x\leq \tcgr{g_j} <\alpha_i \right\} .
$$
With this notation in hand, we proceed to state Proposition \ref{prop646}.
\begin{prop}\label{prop646}
      For any $y\in[\beta_{i-1},\beta_i)$ and $0<\epsilon< (\beta_{i}-y)/2 $, under the inductive hypothesis at $y$, and for any $0 < x < \alpha_i$,
$$
\lim _{r \rightarrow \infty} \frac{\mathbb{E}\left(| \mathcal{G}_x|\right)}{r}=0.
$$
\end{prop}
\begin{proof}[Proof of Proposition \ref{prop646}]
Fix $0<x<\alpha_i$. We define the index set
$$
\mathcal{M}_i^\delta:=\left\{j: v_{j+1}\tcorange{\leq} r(y+\epsilon),  \text { and } \left(j \in \mathcal{R}_i \text { or } e_j>A_\delta \text { or } e_j<\alpha_i \text { or } g_j \geq \alpha_i\right)\right\} .
$$
By construction, $\mathcal{M}_i^\delta$ (i) contains $\mathcal{R}_i$; (ii) contains all indices of items with length less than $\alpha_i$ or greater than $A_\delta$ that lie completely in $\left[r \beta_{i-1}, r(y+\epsilon)\right)$, whose total number is upper bounded by $$\sum_{\ell: \alpha_{\ell}>A_\delta} F_{\ell}(\infty)+ \sum_{\ell<i}\left(F_{\ell}(r(y+\epsilon))-F_{\ell}\left(r \beta_{i-1}\right)\right);$$ and (iii) contains all indices $j$ with the hole $[v_j,u_{j+1})$ having length at least $\alpha_i$ in $[r\beta_{i-1},r(y+\epsilon))$. The number of holes having length at least $\alpha_i$ in $[r\beta_{i-1},r(y+\epsilon))$ is upper bounded by $D_i$. Hence, by Propositions \ref{prop31}, \ref{npr1} and \ref{prop643},
\begin{equation}\label{m1}
    \left|\mathcal{M}_i^\delta\right| \leq\sum_{\ell: \alpha_{\ell}>A_\delta} F_{\ell}(\infty)+\widetilde{o}(r).
\end{equation}
Then, with $\mathcal{M}_i^\delta$ as defined above, we can rewrite $L_x^\delta$ as
$$
\begin{aligned}
L_x^\delta= & \left\{j \mid x \leq \tcgr{g_j}<\alpha_i, \tcor{\exists \ m \in \mathbb{N}_{+}, \text {s.t. } }\ v_{j+m+1}\tcorange{\leq} r(y+\epsilon), \  \forall 0<\ell \leq m, \right. \\
&  \  j+\ell \notin \mathcal{M}_i^{\delta}, \ \tcgr{g_j}+\sum_{\ell=1}^{m-1} H_i\left(\tcgr{e_{j+\ell}+g_{j+\ell}}\right)<\alpha_i, \left. \tcgr{g_j}+\sum_{\ell=1}^m H_i\left(\tcgr{e_{j+\ell}+g_{j+\ell}}\right) \geq \alpha_i\right\}.
\end{aligned}
$$
For any $j \in \mathcal{G}_x$ with $j \notin L_{x}^\delta$, and  $j \notin \mathcal{M}_{i}^\delta$,  let
$$
m(j):=\min\Big(\{\,m\ge1:\ j+m\in\mathcal{M}_i^\delta\,\}\ \cup\ \{\,\tcp{j_{\max}(y+\epsilon)}-j\,\}\Big).
$$
Recall that $\tcp{j_{\max}(y+\epsilon)}$ is the index of the last item fully contained in $[r\beta_{i-1},r(y+\epsilon))$. \tcgr{By the definition of $m(j)$ as the smallest integer $m\ge 1$ such that either $j+m\in\mathcal M_i^\delta$ or $j+m=\tcp{j_{\max}(y+\epsilon)}$, we have $j+1,\ldots,j+m(j)-1\notin \mathcal M_i^\delta$.} \tcp{For any $j+\ell\notin \mathcal M_i^\delta$, we have $\tcgr{e_{j+\ell}}\in[\alpha_i,A_\delta]$ and $\tcgr{g_{j+\ell}}<\alpha_i$ (by the definition of $\mathcal{M}_i^\delta$). Hence, by additivity of $H_i$ when the right-hole increment is strictly less than $\alpha_i$ (see Corollary~\ref{hyx2}), we have
}$$\tcp{
H_i(\tcgr{e_{j+\ell}+g_{j+\ell}})= H_i(\tcgr{e_{j+\ell}})+\tcgr{g_{j+\ell}}\geq  \tcgr{g_{j+\ell}}.
}$$
\tcgr{We now prove by contradiction that for any $j \in \mathcal{G}_x \backslash\left(L_x^\delta \cup \mathcal{M}_{i}^\delta\right)$, 
$$
g_j+\sum_{\ell=1}^{m(j)-1}H_i(e_{j+\ell}+g_{j+\ell})<\alpha_i.
$$
For the sake of contradiction, let us assume that
$$
g_j+\sum_{\ell=1}^{m(j)-1}H_i(e_{j+\ell}+g_{j+\ell})\ge \alpha_i.
$$
Then there exists an integer $m$ with $1\le m\le m(j)-1$ such that
$$
g_j+\sum_{\ell=1}^{m-1}H_i(e_{j+\ell}+g_{j+\ell})<\alpha_i,
$$
\tcor{and
$$
g_j+\sum_{\ell=1}^{m}H_i(e_{j+\ell}+g_{j+\ell})\ge \alpha_i.
$$}
Hence, by the definition of $L_x^{\delta,m}$ in \eqref{eq:L-sets}, we have
$$
j\in L_x^{\delta,m}\subseteq L_x^\delta.
$$
This contradicts the assumption that $j\notin L_x^\delta$ and therefore
$$
g_j+\sum_{\ell=1}^{m(j)-1}H_i(e_{j+\ell}+g_{j+\ell})<\alpha_i.
$$
}\tcgr{For each $1\le \ell\le m(j)-1$, we have $j+\ell\notin \tr{\mathcal{M}_i^\delta}$, and hence $g_{j+\ell}<\alpha_i$ and}
$$
\tcgr{H_i(e_{j+\ell}+g_{j+\ell})=H_i(e_{j+\ell})+g_{j+\ell}\ge g_{j+\ell},}
$$
\tcgr{where the equality follows from Corollary~\ref{hyx2}. Substituting these lower bounds into the previous display yields}
$$
\tcgr{\sum_{\ell=0}^{m(j)-1}g_{j+\ell}
=
g_j+\sum_{\ell=1}^{m(j)-1}g_{j+\ell}
\le
g_j+\sum_{\ell=1}^{m(j)-1}H_i(e_{j+\ell}+g_{j+\ell})
<
\alpha_i.}
$$
\tcp{Thus, for $0\le \ell\le m(j)-1$,   the total hole length is strictly less than $\alpha_i$ for the holes $[v_{j+\ell},u_{j+\ell+1})$. In particular, for $0\le \ell\le m(j)-1$}, at most $\lfloor\alpha_i/x\rfloor+1$ of these holes \tcgr{$[v_{j+\ell},u_{j+\ell+1})$  can have length at least $x$. Consequently, for any $\widetilde{j} \in \mathcal{M}_{i}^\delta \cup \{\tcp{j_{\max}(y+\epsilon)}\}$, there are at most $\left\lfloor\alpha_i / x\right\rfloor+1$ indices $j \in \mathcal{G}_x \backslash\left(L_x^\delta \cup \mathcal{M}_{i}^\delta\right)$ such that $j+m(j)=\widetilde{j}$.}
Therefore,
$$
\left|\mathcal{G}_x \backslash\left(L_x^\delta \cup \mathcal{M}_{i}^\delta\right)\right| \leq\left(\left\lfloor\frac{\alpha_i}{x}\right\rfloor+1\right)\left(\left|\mathcal{M}_{i}^\delta\right|+1\right) .
$$
\tcp{Combining this with the inequality $$|\mathcal{G}_x|\le |L_x^\delta|+|\mathcal{M}_i^\delta|+|\mathcal{G}_x\backslash(L_x^\delta\cup \mathcal{M}_i^\delta)|$$ and Proposition~\ref{prop645}, we may absorb the term $|L_x^\delta|$ into the $\widetilde{o}(r)$ remainder below.} By \eqref{m1}, we obtain 
$$
\left|\mathcal{G}_x\right| \leq\left(\left\lfloor\frac{\alpha_i}{x}\right\rfloor+2\right)\sum_{\ell: \alpha_{\ell}>A_\delta} F_{\ell}(\infty)+\widetilde{o}(r) .
$$

Therefore, for any $\delta>0$,
$$
\begin{aligned}
    \limsup_{r \rightarrow \infty} \frac{\mathbb{E}\left(|\mathcal{G}_x|\right)}{r}\tcor{\leq}&\left(\left\lfloor\frac{\alpha_i}{x}\right\rfloor+2\right) \limsup_{r \rightarrow \infty}\frac{1}{r}\mathbb{E}\left(\sum_{\ell: \alpha_\ell> A_{\delta}}F_{\ell}\left(\infty ; \infty\right)\right) \leq \left(\left\lfloor\frac{\alpha_i}{x}\right\rfloor+2\right) \frac{\delta}{\alpha_1}.
\end{aligned}
$$
Hence, it follows that
$$
\lim _{r \rightarrow \infty} \frac{\mathbb{E}\left(|\mathcal{G}_x|\right)}{r}=0.
$$

\end{proof}
\tcbl{By Proposition~\ref{prop646}, the number of holes with length in $[x,\alpha_i)$ is $\widetilde{o}(r)$. By Proposition~\ref{npr1}, the number of holes with length at least $\alpha_i$ is also $\widetilde{o}(r)$. Therefore, for every fixed $x>0$,  we have the number of holes with length at least $x$ in $[r\beta_{i-1},r(y+\epsilon))$ is $\widetilde{o}(r)$.} {\tcgr{Note that if a type-$k$ item satisfies $H_i(\alpha_k)=h_i(\alpha_k)>0$ and both adjacent holes have lengths less than $(\alpha_i-h_i(\alpha_k))/2$, then after that item departs, successive type-$i$ arrivals leave a hole whose length lies in $[h_i(\alpha_k),\alpha_i)$.}} This observation is the key input in the proof of Corollary \ref{prop647}.
\begin{corollary}\label{prop647}
 For any $y\in[\beta_{i-1},\beta_i)$ and $0<\epsilon< (\beta_{i}-y)/2 $, under the inductive hypothesis at $y$, for any $k>i$ satisfying
$$
H_i\left(\alpha_k\right)=h_i(\alpha_k)>0,
$$
it follows that
$$
\lim _{r \rightarrow \infty} \frac{1}{r} \mathbb{E}\left(F_k(r(y+\epsilon) ; \infty)-F_k\left(r \beta_{i-1} ; \infty\right)\right)=0.
$$
\end{corollary}
\begin{proof}[Proof of Corollary \ref{prop647}]
Fix $k>i$ with $H_i(\alpha_k)=h_i(\alpha_k)>0$. For any $\delta>0$ with $\alpha_k<A_\delta$ and $0<x<\alpha_i$, we define
$$
\widetilde{\mathcal{G}}^{\delta}_{x} := \{j : v_{j+1}\tcorange{\leq} r(y+\epsilon), \alpha_i \leq \tcgr{g_j} < A_{\delta}+2\alpha_i, H_i(\tcgr{g_j})=h_i(\tcgr{g_j})\geq x \}.
$$
\tcp{The index set $\widetilde{\mathcal{G}}^{\delta}_{x}$ collects holes in $[r\beta_{i-1},r(y+\epsilon))$ whose lengths lie in $[\alpha_i,\,A_\delta+2\alpha_i)$ and satisfy \tcgr{$H_i(g_j)=h_i(g_j)\ge x$}.}

We proceed to analyze the dynamics of $|\widetilde{\mathcal{G}}_{h_i\left(\alpha_k\right)}^\delta|$. It suffices to consider one mechanism by which $|\widetilde{\mathcal{G}}_{h_i\left(\alpha_k\right)}^\delta|$ can increase. Consider a size-$\alpha_k$ item $\left[u_j, v_j\right)$ lying completely in $\left[r \beta_{i-1}, r(y+\epsilon)\right)$ with $j \neq \tcp{j_{\min}(y+\epsilon)}$ and \tcp{$\tcor{j \neq} j_{\max}(y+\epsilon)$}, and let the two holes adjacent to this item have length less than $\left(\alpha_i-h_i\left(\alpha_k\right)\right) / 2$. We have
$$
h_i\left(\alpha_k\right) \leq h_i\left(\tcgr{e_j}\right)+\tcgr{g_j}+\tcgr{g_{j-1}}<\alpha_i .
$$
\tcp{Here $\tcgr{e_j}=\alpha_k$, and the upper bound follows from the assumption that both adjacent holes have length less than $(\alpha_i-h_i(\alpha_k))/2$.} Upon the departure of this size-$\alpha_k$ item, by Corollary \ref{hyx2}, the newly formed hole $\left[v_{j-1}, u_{j+1}\right)$ satisfies
$$
\begin{aligned}
H_i\left(\tcgr{g_{j-1}+e_{j}+g_{j}}\right) & =H_i\left(\tcgr{e_j}\right)+\tcgr{g_j}+\tcgr{g_{j-1}} \\
& =h_i\left(\tcgr{e_j}\right)+\tcgr{g_j}+\tcgr{g_{j-1}}\\&=h_i\left(\tcgr{g_{j-1}+e_{j}+g_{j}}\right).
\end{aligned}
$$
Thus, \tcp{after the departure of this size-$\alpha_k$ item, the index $j-1$ satisfies} the defining conditions of $\widetilde{\mathcal{G}}_{h_i\left(\alpha_k\right)}^\delta$ and therefore $|\widetilde{\mathcal{G}}_{h_i\left(\alpha_k\right)}^\delta|$ increases by $1$. \tcp{\tcgr{To lower-bound the number of departures of type-$k$ items that increase $|\widetilde{\mathcal G}_{h_i(\alpha_k)}^\delta|$ by $1$, we begin with the total number of type-$k$ items lying completely in $\left[r\beta_{i-1},\,r(y+\epsilon)\right)$ and then exclude those type-$k$ items that do not satisfy the conditions stated below.} Specifically, we exclude the boundary items \tcp{$j_{\min}(y+\epsilon)$} and \tcp{$j_{\max}(y+\epsilon)$}, as well as any item that fails the adjacent hole length requirement (i.e., any item bounded by at least one hole of length at least $(\alpha_i-h_i(\alpha_k))/2$). Since each such hole is adjacent to at most two items, the number of items failing this \tcp{adjacent hole} requirement is upper bounded by $2|\mathcal{G}_{(\alpha_i-h_i(\alpha_k))/2}|+2D_i+2$, which by Proposition \ref{prop646} and Proposition  \ref{npr1} is $\widetilde{o}(r)$.}
Therefore, the rate of \tcp{the increase of }$|\widetilde{\mathcal{G}}^{\delta}_{h_i(\alpha_k)}|$ caused by departures of type-$k$ items is at least
$$
F_k(r(y+\epsilon))-F_k\left(r \beta_{i-1} \right)+\widetilde{o}(r).
$$
On the other hand, there are two scenarios by which $|\widetilde{\mathcal{G}}^{\delta}_{h_i(\alpha_k)}|$ decreases. Firstly, a new item may arrive and occupy the hole in a packing configuration belonging to $\widetilde{\mathcal{G}}^{\delta}_{h_i(\alpha_k)}$, with rate upper bounded by $
r \cdot I\left(|\widetilde{\mathcal{G}}^{\delta}_{h_i\left(\alpha_k\right)}|>0\right) .
$ Secondly, \tcp{a decrease in }$|\widetilde{\mathcal{G}}^{\delta}_{h_i(\alpha_k)}|$ may occur when one of the adjacent items in the packing configuration of $\widetilde{\mathcal{G}}^{\delta}_{h_i(\alpha_k)}$ departs, with rate bounded above by \tcp{$
2 | \widetilde{\mathcal{G}}^{\delta}_{h_i\left(\alpha_k\right)}|\leq 2D_i=\widetilde{o}(r).
$} This gives the following generator inequality
\begin{equation}\label{gh1}
\mathcal{A}\left(|\widetilde{\mathcal{G}}^{\delta}_{h_i\left(\alpha_k\right)}|\right) \geq F_k(r(y+\epsilon))-F_k\left(r \beta_{i-1}\right)-r \cdot I\left(| \widetilde{\mathcal{G}}^{\delta}_{h_i\left(\alpha_k\right)}|>0\right)+\widetilde{o}(r) .
\end{equation}
Similarly, for $|\mathcal{G}_{h_i(\alpha_k)}|$, it suffices to consider one scenario that increases $|\mathcal{G}_{h_i(\alpha_k)}|$: \tcp{when $j_i^* \in\widetilde{\mathcal{G}}^{\delta}_{h_i\left(\alpha_k\right)}$ and $\tcgr{FH_i}<2\alpha_i$, exactly one arriving size-$\alpha_i$ item fits into the hole $[v_{j_i^*},u_{j_i^*+1})$, leaving a \tcgr{remaining} hole of length $(\tcgr{FH_i})-\alpha_i=h_i(\tcgr{FH_i})\in [h_i(\alpha_k),\alpha_i)$ that constitutes a new element of $\mathcal{G}_{h_i\left(\alpha_k\right)}$.} Therefore, \tcp{the rate of increase of }$|\mathcal{G}_{h_i(\alpha_k)}|$ is at least
$$
p_i r \cdot I\left(j_i^* \in \widetilde{\mathcal{G}}^{\delta}_{h_i\left(\alpha_k\right)},\, FH_i<2 \alpha_i\right)+\widetilde{o}(r).
$$
Note that \tcp{a decrease in }$|\mathcal{G}_{h_i\left(\alpha_k\right)}|$ may occur via two scenarios. Firstly, for any $j \in \mathcal{G}_{h_i\left(\alpha_k\right)}$, one of the two items adjacent to the hole $\left[v_j, u_{j+1}\right)$ can depart. Secondly, an item of size less than $\alpha_i$ can arrive in $\left[r \beta_{i-1}, r(y+\epsilon)\right)$ and occupy the hole $\left[v_j, u_{j+1}\right)$. By Propositions  \ref{prop31} and \ref{prop646}, \tcp{the aggregate rate of decrease of }$|\mathcal{G}_{h_i\left(\alpha_k\right)}|$ is bounded above by $2 |\mathcal{G}_{h_i\left(\alpha_k\right)}|+r\cdot \sum_{j=1}^{i-1} I\left(G_j=0\right)=\widetilde{o}(r)$. Thus the generator inequality is
\begin{equation}\label{gh2}
 \mathcal{A}\left(| \mathcal{G}_{h_i\left(\alpha_k\right)}|\right) \geq  p_ir\cdot I\left(j_i^* \in \widetilde{\mathcal{G}}^{\delta}_{h_i\left(\alpha_k\right)},\, FH_i<2 \alpha_i\right)+\widetilde{o}(r).
\end{equation}

Recall in \eqref{ldn3} we obtained an upper bound on $\mathbb{P}\left(\left|\widetilde{L}_x^{\delta, n}\right|>0, D_i \leq N\right)$ in terms of $\mathbb{P}\left(j_i^* \in \widetilde{L}_x^{\delta, n},\, FH_i<2 \alpha_i\right)$. Repeating verbatim the argument leading to \eqref{ldn3} yields an analogous estimate for $|\widetilde{\mathcal{G}}_{h_i\left(\alpha_k\right)}^\delta|$. Namely, for any $N\in\mathbb{N}_+$
\begin{equation}\label{gh3}
\begin{aligned}
\mathbb{P}\left(|\widetilde{\mathcal{G}}_{h_i\left(\alpha_k\right)}^\delta|>0, D_i \leq N\right) \leq &  \tcp{(N+1)\left(\left\lfloor\frac{A_\delta}{\alpha_i}\right\rfloor+2\right)\left(\frac{1}{c_{i,y,\epsilon}}\right)^{N+\left\lfloor \frac{A_\delta}{\alpha_i}\right\rfloor+2}} \mathbb{P}\left(j_i^* \in \widetilde{\mathcal{G}}_{h_i\left(\alpha_k\right)}^\delta,\, FH_i<2 \alpha_i\right)+\tcp{o(1)} \\
\end{aligned}
\end{equation}
Taking expectations  of \eqref{gh2}, we obtain
$$
\mathbb{P}\left(j_i^* \in \widetilde{\mathcal{G}}_{h_i\left(\alpha_k\right)}^\delta, \, FH_i<2 \alpha_i\right)=o(1).
$$
Thus by \eqref{gh3}, for any $N\in\mathbb{N}_+$,
$$
\mathbb{P}\left(|\widetilde{\mathcal{G}}_{h_i\left(\alpha_k\right)}^\delta|>0, D_i \leq N\right)=o(1).
$$
Then taking expectations of \eqref{gh1} and utilizing the previous display, we obtain
\begin{equation}\label{ad4}
    \begin{aligned}
        \frac{1}{r}\mathbb{E}\left(F_k(r(y+\epsilon) ; \infty)-F_k\left(r \beta_{i-1} ; \infty\right)\right)\leq \mathbb{P}\left(|\widetilde{\mathcal{G}}_{h_i\left(\alpha_k\right)}^\delta|>0\right)+{o}(1) \leq \mathbb{P}(D_i>N)+{o}(1) .
    \end{aligned}
\end{equation}
By Proposition  \ref{npr1}, inequality \eqref{ad4} implies that
$$
\lim _{r \rightarrow \infty}\frac{1}{r} \mathbb{E}\left(F_k(r(y+\epsilon) ; \infty)-F_k\left(r \beta_{i-1} ; \infty\right)\right)=0.
$$
This completes the proof.
\end{proof}
The following proposition says that the hydrodynamic-scaled total length of holes in $[r\beta_{i-1},\, r(y+\epsilon))$ vanishes as $r\rightarrow\infty$.
\begin{prop}\label{prop648}
       For any $y\in[\beta_{i-1},\beta_i)$ and $0<\epsilon< (\beta_{i}-y)/2 $, under the inductive hypothesis at $y$,
      we have
      \begin{equation*}
\lim_{r\rightarrow\infty}\frac{1}{r}\mathbb{E}\left[\sum_{\ell=1}^\infty \alpha_{\ell}\left(F_{\ell}(r(y+\epsilon) ; \infty)-F_{\ell}\left(r \beta_{i-1} ; \infty\right)\right)\right]=\left(y+\epsilon-\beta_{i-1}\right).
      \end{equation*}
      \end{prop}
\begin{proof}[Proof of Proposition \ref{prop648}]
\tcgr{By definition, the total hole length in $[r\beta_{i-1},\,r(y+\epsilon))$ equals}
$$
\tcp{
r\left(y+\epsilon-\beta_{i-1}\right)
-\sum_{\ell=1}^{\infty}\alpha_\ell\left(F_\ell(r(y+\epsilon))-F_\ell\left(r\beta_{i-1}\right)\right)
+\widetilde o(r),
}
$$
\tcgr{where the $\widetilde o(r)$ term accounts for the difference between the occupied length in $[r\beta_{i-1},\,r(y+\epsilon))$ and the sum above;  this difference can come only from items that intersect one of the two endpoints of the interval, namely, an item with
$u_j \le r\beta_{i-1} < v_j$
or an item with
$u_j \le r(y+\epsilon) < v_j$}. \tcblue{For every $\delta>0$, the absolute value of this difference is bounded above by
$$
2A_\delta
+
\sum_{\ell:\ \alpha_\ell>A_\delta}
\alpha_\ell F_\ell(\infty).
$$
Since $\mathbb E\left(F_\ell(\infty;\infty)\right)=p_\ell r$,
$$
\limsup_{r\rightarrow\infty}
\frac{1}{r}
\mathbb E\left[
2A_\delta
+
\sum_{\ell:\ \alpha_\ell>A_\delta}
\alpha_\ell F_\ell(\infty;\infty)
\right]
\leq
\sum_{\ell:\ \alpha_\ell>A_\delta}
\alpha_\ell p_\ell
<\delta.
$$
Letting $\delta\rightarrow0$ shows that the absolute value of this difference is $\widetilde o(r)$.}
    \tcp{We fix $\tcp{x\in (0,\alpha_i)}$ and split this total hole length into holes \tcgr{with length less} than $x$ and holes of length at least $x$. Since any two holes \tcgr{with length less} than $x$ must be separated by an item of length at least $\alpha_1$, their total \tcgr{length} is at most
$$
r\left(y+\epsilon-\beta_{i-1}\right)\frac{x}{x+\alpha_1}+x.
$$
\tcgr{We now upper bound the total length of holes with length at least $x$ in $[r\beta_{i-1},r(y+\epsilon))$. Each hole with length in $[x,\alpha_i)$ contributes at most $\alpha_i$, and there are at most $|\mathcal{G}_x|+2$ such holes (including two possible boundary holes) in $[r\beta_{i-1},r(y+\epsilon))$. Hence their total length is at most}
$
\tcgr{\alpha_i\left(\left|\mathcal{G}_x\right|+2\right).}
$
\tcgr{If a hole has length $L\geq \alpha_i$, then}
$
\tcgr{L<2\alpha_i\left\lfloor \frac{L}{\alpha_i}\right\rfloor.}
$
\tcgr{Summing over all holes with length at least $\alpha_i$, and using the definition of $D_i$, their total length is at most}
$
\tcgr{2\alpha_iD_i.}
$
\tcgr{Therefore the total length of holes with length at least $x$ is at most}
$$
\tcgr{\alpha_i\left(\left|\mathcal{G}_x\right|+2\right)+2\alpha_iD_i,}
$$
\tcgr{which is $\widetilde{o}(r)$ by Propositions~\ref{prop646} and~\ref{npr1}.} Therefore, for any $x>0$,}\begin{equation*}
\left[r\left(y+\epsilon-\beta_{i-1}\right)-\sum_{\ell=1}^\infty \alpha_{\ell}\left(F_{\ell}(r(y+\epsilon))-F_{\ell}\left(r \beta_{i-1}  \right)\right)\right]+\widetilde{o}(r) \leq r\left(y+\epsilon-\beta_{i-1}\right) \frac{x}{x+\alpha_1}+x.
\end{equation*}
Thus, by \tcor{Propositions} \ref{npr1} and \ref{prop646}, for any $x>0$,
\begin{equation*}
    \limsup_{r\rightarrow\infty} \mathbb{E}\left[\left(y+\epsilon-\beta_{i-1}\right)-\frac{1}{r}\sum_{\ell=1}^{\infty} \alpha_{\ell}\left(F_{\ell}(r(y+\epsilon) ;\infty)-F_{\ell}\left(r \beta_{i-1};\infty\right)\right)\right]\leq \left(y+\epsilon-\beta_{i-1}\right) \cdot \frac{x}{x+\alpha_1}+x.
\end{equation*}
Letting $x\rightarrow 0$ yields
$$
\lim_{r\rightarrow\infty}\frac{1}{r}\mathbb{E}\left[\sum_{\ell=1}^\infty \alpha_{\ell}\left(F_{\ell}(r(y+\epsilon) ; \infty)-F_{\ell}\left(r \beta_{i-1} ; \infty\right)\right)\right]=\left(y+\epsilon-\beta_{i-1}\right).
$$
\end{proof}      
Combining Corollaries \ref{corr644}--\ref{prop647}, we obtain that, for any $k>i$ with $H_i(\alpha_k)>0$, the hydrodynamic-scaled \tcor{number} of type-$k$ items within $[r\beta_{i-1}, r(y+\epsilon))$ vanishes as $r\to\infty$. Any remaining size-$\alpha_k$ items must satisfy $H_i(\alpha_k)=0$, which in turn implies that $\alpha_k$ is an integer multiple of \tcp{$\alpha_i$. By definition of $h_i$ in \eqref{def:hi} and the fact that $H_i(x) \geq h_i(x)$, $H_i\left(\alpha_k\right)=0$ forces $h_i\left(\alpha_k\right)=0$, hence $\alpha_k \in \alpha_i \mathbb{Z}$}. We now proceed to Proposition \ref{prop649}, which shows that, for any $k>i$, the hydrodynamic-scaled number of type- $k$ items with $H_i\left(\alpha_k\right)=0$ also vanishes as $r \rightarrow \infty$.
\begin{prop}\label{prop649}
    Fix $k>i$ with $H_i\left(\alpha_k\right)=0$.  For any $y\in[\beta_{i-1},\beta_i)$ and $0<\epsilon< (\beta_{i}-y)/2 $, under the inductive hypothesis at $y$,
$$
\lim _{r \rightarrow \infty}  \frac{1}{r} \mathbb{E}\left(F_k(r(y+\epsilon) ; \infty)-F_k\left(r \beta_{i-1} ; \infty\right)\right)=0.
$$
\end{prop}
\begin{proof}[Proof of Proposition \ref{prop649}]
Define
$$
\mathcal{G}^{\left(2 \alpha_i\right)}:=\left\{j: \tcblue{r\beta_{i-1}\leq u_j,}\ \tcgr{g_j} \geq 2 \alpha_i, v_{j+1}\tcorange{\leq} r(y+\epsilon), H_i\left(\tcgr{g_j}\right)=h_i\left(\tcgr{g_j}\right)\right\}, 
$$
as an index set of holes $\left[v_j, u_{j+1}\right)$ of length at least $2 \alpha_i$ whose length satisfies $H_i(\cdot)=h_i(\cdot)$, and which lie entirely in $\left[r \beta_{i-1}, r(y+\epsilon)\right)$.
We then define
$$
Q_i:=\tcblue{\sum_{\substack{j:\ r\beta_{i-1}\leq u_j,\ v_{j+1}\tcorange{\leq} r(y+\epsilon)\\ g_j\geq2\alpha_i}}}\left(\tcgr{g_j}-h_i\left(\tcgr{g_j}\right)\right).
$$
\tcp{Since $x-h_i(x)$ is the largest multiple of $\alpha_i$ not exceeding $x$, the quantity $Q_i/\alpha_i$ is the number of further type-$i$ arrivals that can still be inserted into the holes \tcblue{included in the sum defining $Q_i$} before the remaining length of each such hole becomes less than $\alpha_i$.} In this proof, we shall use
\begin{equation}\label{eq:Qi-Lyap}
Q_i+\sum_{\ell>i: H_i\left(\alpha_{\ell}\right)=0}\alpha_{\ell}
\left(F_{\ell}(r(y+\epsilon))-F_{\ell}(r\beta_{i-1})\right)
\end{equation} as the Lyapunov functional. 
When $j_i^* \in \mathcal{G}^{(2\alpha_i)}$ and a size-$\alpha_i$ item arrives and occupies the hole $\left[v_{j_i^*}, u_{j_i^*+1}\right)$, if $\tcgr{FH_i} \geq 3 \alpha_i$ the change in $Q_i$ is
$$
\Delta Q_i=\left[\left(\tcgr{FH_i}-\alpha_i\right)-h_i\left(\tcgr{FH_i}-\alpha_i\right)\right]-\left[\tcgr{FH_i}-h_i\left(\tcgr{FH_i}\right)\right].
$$
\tcp{Since $h_i(\cdot)$ is the remainder modulo $\alpha_i$, subtracting $\alpha_i$ does not change it, i.e., $h_i(x-\alpha_i)=h_i(x)$ for all $x\ge \alpha_i$. Thus if $\tcgr{FH_i}\ge 3\alpha_i$, \tcgr{the length of the hole after \tcor{the} arrival of size-$\alpha_i$ item is at least $2\alpha_i$}, so the corresponding summand in $Q_i$ decreases by exactly $\alpha_i$, i.e., $\Delta Q_i=-\alpha_i$. If instead $2\alpha_i\le \tcgr{FH_i}<3\alpha_i$, then after inserting one size-$\alpha_i$ item \tcgr{the length of the hole that remains after the arrival} lies in $[\alpha_i,2\alpha_i)$ and \tcblue{does not contribute to $Q_i$}; the summand \tcgr{decreases} from $\tcgr{FH_i}-h_i(\tcgr{FH_i})=2\alpha_i$ to $0$, so $\Delta Q_i=-2\alpha_i$. In either case, $Q_i$ decreases by at least $\alpha_i$.}
\tcblue{On the event
$B^{\min}_{i,y,\epsilon}
\tcorange{\cap\{j_i^*\in\mathcal G^{(2\alpha_i)}\}}$,
\tcorange{where $B^{\min}_{i,y,\epsilon}$ is}
defined in \eqref{eq:bmin}, the first-fit rule places every
arriving size-$\alpha_i$ item in
$\left[v_{j_i^*},u_{j_i^*+1}\right)$.
Therefore, the contribution of size-$\alpha_i$ arrivals
\tcorange{on this event}
to the generator of the Lyapunov functional in
\eqref{eq:Qi-Lyap} is upper bounded by}
\begin{equation}\label{eq:Qi-size-i-arrival}
-\alpha_i p_i r I\left(j_i^*\in\mathcal G^{(2\alpha_i)}\right)
\tcblue{
+\alpha_i p_i r I\left(\widetilde G_i>0\right)
+\alpha_i p_i r I\left(G^{\min}_{i,y,\epsilon}>0\right)
}.
\end{equation}

There are \tcblue{several} possible contributions to \tcp{the rate of increase of }\tcblue{the }\tcp{Lyapunov functional in \eqref{eq:Qi-Lyap}}.
\tcblue{We first consider departures for which the departing item and its two adjacent holes are all completely contained in $[r\beta_{i-1},r(y+\epsilon))$. For such a departure, let $s$ and $t$ denote the lengths of the two adjacent holes.}

Firstly, for $\ell$ with $H_i\left(\alpha_{\ell}\right)=0$, \tcblue{we distinguish the cases $\ell=i$ and $\ell>i$.} \tcgr{For the departures of size-$\alpha_i$ items, if $s<\alpha_i/2$ and $t<\alpha_i/2$, then the merged hole has length $s+\alpha_i+t<2\alpha_i$, and therefore it \tcblue{does not contribute to $Q_i$}. Hence in this case the $Q_i$ term in \eqref{eq:Qi-Lyap} does not increase. Thus a positive increase can occur only if the departing size-$\alpha_i$ item has an adjacent hole of length at least $\alpha_i/2$. For such a departure, the contribution of the new merged hole to $Q_i$ is at most $s+\alpha_i+t$. Summing over all such departures, the contribution of the $\alpha_i$ term to the rate of increase of the Lyapunov functional in \eqref{eq:Qi-Lyap} is at most}
$$
\tcgr{\alpha_i\bigl(2|\mathcal G_{\alpha_i/2}|+2D_i\bigr)=\widetilde{o}(r).}
$$

\tcgr{We next consider departures of size-$\alpha_\ell$ items with $\ell>i$ and $H_i(\alpha_\ell)=0$. Let $\Phi_i$ denote the Lyapunov functional in \eqref{eq:Qi-Lyap}. Since $H_i(\alpha_\ell)=0$, we have $h_i(\alpha_\ell)=0$, and therefore $\alpha_\ell$ is an integer multiple of $\alpha_i$. The second part of \eqref{eq:Qi-Lyap} decreases by exactly $\alpha_\ell$. The contribution of the new merged hole to $Q_i$ is at most $s+\alpha_\ell+t$. Letting $\Delta\Phi_i$ denote the change} \tcblue{in $\Phi_i$, we have}
$$
\tcgr{\Delta\Phi_i\le s+t.}
$$

Secondly, the contribution of the departures of size-$\alpha_\ell$ items with $H_i\left(\alpha_{\ell}\right) \neq 0$ \tcblue{to the rate of increase of the Lyapunov functional in \eqref{eq:Qi-Lyap} is bounded above as follows. The positive change of the Lyapunov functional caused by each such departure is at most}
$$
\tcblue{s+\alpha_\ell+t.}
$$

\tcgr{Summing }\tcblue{the terms involving $s$ and $t$ over the departures considered in the first and second cases}\tcgr{, each hole in $[r\beta_{i-1},r(y+\epsilon))$ is counted at most twice, once as a left adjacent hole and once as a right adjacent hole. Hence the total rate of increase of the functional in \eqref{eq:Qi-Lyap} coming from the $s+t$ terms is at most twice the total hole length in $[r\beta_{i-1},r(y+\epsilon))$, which by Proposition~\ref{prop648} is $\widetilde{o}(r)$.}

\tcblue{For any $\delta>0$, the terms involving $\alpha_\ell$ for departures with $H_i(\alpha_\ell)\neq0$ are bounded above by}
\begin{equation}\label{eq:Qi-departure-bound}
\tcblue{
\begin{aligned}
&\sum_{\substack{\ell:\ H_i(\alpha_\ell)\neq0,\ \alpha_\ell\leq A_\delta}}
\alpha_\ell\left(F_\ell(r(y+\epsilon))-F_\ell(r\beta_{i-1})\right)\\
&\qquad+
\sum_{\ell:\ \alpha_\ell>A_\delta}
\alpha_\ell F_\ell(\infty).
\end{aligned}
}
\end{equation}
\tcblue{For fixed $\delta$, only finitely many item types satisfy $\alpha_\ell\leq A_\delta$. For $\ell<i$, Proposition~\ref{prop31} shows that the number of type-$\ell$ items completely contained in $[r\beta_{i-1},r(y+\epsilon))$ is $\widetilde{o}(r)$. The first sum in \eqref{eq:Qi-departure-bound} may include the contribution of at most one item intersecting $r\beta_{i-1}$, and this contribution is at most $A_\delta$. For $\ell>i$, Corollary~\ref{corr644} applies when $H_i(\alpha_\ell)>h_i(\alpha_\ell)$, and Corollary~\ref{prop647} applies when $H_i(\alpha_\ell)=h_i(\alpha_\ell)>0$. Therefore, the first sum in \eqref{eq:Qi-departure-bound} is $\widetilde{o}(r)$.}

\tcblue{Moreover, the Lyapunov functional in
\eqref{eq:Qi-Lyap} does not increase if the departing item is the
first or last item completely contained in
$[r\beta_{i-1},r(y+\epsilon))$, if it intersects either endpoint,
or if it lies completely outside this interval.}

Thirdly, we consider \tcblue{\tcorange{the remaining} arrivals. Suppose first that a size-$\alpha_\ell$ item is placed completely in a hole $[v_j,u_{j+1})$ whose two adjacent items lie completely in $[r\beta_{i-1},r(y+\epsilon))$. If $H_i(\alpha_\ell)=0$, then $\alpha_\ell$ is an integer multiple of $\alpha_i$. Distinguishing the three cases}
$$
\tcblue{
g_j-\alpha_\ell<\alpha_i,\qquad
\alpha_i\leq g_j-\alpha_\ell<2\alpha_i,
\qquad\text{or}\qquad
g_j-\alpha_\ell\geq2\alpha_i,
}
$$
\tcblue{we obtain}
$$
\tcblue{
\begin{aligned}
&\alpha_\ell I(\ell>i)
+\left(g_j-\alpha_\ell-h_i(g_j-\alpha_\ell)\right)
I(g_j-\alpha_\ell\geq2\alpha_i)\\
&\qquad-
\left(g_j-h_i(g_j)\right)I(g_j\geq2\alpha_i)
\leq0.
\end{aligned}
}
$$
\tcblue{The left-hand side is the change in the Lyapunov functional in \eqref{eq:Qi-Lyap}. If $H_i(\alpha_\ell)\neq0$, the second term in \eqref{eq:Qi-Lyap} does not change. By \eqref{def:hi}, the function}
$$
\tcblue{\left(z-h_i(z)\right)I(z\geq2\alpha_i)}
$$
\tcblue{is nondecreasing in $z$, and hence}
$$
\tcblue{
\left(g_j-\alpha_\ell-h_i(g_j-\alpha_\ell)\right)
I(g_j-\alpha_\ell\geq2\alpha_i)
\leq
\left(g_j-h_i(g_j)\right)I(g_j\geq2\alpha_i).
}
$$
\tcblue{Thus, arrivals into holes whose two adjacent items lie completely in $[r\beta_{i-1},r(y+\epsilon))$ do not increase the Lyapunov functional in \eqref{eq:Qi-Lyap}.}

\tcblue{An arrival intersecting $r\beta_{i-1}$ can increase the Lyapunov functional in \eqref{eq:Qi-Lyap} only when $H_i(\alpha_\ell)=0$. In that case, $\alpha_\ell$ is an integer multiple of $\alpha_i$, and the arrival requires an empty interval of length $\alpha_i$ in $[0,r\beta_{i-1}+\alpha_i)$. Hence, for every $\delta>0$, the contribution of these arrivals is upper bounded by}
\begin{equation}\label{eq:Qi-left-arrival}
\tcblue{
rA_\delta I(\widetilde G_i>0)
+r\sum_{\ell:\ \alpha_\ell>A_\delta}\alpha_\ell p_\ell.
}
\end{equation}
\tcblue{Arrivals lying completely to the left of $r\beta_{i-1}$ or completely to the right of $r(y+\epsilon)$, and arrivals intersecting $r(y+\epsilon)$, do not increase the Lyapunov functional in \eqref{eq:Qi-Lyap}.}

\tcblue{It remains to consider arrivals placed completely in the at most two maximal empty intervals in $[r\beta_{i-1},r(y+\epsilon))$ which are not between two items both completely contained in this interval. Such an arrival increases the second term in \eqref{eq:Qi-Lyap} by at most the item length and may leave at most one new hole counted by $Q_i$. The item length plus the contribution of this new hole is at most the decrease in the total length of these maximal empty intervals caused by the arrival. Hence, the positive contribution of all such arrivals is upper bounded by the rate at which this total length decreases. Since this total length is bounded above by $r(y+\epsilon-\beta_{i-1})$, its expected rate of decrease equals its expected rate of increase in stationarity.}

\tcblue{This total length can increase only when the first or last item completely contained in $[r\beta_{i-1},r(y+\epsilon))$ departs, or when an item intersecting one of the two endpoints departs. Thus, at most four items are involved. The total length of the involved items of size
at most $A_\delta$ is at most $4A_\delta$, and the total length of the involved items of size greater than $A_\delta$ is at most $\sum_{\ell:\alpha_\ell>A_\delta}
\alpha_\ell F_\ell(\infty)$. The total length within $[r\beta_{i-1},r(y+\epsilon))$ of the holes adjacent to these items is $\widetilde{o}(r)$ by Proposition~\ref{prop648}. Therefore, the rate at which this total length increases is upper bounded by}
\begin{equation}\label{eq:Qi-end-empty-intervals}
\tcblue{
4A_\delta
+
\sum_{\ell:\alpha_\ell>A_\delta}
\alpha_\ell F_\ell(\infty)
+\widetilde{o}(r).
}
\end{equation}

Putting \tcblue{\eqref{eq:Qi-size-i-arrival}--\eqref{eq:Qi-end-empty-intervals}} together, \tcblue{taking expectations of the generator of the Lyapunov functional in \eqref{eq:Qi-Lyap} and using \eqref{tildeg} and \eqref{eq:gmin} for the three indicator terms in \eqref{eq:Qi-size-i-arrival} and \eqref{eq:Qi-left-arrival},} we obtain
\begin{equation}\label{eq:Qi-generator}
\tcblue{
\begin{aligned}
0\leq{}&
-\alpha_i p_i r\,\mathbb P\left(j_i^*\in\mathcal G^{(2\alpha_i)}\right)\\
&+\mathbb E\left[
\sum_{\substack{\ell:\ H_i(\alpha_\ell)\neq0,\ \alpha_\ell\leq A_\delta}}
\alpha_\ell\left(F_\ell(r(y+\epsilon);\infty)-F_\ell(r\beta_{i-1};\infty)\right)
\right]\\
&+3r\sum_{\ell:\ \alpha_\ell>A_\delta}\alpha_\ell p_\ell
+4A_\delta+o(r).
\end{aligned}
}
\end{equation}
\tcblue{The expectation in the second line of \eqref{eq:Qi-generator} is $o(r)$ by Proposition~\ref{prop31} and Corollaries~\ref{corr644} and~\ref{prop647}. By the choice of $A_\delta$, $\sum_{\ell:\ \alpha_\ell>A_\delta}\alpha_\ell p_\ell<\delta$, and hence}
$$
\tcblue{
\limsup_{r\rightarrow\infty}
\mathbb P\left(j_i^*\in\mathcal G^{(2\alpha_i)}\right)
\leq
\frac{3\delta}{\alpha_i p_i}.
}
$$
Therefore,
\begin{equation}\label{eq:Qi-first-hole}
\lim_{r\rightarrow\infty}
\mathbb P\left(j_i^*\in\mathcal G^{(2\alpha_i)}\right)=0.
\end{equation}

We now fix $k>i$ with $H_i\left(\alpha_k\right)=0$, so that  $\alpha_k\geq 2\alpha_i$. Consider a size-$\alpha_k$ item $\left[u_j, v_j\right)$ lying completely in $\left[r \beta_{i-1}, r(y+\epsilon)\right)$ with $j \neq \tcp{j_{\min}(y+\epsilon)}$ and \tcp{$\tcor{j \neq} j_{\max}(y+\epsilon)$}, and both holes adjacent to this item have length less than $\alpha_i / 2$.  Upon the departure of this item, the \tcor{newly} generated hole $[v_{j-1},u_{j+1})$ satisfies $\tcgr{g_{j-1}+e_{j}+g_{j}}\geq 2\alpha_i$, and \tcblue{by Corollary~\ref{hyx2}, we have}
$$
\begin{aligned}
H_i(\tcgr{g_{j-1}+e_{j}+g_{j}})=&H_i(\tcgr{e_j})+\tcgr{g_j}+\tcgr{g_{j-1}}\\=&\tcgr{g_j}+\tcgr{g_{j-1}}= h_i(\tcgr{g_{j-1}+e_{j}+g_{j}}),
\end{aligned}
$$
\tcp{where in the above display we use that $\tcgr{e_j}=\alpha_k$, $H_i(\alpha_k)=0$, and that $(\tcgr{g_j})+(\tcgr{g_{j-1}})<\alpha_i$ so that $h_i(\tcgr{g_{j-1}+e_{j}+g_{j}})=\tcgr{g_j}+\tcgr{g_{j-1}}$.} Thus $j-1$ becomes a new element of $\mathcal{G}^{(2\alpha_i)}$ and hence $|\mathcal{G}^{(2\alpha_i)}|$ increases by 1 . \tcblue{The difference $F_k(r(y+\epsilon))-F_k(r\beta_{i-1})$ may count at most one item intersecting $r\beta_{i-1}$. We also exclude any size-$\alpha_k$ item that is the first or last item completely contained in $[r\beta_{i-1},r(y+\epsilon))$. Every hole of length at least $\alpha_i/2$ is adjacent to at most two items, the holes with lengths in $[\alpha_i/2,\alpha_i)$ are counted by $|\mathcal G_{\alpha_i/2}|$, and the number of holes of length at least $\alpha_i$ is at most $D_i$.} Therefore, the rate  of increase of $|\mathcal{G}^{(2\alpha_i)}|$ is at least
$$
\begin{aligned}
    &F_k(r(y+\epsilon))-F_k\left(r \beta_{i-1}\right)-2 |\mathcal{G}_{\alpha_i / 2}|-\tcp{2D_i}\tcblue{-3}
    \\=&F_k(r(y+\epsilon))-F_k\left(r \beta_{i-1}\right)+\widetilde{o}(r).
\end{aligned}
$$
\tcp{A decrease in }\tcgr{$|\mathcal{G}^{(2\alpha_i)}|$} \tcp{may occur via two scenarios: (a) a new item arrives and occupies a hole in a packing configuration belonging to }$\mathcal{G}^{(2\alpha_i)}$ \tcp{with rate upper bounded by }$r\cdot I\left(| \mathcal{G}^{(2\alpha_i)}|>0\right)$; \tcp{and (b) one of the items adjacent to the hole in a packing configuration belonging to }$\mathcal{G}^{(2\alpha_i)}$ \tcp{departs, with rate at most }$\tcp{2|\mathcal{G}^{(2\alpha_i)}|\leq 2D_i =\widetilde{o}(r)}$.
Hence, we obtain the generator inequality
\begin{equation}\label{eqq}
   \mathcal{A}\left(| \mathcal{G}^{(2\alpha_i)}|\right) \geq F_k(r(y+\epsilon))-F_k\left(r \beta_{i-1}\right)-r\cdot I\left(|\mathcal{G}^{(2\alpha_i)}|>0\right)+\widetilde{o}(r). 
\end{equation}

Again, recall that in \eqref{ldn4} we obtained an upper bound of $
\mathbb{P}\left(\left|\widetilde{L}_x^{\delta, n}\right|>0, D_i \leq N\right)
$ in terms of $\mathbb{P}\left(j_i^* \in \widetilde{L}_x^{\delta, n}\right)$. Repeating  verbatim the argument leading to \eqref{ldn4} yields an analogous estimate for $\mathbb{P}\left(|\mathcal{G}^{(2\alpha_i)}|>0, D_i\leq N\right)$. Namely, for any $N\in \mathbb{N}_+$,  
\begin{equation}\label{ad5}
    \mathbb{P}\left(|\mathcal{G}^{(2\alpha_i)}|>0, D_i\leq N\right)\leq \tcp{(N+1)\left(\frac{1}{c_{i,y,\epsilon}}\right)^{N}} \mathbb{P}\left(j^*_i\in\mathcal{G}^{(2\alpha_i)}\right).
\end{equation}
Hence, combining \eqref{ad5} with the expectation of \eqref{eqq} and applying Proposition \ref{npr1}, it follows that, for any $N$ and \tr{$\widetilde\delta>0$},
\begin{equation*}\tcp{
\begin{aligned}
        &\limsup_{r\rightarrow\infty}\frac{1}{r}\mathbb{E}\left(F_k(r(y+\epsilon) ; \infty)-F_k\left(r \beta_{i-1} ; \infty\right)\right)\leq \limsup_{r\rightarrow\infty}\mathbb{P}\left(|\mathcal{G}^{(2\alpha_i)}|>0\right) \\
        \leq& \tcor{(N+1)}\left(\frac{1}{c_{i,y,\epsilon}}\right)^{N}  \limsup_{r\rightarrow\infty}\mathbb{P}\left(j^*_i\in\mathcal{G}^{(2\alpha_i)}\right)+\limsup_{r\rightarrow\infty}\mathbb{P}(D_i>N) \\
        \tcorange{\leq}& \frac{p_i+z_{y,\epsilon}}{p_i-z_{y,\epsilon}}\frac{B_{\widetilde\delta,i}+1}{N}+\frac{6\widetilde\delta}{\alpha_1(p_i-z_{y,\epsilon})}.
\end{aligned}
}\end{equation*}
Thus,
\begin{equation*}
    \lim_{r\rightarrow\infty}\frac{1}{r}\mathbb{E}\left(F_k(r(y+\epsilon) ; \infty)-F_k\left(r \beta_{i-1} ; \infty\right)\right)=0,
\end{equation*}
completing the proof.
\end{proof}

Equipped with Propositions \ref{prop642}--\ref{prop649} and Corollaries \ref{corr644}--\ref{prop647}, the remaining part of the proof of Theorem \ref{thm-ff2} now becomes straightforward. Recall that by Proposition \ref{prop648}, we have
$$
\lim_{r\rightarrow\infty}\frac{1}{r}\mathbb{E}\left[\sum_{\ell=1}^\infty \alpha_{\ell}\left(F_{\ell}(r(y+\epsilon) ; \infty)-F_{\ell}\left(r \beta_{i-1} ; \infty\right)\right)\right]-(y+\epsilon-\beta_{i-1})=0,
$$
and by Corollaries \ref{corr644}--\ref{prop647} together with Propositions \ref{prop31} and \ref{prop649}, we have for any $\ell\neq i$,
\begin{equation*}
    \tcp{
\lim_{r\rightarrow\infty}\frac{1}{r} \mathbb{E}\left(F_\ell(r(y+\epsilon) ; \infty)-F_\ell\left(r \beta_{i-1} ; \infty\right)\right)=0.
}
\end{equation*}
Hence, for any $y\in[\beta_{i-1},\beta_i)$ and $0<\epsilon< (\beta_{i}-y)/2$, and for any $\delta>0$, 

\begin{equation*}
\begin{aligned}
        &\liminf_{r\rightarrow\infty}\frac{\alpha_i}{r}\mathbb{E}\left(F_i(r(y+\epsilon);\infty)-F_i(r\beta_{i-1};\infty)\right)-(y+\epsilon-\beta_{i-1})\\\geq &\liminf_{r\rightarrow\infty}\frac{1}{r}\mathbb{E}\left[\sum_{\ell=1}^\infty \alpha_{\ell}\left(F_{\ell}(r(y+\epsilon) ; \infty)-F_{\ell}\left(r \beta_{i-1} ; \infty\right)\right)\right]-(y+\epsilon-\beta_{i-1})\\&-\limsup_{r\rightarrow\infty}\frac{1}{r}\mathbb{E}\left[\sum_{\ell\neq i, \ \alpha_\ell\leq A_\delta} \alpha_{\ell}\left(F_{\ell}(r(y+\epsilon) ; \infty)-F_{\ell}\left(r \beta_{i-1} ; \infty\right)\right)\right]\\
        &-\limsup_{r\rightarrow\infty}\frac{1}{r}\mathbb{E}\left[\sum_{\alpha_\ell>A_\delta} \alpha_{\ell}\left(F_{\ell}(r(y+\epsilon) ; \infty)-F_{\ell}\left(r \beta_{i-1} ; \infty\right)\right)\right]\\ \geq &-\delta.
\end{aligned}
\end{equation*}Letting $\delta\rightarrow 0$ yields 
$$
\lim_{r\rightarrow\infty}\frac{\alpha_i}{r}\mathbb{E}\left(F_i(r(y+\epsilon);\infty)-F_i(r\beta_{i-1};\infty)\right)=(y+\epsilon-\beta_{i-1}).
$$\tcp{
Note that almost surely,
}$$\tcp{
0\le \frac{\alpha_i}{r}\left(F_i(r(y+\epsilon))-F_i(r\beta_{i-1})\right)\le (y+\epsilon-\beta_{i-1})+ \tr{\frac{\alpha_i}{r}},
}$$\tcp{
since each $i$-item entirely in $[r\beta_{i-1},\,r(y+\epsilon))$ occupies length $\alpha_i$. \tcblue{Therefore,
$$
\frac{\alpha_i}{r}
\left(F_i(r(y+\epsilon);\infty)-F_i(r\beta_{i-1};\infty)\right)
\longrightarrow y+\epsilon-\beta_{i-1}
$$
in $L^1$, and hence in probability.}
\tcor{Also, by the inductive hypothesis at $y$, we have
$$
\frac{1}{r}F_i(r\beta_{i-1};\infty)\xrightarrow{P}0.
$$
}
Consequently, for $y\in[\beta_{i-1},\beta_i)$ and $0<\epsilon< (\beta_{i}-y)/2 $, 
}$$
\frac{1}{r} F_i\left(r\left(y+\epsilon\right) ; \infty\right) \xrightarrow{P}\frac{y+\epsilon-\beta_{i-1}}{\alpha_i}.
$$

This completes the induction step from $y$ to $y+\epsilon$ for sufficiently small $\epsilon$, thereby completing the proof of Theorem~\ref{thm-ff2}.
\medskip   

\textbf{Acknowledgments} We dedicate this paper to our colleague, mentor, and friend, Professor Larry Shepp (1936-2013). Professor Shepp brought the open problem from Coffman, Kadota, and Shepp \cite{CKS85} to our attention. We also wish to thank Professor Ed Coffman and Dr. Quan Zhou for help with simulations. Finally, we are most grateful for the insights of two anonymous referees, whose reports helped greatly improve the quality of this paper.


 \newpage

 \newpage

 \section{Appendix}

 \begin{proof}[Proof of Corollary \ref{cor:gmin}]

\tcblue{We define $L^{\min}_{i,y,\epsilon}$ to be the distance from $r\beta_{i-1}$ to the smallest left endpoint among all items whose left endpoints lie in $[r\beta_{i-1},\,r(y+\epsilon))$. If there is no such item, we set
$$
L^{\min}_{i,y,\epsilon}:=r(y+\epsilon-\beta_{i-1}).
$$
By definition,
$$
0\leq L^{\min}_{i,y,\epsilon}\leq r(y+\epsilon-\beta_{i-1}),
$$
so $L^{\min}_{i,y,\epsilon}$ is a bounded function of the process state projection $S_{r(y+\epsilon)}$. If $G^{\min}_{i,y,\epsilon}>0$, the smallest left endpoint in the definition of $L^{\min}_{i,y,\epsilon}$ is $u_{j_{\min}(y+\epsilon)}$.}
We consider the dynamics of $L^{\min}_{i,y,\epsilon}$.
If $G^{\min}_{i,y,\epsilon}>0$, there exists an interval of length $\alpha_i$ completely contained in the available empty space in
$$
[r\beta_{i-1},\,u_{j_{\min}(y+\epsilon)}).
$$
Therefore, whenever $G^{\min}_{i,y,\epsilon}>0$ and $\widetilde{G}_i=0$, a size-$\alpha_i$
arrival decreases $L^{\min}_{i,y,\epsilon}$ by at least $\alpha_i$. Hence, by \eqref{tildeg}, the rate of decrease of $L^{\min}_{i,y,\epsilon}$ is at least
$$
\alpha_i p_i r\cdot I\left(G^{\min}_{i,y,\epsilon}>0\right)+\tcorange{\widetilde o(r)}.
$$
We next upper bound the rate of increase of $L^{\min}_{i,y,\epsilon}$.
\tcblue{If no item has a left endpoint in $[r\beta_{i-1},\,r(y+\epsilon))$, then $L^{\min}_{i,y,\epsilon}$ does not increase. Otherwise, $L^{\min}_{i,y,\epsilon}$ can increase only when the item with the smallest left endpoint in $[r\beta_{i-1},\,r(y+\epsilon))$ departs. In that case, $L^{\min}_{i,y,\epsilon}$ increases by at most the length of the departing item plus the length within $[r\beta_{i-1},\,r(y+\epsilon))$ of the empty interval immediately to its right.
Since an empty interval of length $x$ can contain at most $\lfloor x/\alpha_i\rfloor$ size-$\alpha_i$ items, the length within $[r\beta_{i-1},\,r(y+\epsilon))$ of the empty interval immediately to the right of the departing item is less than $\alpha_i(D_i+1)$. The length of the departing item is bounded above by}
$$
\tcblue{A_\delta+\sum_{\ell:\alpha_\ell>A_\delta}\alpha_\ell F_\ell(\infty;\infty).}
$$
Therefore, the rate of increase of $L^{\min}_{i,y,\epsilon}$ is at most
$$
A_\delta+\alpha_i(D_i+1)+\sum_{\ell:\alpha_\ell>A_\delta}\alpha_\ell F_\ell(\infty;\infty).
$$
Since $\mathbb E\bigl(\mathcal A L^{\min}_{i,y,\epsilon}\bigr)=0$, it follows that
$$
0
\le
A_\delta+\alpha_i\,\mathbb E(D_i+1)
+\mathbb E\!\left(\sum_{\ell:\alpha_\ell>A_\delta}\alpha_\ell F_\ell(\infty;\infty)\right)
-\alpha_i p_i r\cdot \mathbb P\left(G^{\min}_{i,y,\epsilon}>0\right)+o(r).
$$
Dividing by $r$ and employing Proposition~\ref{npr1},
\tcorange{first letting $r\to\infty$ for fixed $\delta>0$
and then letting $\delta\downarrow0$, we obtain}
$$
\limsup_{r\to\infty}\mathbb P\left(G^{\min}_{i,y,\epsilon}>0\right)=0,
$$
proving \eqref{eq:gmin}.
\end{proof}

\end{document}